\documentclass[11pt]{article}

\usepackage[latin1]{inputenc}
\usepackage[all]{xy}
\usepackage{amsfonts}
\usepackage{amsmath}
\usepackage{amssymb}
\usepackage{amsthm}

\newcommand{\cqfd}{\hfill\blbx \\} 
\def\blbx{\hbox{\vrule height 5pt width 5pt depth 0pt}\medskip}

\newtheorem{theo}{Theorem}[section]
\newtheorem{prop}[theo]{Proposition}
\newtheorem{lem}[theo]{Lemma}
\newtheorem{cor}[theo]{Corollary}

\newtheorem{defi}[theo]{Definition}

\numberwithin{equation}{section}

\newcommand{\CC}{\mathbb{C}}

\newcommand{\EE}{\mathbb{E}}

\newcommand{\HH}{\mathbb{H}}

\newcommand{\LL}{\mathbb{L}}

\newcommand{\NN}{\mathbb{N}}
\newcommand{\PP}{\mathbb{P}}
\newcommand{\QQ}{\mathbb{Q}}
\newcommand{\RR}{\mathbb{R}}

\newcommand{\Aa}{ {\cal A }}
\newcommand{\Ba}{ {\cal B }}

\newcommand{\Ea}{ {\cal E }}
\newcommand{\Sa}{ {\cal S }}

\newcommand{\Fa}{ {\cal F }}
\newcommand{\Ga}{ {\cal G }}

\newcommand{\Ia}{ {\cal I }}

\newcommand{\Ma}{ {\cal M }}
\newcommand{\Ta}{ {\cal T}}

\newcommand{\Pa}{ {\cal P }}

\newcommand{\Qb}{ {\bf Q}}
\newcommand{\Eb}{ {\bf E}}

\newcommand{\tcoprod}{ \mbox{\rm\tiny  $\coprod$}}
\newcommand{\scoprod}{ \mbox{\rm\small  $\coprod$}}
\def\blbx{\hbox{\vrule height 5pt width 5pt depth 0pt}\medskip}

\def \PP{\mathbb{P}}
\def \RR{\mathbb{R}}
\def \EE{\mathbb{E}}
\def \QQ{\mathbb{Q}}
\def \CC{\mathbb{C}}
\def \LL{\mathbb{L}}

\begin{document}

\title{Coalescent tree based functional representations
for some  Feynman-Kac particle models }
\author{P. Del Moral, F. Patras\footnote{The second author was supported by the ANR grant AHBE 05-42234},
S. Rubenthaler\footnote{CNRS UMR 6621, 
Universit\'e de Nice, Laboratoire de Math\'ematiques J. Dieudonn\'e,  
Parc Valrose, 
06108 Nice Cedex 2, 
France}}

\maketitle

\begin{abstract}
We design a theoretic tree-based functional
representation of a class of Feynman-Kac particle
distributions, including an extension of the Wick
product formula to interacting particle systems.  
These weak expansions rely on
an original combinatorial, and permutation group analysis of
a special class of forests.  They 
provide refined non asymptotic propagation
of chaos type properties, as well as  sharp $\LL_p$-mean error
bounds, and laws of large numbers for $U$-statistics.
Applications to particle interpretations of the top eigenvalues, and the ground states of Schrödinger semigroups are also discussed.   
\\
{\em Mathematics Subject of Classification :} Primary: 47D08, 60C05, 60K35, 65C35; Secondary: 31B10, 60J80, 65C05, 92D25.
\end{abstract}

\tableofcontents

\section{Introduction and formulation of the main results}
\subsection{Introduction}
The field of Feynman-Kac semigroups, and their interacting
particle models is one of the most active contact points
between probability,  theoretical chemistry, quantum physics, and  engineering sciences, 
including rare event analysis, and advanced signal processing.
 Whereas it is clearly out of the scope
of the present article to review these applications, some 
of them are briefly pointed out in the further development of this introduction.   For a rather thorough discussion, the interested reader is recommended to consult the pair of books~\cite{fk,doucet}, and the references therein. 

The common aspect of all these applications is that
they represent path integration of free Markov particle evolutions, weighted by some potential functions. The models we consider belong to the general class of branching and interacting particle 
systems. Particles explore the state space as a free 
Markov evolution; during their exploration particles with low potential
are killed, while the ones with high potential value duplicate. From a more probabilistic point of view, these 
evolutionary genetic type particle models can be interpreted as
stochastic linearization, or as  
sophisticated interacting acceptation-rejection type simulation techniques.

During the last two decades, the asymptotic analysis of these 
models has been developed in various directions, including
propagation of chaos analysis,
$\LL_p$-mean error estimates, central limit type theorems, and 
large deviation principles. 
The purpose of the present work is to
 develop algebraic tree-based functional
representations of particle block
distributions, stripped of all analytical superstructure, and 
probabilistic irrelevancies. These weak expansions rely on
an original combinatorial, and permutation group analysis on
a special class of classical and colored  forests that parametrize naturally the trajectories of interacting particle systems. In order to modelize and compute the corresponding path integrals, we develop
an original differential and combinatorial calculus, and take advantage of the natural permutation group action on particle configurations to handle the combinatorial complexity of the problem. Wreath products of permutation groups appear in the process, and new 
closed algebraic formulae for labelled forests enumeration are obtained.
We also discuss Hilbert series techniques for counting the forests relevant to our analysis, that is with 
prescribed number of vertices at each level, or with given coalescence degrees. 

Let us mention that forests, colored graphs and their combinatorics have appeared recently in various fields such as in theoretical physics and gaussian matrix integral models~\cite{alice},
renormalization theory in high-energy physics or Runge-Kutta methods, two fields where the structure and complexity of perturbative expansions has required the development of new tools \cite{brouder,conkre98}. To the best of our knowledge, their introduction in the analysis of Feynman-Kac and interacting particle models is new.
The Laurent type integral representations presented in this article 
do not only reveal the
combinatorial essence of interacting particle distributions. They also 
provide refined non asymptotic propagation
of chaos type properties, as well as explicit formulae for $\LL_p$-mean error
bounds, and laws of large numbers for $U$-statistics. 

The article is divided into four main parts, devoted respectively to
the precise description of Feynman-Kac particle models and 
tree-based integral expansions, to the proof of these formulae, to the combinatorial and group-theoretic analysis, and to the extension to path space models with applications to propagation
of chaos type properties .

In section~\ref{secprel}, we describe a general class of Feynman-Kac semigroups associated
with some free Markov evolution, and with a collection of potential functions. We motivate these abstract mathematical models with a brief discussion on a particle absorption model arising in 
the spectral analysis of Schrödinger operators. The mean field particle models are presented in section~\ref{meanfield}. Section~\ref{particleblock} is concerned with the description of the corresponding particle block distributions. The main results of the
article are housed in section~\ref{forestexpansion}. 

Section~\ref{sectroot}
is essentially concerned with the proof of the tree-based 
Laurent expansions presented in section~\ref{forestexpansion}.

Section~\ref{treesfsec} is devoted to the combinatorial analysis of the class of forests 
involved in the tree-based Laurent expansions. The analysis includes some exotic 
planar forests, the branchs of which are allowed to cross, referred to in the following as jungles. 
The first subsection,
section~\ref{termino}, starts with a short review on the vocabulary of trees
and forests. In section~\ref{algb}, and section~\ref{autojung}, we provide an algebraic 
representation of 
forests and jungles in terms of sequences on maps. Sections~\ref{orbitsec},~\ref{closed} and \ref{wreathsec} introduce the algebraic tools relevant to the combinatorial analysis. Permutation group actions on jungles are introduced. From the Feynman-Kac modelling point of view, they correspond to the action of the permutation group acting naturally on particle systems. This action is studied in detail. For example, a wreath product representation of the stabilizers is obtained. Closed inductive formulas for interacting particle models path enumeration are derived in section~\ref{closed}. Finally, in section~\ref{hseries}, we design
an  Hilbert series method for counting forests and the other famillies of relevant combinatorial and geometrical objects. As an application of these results on trees and forests, some 
explicit calculations of the first order terms of the expansions obtained in section~\ref{sectroot} are given in
section~\ref{someforests}, and a generalized Wick formula for forests 
is obtained in section~\ref{wicksetcion}. 

Section~\ref{extensions} is essentially concerned with 
the extension  of
the previous analysis to Feynman-Kac particle models on path 
spaces. The corresponding tree-based expansions involve a new class of forests, with two distinguished colored  vertices. The combinatorial analysis of these objects
is slightly more involved, but essentially follows the same
line of development as before. This study is  summarized 
in section~\ref{tropsect}, and section~\ref{autotrop}. 
The next two sections, section~\ref{fkpath}, and 
section~\ref{unnpath},
are devoted respectively to unnormalized Feynman-Kac semigroups 
on path spaces, and to the corresponding tree-based
Laurent type expansions.
In the final section, section~\ref{finalsec},  we use the differential 
forest expansion machinery developed previously in this article to derive
precise propagations of chaos type properties of interacting particle models, including strong expansions of the particle block distributions
with respect to Zolotarev type seminorms,
sharp $\LL_p$-mean error estimates,  and a law of large numbers
for $U$-statistics associated to interacting particle models.

\subsection{Feynman-Kac semigroups, and some application areas}\label{secprel} 
To avoid any state space restrictions, and any unnecessary topological assumptions,   
this article is concerned with abstract mathematical models
in general measurable state spaces. The main advantage of this 
general set-up comes from the fact that it 
applies without further work to the analysis of 
path space Feynman-Kac models, and their genealogical tree based interpretations.

Before we continue with the body of the introduction, 
we already present a series of definitions, the benefits of which are reaped in the subsequent 
sections. We let $(E_{n},\mathcal{E}_{n})_{n\geq0}$ be a collection of
measurable state spaces. We denote respectively
by $\mathcal{M}(E_n)$,   $\mathcal{P}(E_n)$, and $\mathcal{B}_{b}(E_{n})$, the set
of all finite signed measures on $(E_{n},\mathcal{E}_{n})$,  the convex subset
of all probability measures,
and
 the Banach
space of all bounded and measurable functions $f$ on $E_{n}$, equipped with
the uniform norm $\Vert f\Vert=\sup_{x_{n}\in E_{n}}{|f(x_{n})|}$.
We let $
\mu_n(f_n)=\int_{E_{n}}~\mu_{n}(dx_n)~f_{n}(x_{n})
$, be  the
 Lebesgue integral of a function  $f_n\in\mathcal{B}_{b}(E_{n})$, with respect to
a measure $\mu_n\in \mathcal{M}(E_n)$. We equip $\Ma(E_n)$ with the total variation norm
$\|\mu_n\|_{\rm\tiny tv}=\sup_{f\in\Ba_b(E_n):\Vert f\Vert\leq 1}~|\mu_n(f)|$.  We recall that a bounded integral operator $Q_n$ from $E_{n-1}$
into $E_n$,  such that, for any $f_n\in\mathcal{B}_{b}(E_n)$,
 the functions
$$
Q_n(f_n)~:~x_{n-1}\in E_{n-1}\mapsto Q_n(f_n)(x_{n-1})=
\int_{E_{n-1}}~Q_n(x_{n-1},dx_n)~f_n(x_n)~\in \RR
$$
are $\Ea_{n-1}$-measurable, and bounded, generates a dual operator $\mu_{n-1}\mapsto \mu_{n-1}Q_n$
from $\mathcal{M}(E_{n-1})$ into  $\mathcal{M}(E_{n})$, and defined by $(\mu_{n-1} Q_{n})(f_{n}):=\mu_{n-1} (Q_{n}(f_{n}))$.

For a bounded integral
operator $Q_1$ from $E_0$ into $E_1$, and
an operator $Q_2$ from $E_1$ into $E_2$, we
denote by $Q_1Q_2$ the composition operator from
$E_0$ into $E_2$, defined for any $f_2\in\Ba_b(E_2)$ by
$(Q_1Q_2)(f_2):=Q_1(Q_2(f_2))$. The tensor power $Q^{\otimes q}_n$
represents
the bounded integral operator  on $E^q_n$, defined  for any $F\in\mathcal{B}_{b}(E_n^q)$ by
 $$
 Q_n^{\otimes q}(F)(x^1_n,\ldots,x^q_n)=\int_{E^q_n}~\left[Q_n(x^1_n,dy^1_n)~\ldots~Q_n(x^q_n,dy^q_n)\right]~F(y^1_n,\ldots,y_n^q)
 $$
We consider a distribution $\eta_{0}$ on $E_{0}$,  a collection of Markov
transitions $M_{n}(x_{n-1},dx_{n})$ from $E_{n-1}$ into $E_{n}$,
and a collection of $\mathcal{E}_{n}$-measurable, and bounded potential functions $G_{n}$ on the state spaces
$E_{n}$.
To simplify
the presentation, and avoid unnecessary technical discussion, we shall suppose
that the potential functions are chosen such that
\begin{equation}
0<\inf_{x_{n}\in E_{n}}G_{n}(x_{n})\leq \sup_{x_{n}\in E_{n}}G_{n}(x_{n})<\infty
\label{potfun}%
\end{equation}
We associate to these objects the
Feynman-Kac measures defined for any function
$f_{n}\in\mathcal{B}_{b}(E_{n})$ by the
following formulae
\begin{equation}
\begin{array}
[c]{c}%
\eta_{n}(f_{n})=\gamma_{n}(f_{n})/\gamma_{n}(1)\quad\mbox{\rm with}\quad
\gamma_{n}(f_{n})=\mathbb{E}[f_{n}(X_{n})\prod_{0\leq k<n}G_{k}(X_{k})]
\end{array}
\label{fk}%
\end{equation}
In (\ref{fk}),
$(X_{n})_{n\geq0}$ represents a Markov chain, taking values
in the state spaces $(E_{n})_{n\geq0}$, with initial distribution
$\eta_{0}$ on $E_0$, and elementary transitions $M_{n}$ from
$E_{n-1}$ into $E_n$. The choice of non homogeneous state spaces
$E_n$ is not innocent. In several application areas the underlying Markov 
model is a path-space Markov chain
\begin{equation}\label{pathmodel}
X_n=(X'_0,\ldots,X'_n)\in E_n=(E'_0\times\ldots\times E'_n)
\end{equation}
The elementary prime variables $X'_n$ represent an
elementary Markov chain taking values in some measurable
spaces $(E'_n,\Ea'_n)$. For instance, in macro-molecular analysis, 
the elementary variables $X'_n$ represents the monomers
in a directed chain $X_n$, with polymerization degree $n$.
In this context, the potential functions $G_n$ reflect the 
intermolecular attraction or repulsion forces between the 
monomers, and the corresponding Feynman-Kac model represents 
the distribution of the polymer configurations in a solvent. 
In signal processing, and more particularly in filtering problems,
the chain $X_n$ represent the paths of a given signal, and the
potential functions $G_n$ reflect the likelihood of their terminal values
$X'_n$, with respect to the sequence of observations delivered by the sensors. In this context, the corresponding Feynman-Kac model represents 
the  conditional distribution of the path of the signal, given the sequence of observations. 
To motivate this article, let us also mention an important, and more classical physical interpretation
of these Feynman-Kac flows.
In theoretical particle physics,
these models
represents the distribution of a particle
evolving in an absorbing medium, with obstacles
related to potential functions $G_n$, taking values in $[0,1]$.
In this context, the particle $X^c_n$
evolves according to two separate
mechanisms. Firstly, it moves from a site $x_{n-1}\in E_{n-1}$, to another $x_n\in E_n$
 according to elementary transitions  $M_n(x_{n-1},dx_n)$.
 Then, it is absorbed with a probability $1-G_{n}(x_{n})$, and placed
 in an auxiliary cemetery state $X^c_n=c$; otherwise it remains
 in the same site. If we let $T$ be the random
 absorption time, then it is not difficult to check that
 \begin{equation}\label{absoreq}
 \gamma_n(f_n)=\mathbb{E}[f_{n}(X^c_{n})~1_{T\geq n}]\quad
 \mbox{\rm and}\quad
 \eta_n(f_n)=\mathbb{E}[f_{n}(X^c_n)~|~{T\geq n}]
\end{equation}
In time homogeneous settings ($G_n=G$ and $M_n=M$), we have
$$
 \gamma_n(1)=\PP\left(T\geq n\right)\sim e^{-\lambda~n}
$$
The positive constant $\lambda$ is a measure of the strength and trapping effects of the obstacle. By a lemma of Varadhan's,
 $\lambda$ coincide with
the logarithmic Lyapunov
exponent $\lambda_0$ of the transition operator $Q(x,dy)=G(x)~M(x,dy)$. Whenever it exists, the corresponding
 eigenfunction $h$ of $Q$ represents the ground state of the operator $Q$.
 
These Feynman-Kac models
 can also be regarded as a natural time discretization of
  continuous time particle absorption models associated with  some
Schr\"odinger  operator.
To be more precise, let us suppose that
 $M(x,dy)=M^{\Delta}(x,dy)$ represent the transitions of a continuous time Markov process with infinitesimal
 generator $L$, during a short period of time $\Delta$; in the sense that
$$
L(f)(x)=\lim_{\Delta\rightarrow 0}\frac{1}{\Delta}~\left[M^{\Delta}(f)(x)-f(x)\right]
$$
for sufficiently regular functions $f$ on $E$. Let us also suppose that the potential functions have the form $G(x)=G^{\Delta}(x)=e^{-V(x)\Delta}$, for
 some non negative energy function $V$. In this case, it is readily checked that the
infinitesimal generator of the particle absorption transitions
$$
Q^{\Delta}(x,dy)=e^{-V(x)\Delta}~M^{\Delta}(x,dy)
$$
is given by  the Schrodinger  operator
$$
L(f)(x)-V(x)f(x)=\lim_{\Delta\rightarrow 0}\frac{1}{\Delta}~\left[Q^{\Delta}(f)(x)-f(x)\right]
$$
Under reasonably weak assumptions, the pair $(\lambda,h)$
coincide respectively
with the top eigenvalue, and the ground state
of  the Schrodinger  operator $L^V=L-V$. The estimation
of these spectral quantities has recently received a lot of attention in modern numerical physics, and quantum
chemistry.

A natural, and common key idea in most of the
Monte-Carlo approximation models is to simulate
a population of walkers or particles mimicking
the evolution of a particle in an absorbing medium.
These models can be interpreted as
a natural mean field particle approximation of the
distribution flow $\eta_n$. We end this section with a
precise description of these nonlinear models.
By the Markov property and the
multiplicative structure of (\ref{fk}), it is easily checked that the flow
$(\eta_{n})_{n\geq0}$ satisfies the following equation
\begin{equation}
\eta_{n+1}=\Phi_{n+1}(\eta_{n})\label{eq}
\end{equation}
The transformations $\Phi_{n+1}:\Pa(E_{n})\rightarrow \Pa(E_{n+1})$  are defined for any pair $(\eta_{n},f_{n+1})\in(\Pa(E_{n})\times\mathcal{B}_{b}(E_{n+1}))$ as follows
$$
\Phi_{n+1}(\eta_{n})(f_{n+1})=
\frac{ \eta_{n}( Q_{n+1}(f_{n+1}) )}{\eta_{n}(Q_{n+1}(1))}
\quad\mbox{\rm 
with}\quad
Q_{n+1}(x_{n},dx_{n+1})=G_{n}(x_{n})\times M_{n+1}(x_{n},dx_{n+1})
$$
Notice that the evolution semigroup $(Q_{p,n})_{0\leq p\leq n}$
of the unnormalized distribution
flow is given by
\begin{equation}\label{eqgamma}
\forall 0\leq p< n\qquad \gamma_n=\gamma_pQ_{p,n}\quad \mbox{\rm with}\quad
Q_{p,n}=Q_{p,n-1}Q_{n}
\end{equation}
and the unnormalized Feynman-Kac measures $\gamma_n$ satify
$\gamma_n(f_n)=\gamma_{n-1}(Q_n(f_n))$.
Equivalently, they can be expressed in terms of the flow $(\eta_p)_{0\leq p\leq n}$ with
the formulae
\begin{equation}\label{gamma}
\gamma_{n}(f_{n})=\eta_{n}(f_{n})~\times~\gamma_{n}(1)\quad\mbox{\rm with}\quad \gamma_{n}(1)=\prod_{0\leq p<n}\eta_{p}(G_{p})
\end{equation}
 for any $f_{n}\in\mathcal{B}_{b}(E_{n})$. 
Returning to the physical interpretation, this product formula shows that the
 Lyapunov exponent $\lambda$ is related to the asymptotic
 mean logarithm
 value of the potential $G_n=G$ with respect to the measures
 $\eta_n$
 \begin{equation}\label{lam}
 \frac{1}{n}\sum_{0\leq p\leq n}\log{\eta_p(G)}\simeq -\lambda
 \end{equation}
 Much more is true. In the case where the
 homogeneous transitions $M_n=M$ are reversible
 with respect to some measure $\mu$, the ground state
 represents the limiting fixed point
 measure of the distribution flow $\eta_n$; in the sense that
 for any bounded measurable function $f$, we have that
 \begin{equation}\label{hm}
 \frac{1}{n}\sum_{0\leq p\leq n}{\eta_p(Gf)}/{\eta_p(G)}\simeq
 {\mu(hf)}/{\mu(h)}
 \end{equation}

\subsection{Mean field interacting particle models} \label{meanfield}
A natural mean field
particle model associated with the nonlinear
Feyman-Kac flow (\ref{eq}) is the
$E_{n}^{N}$-valued Markov chain $\xi_{n}^{(N)}=(\xi_{n}^{i,N})_{1\leq i\leq N}$
  with 
elementary transitions defined for any $F_n^N\in \Ba_b(E_{n}^{N})$ by
\begin{equation} 
\EE\left(F_n^N(\xi_{n}^{(N)})~|~\xi_{n-1}^{(N)}
\right)=\Phi_{n}\left(
m\left(\xi_{n-1}^{(N)}\right)\right)^{\otimes N}(F_n^N)
\quad\mbox{\rm with}\quad
m\left(\xi_{n-1}^{(N)}\right)=\frac{1}{N}\sum
_{i=1}^{N}\delta_{\xi_{n-1}^{i,N}}
\label{gent}%
\end{equation}
In other terms, given the
configuration $\xi_{n-1}^{(N)}$ at rank $(n-1)$,
the particle system  $\xi_{n}^{(N)}$ at rank $n$,
 consists of $N$ independent
and identically distributed random variables with common distribution $\Phi_{n}\left(
m\left(\xi_{n-1}^{(N)}\right)\right)$.
The initial
configuration $\xi_{0}^{(N)}$ consists of $N$ independent
and identically distributed random variables with distribution $\eta_{0}$.  
Although the dependency of the model on $N$ is strong, due to the mean-field nature of the model, we abbreviate $(\xi_{n}^{i,N})_{1\leq i\leq N}$ to $(\xi_{n}^{i})_{1\leq i\leq N}$ when $N$ is fixed, excepted when we want to emphasize explicitely the dependency on $N$ of the model, e.g. as in the end of the present section.

Notice that
$$
\Phi_{n}\left(
m\left(\xi_{n-1}^{(N)}\right)\right)(dx_n)=\frac{1}{
\sum_{j=1}^NG_{n-1}(\xi_{n-1}^{j})
}\left(\sum_{i=1}^NG_{n-1}(\xi_{n-1}^{i})~M_n(\xi_{n-1}^{i},dx_n)
\right)$$
so that the particle model evolves as a genetic type
model with proportional selections,
and mutation transitions
dictated by the pair of potential-transition  $(G_{n-1},M_n)$.
In numerical physics, the selection transition is often called reconfiguration of the  population of walkers.
In the context of path space Markov chain model (\ref{pathmodel}),
it is important to notice that the particle interpretation model
is again a path particle model. During the selection stage, 
we select a path particle according to its fitness, and the
mutation stage consists in extending the selected path with
an elementary transition dictated by the Markov transitions
of the elementary prime chain. In this context, the particle model
can be interpreted as a genealogical particle evolution model.

The approximation measures $(\eta^N_n,\gamma^N_n)$ associated with the pair of Feynman-Kac measures $(\eta_n,\gamma_n)$ are defined by the empirical occupation measures
$$
\eta^N_n=m\left(\xi_{n}^{(N)}\right)
$$
and the unnormalized particle distributions $\gamma^N_n$ are defined
for any $f_{n}\in\mathcal{B}_{b}(E_{n})$ by
\begin{equation}\label{gammaN}
\gamma^N_n(f_n)=\eta^N_n(f_n)~\times~\gamma^N_{n}(1)\quad\mbox{\rm with}\quad \gamma_{n}^N(1)=\prod_{0\leq p<n}\eta^N_{p}(G_{p})
\end{equation}
For a rather complete asymptotic analysis of these measures
 we again refer the reader to~\cite{fk}, and references therein.
 In particular, for any $f_{n}\in\mathcal{B}_{b}(E_{n})$, we have the following almost convergence results
 $$
\lim_{N\rightarrow\infty}\eta^N_n(f_n)=\eta_n(f_n)\quad\mbox{\rm and}\quad
\lim_{N\rightarrow\infty}\gamma^N_n(f_n)=\gamma_n(f_n)
$$
 Therefore, the particle interpretations of the pair objects $(\lambda,h)$
 introduced in (\ref{lam}) and (\ref{hm}), are defined by simply with replacing $\eta_n$ by the particle occupation measure $\eta^N_n$.  It is also well known that $\gamma^N_n$ is an
 unbias
 approximation measure of $\gamma_n$, in the sense that
$$
\EE(\gamma^N_n(f_n))=\gamma_n(f_n)
$$
However, it turns out that for any $n\geq 1$
$$
\EE(\eta^N_n(f_n))=\EE(f_n(\xi^{1,N}_n))\not=\eta_n(f_n)
$$
This means that the mean field particle interpretation
model is not an exact sampling algorithm of the distributions $\eta_n$.   In practice, it is clearly important to analyze these quantities, and more generally the distribution
of particle blocks of any size $q\leq N$
$$
\PP_{n,q}^{N}=\mbox{\rm Law}\left(\xi^{1,N}_n,\ldots,\xi^{q,N}_n\right)\in \Pa(E_n^{q})
$$
\subsection{Particle block distributions}\label{particleblock}

For any pair of integers $(l,m)\in(\NN^{\star})^2$,
we set $[l]=\{1,\ldots,l\}$, and $[l]^{[m]}$ the set of mappings $a$ from $[m]$ into $[l]$. By $|a|$, we denote the cardinality of the set
$a([m])$, and for any $1\leq p\leq l$ we set
$$
[l]^{[m]}_p=\{a\in [l]^{[m]}~:~|a|=p\}
$$
The number $m-|a|$ will be refered to as the coalescence number of $a$.
To simplify the presentation, when $l>m$, we also denote by
$\langle m,l\rangle\left(=[l]^{[m]}_{m}\right)$ the set of all 
$$(l)_{m}:=\frac{l!}{(l-m)!}$$
one to one mappings from $[m]$ into $[l]$, and by
$\Ga_l=\langle l,l\rangle$
the symmetric group of all permutations of $[l]$.

From the pure mathematical point of view,
the measures $\PP_{n,q}^{N}$ are better understood when they are connected
with the
$q$-tensor product occupation measures (resp. the restricted $q$-tensor product occupation measures) on $E^q_n$ defined by
$$
(\eta^N_n)^{\otimes q}=\frac{1}{N^q}\sum_{a\in [N]^{[q]}}~
\delta_{(\xi^{a(1)}_n,\ldots,\xi^{a(q)}_n)}\quad\mbox{\rm and}\quad
(\eta^N_n)^{\odot q}=\frac{1}{(N)_q}\sum_{a\in\langle q,N\rangle }~
\delta_{(\xi^{a(1)}_n,\ldots,\xi^{a(q)}_n)}
$$

Notice that these measures are, by construction, symmetry-invariant. That is, for any $F\in \Ba_b(E_n^q)$ 
$$\left(\eta^N_n\right)^{\otimes q}(F)=\left(\eta^N_n\right)^{\otimes q}((F)_{\rm\tiny sym })\quad\mbox{
and}
\quad
\left(\eta^N_n\right)^{\odot q}(F)=\left(\eta^N_n\right)^{\odot q}((F)_{\rm\tiny sym })
$$
with the symmetrization operator
$$
\forall q\geq 1\qquad
F\in \Ba_b(E_n^{q})\mapsto
(F)_{\rm\tiny sym }:=\frac{1}{
p!}
\sum_{\sigma\in\Ga_{q}}
D_{\sigma}F
$$
In the above displayed definition, the operator $D_{b}$ stands for the
Markov transition from $E^{q}_n$ into itself, associated
with a  mapping $b\in [q]^{[q]}$, and
defined by
$$
D_{b}(F)(x^1_n,\ldots,x^{q}_n)=F(x_n^{b(1)},\ldots,x_n^{b(q)})
$$
for any $F\in \mathcal{B}_{b}(E^{q}_n)$, and
$(x^1,\ldots,x^{q})\in E^{q}_n$. 
 In particular, we may assume, without restriction, in our forthcoming computations on $q$-tensor product occupation measures, that $F$ is a symmetric function. That is, for all $s\in \Ga_q$, $D_{s}F=F$. We write from now on, $\Ba_b^{sym}(E_n^q)$ for the set of all symmetric functions in $\Ba_b(E_n^q)$.

 Notice also the following symmetry property, essential in view of all the forthcoming computations. Since, conditional to $\xi_{n-1}^{(N)}$, the $\xi_n^{i}$ are i.i.d., for any $F\in \Ba_b(E_n^q)$ and any $a ,b\in \langle q,N\rangle ^2$, we have
 $$\EE\left[F(\xi^{a(1)}_n,\ldots,\xi^{a(q)}_n)\right]=\EE\left[F(\xi^{b(1)}_n,\ldots,\xi^{b(q)}_n)\right].$$
 For instance, we have that
 \begin{eqnarray}
 \EE\left(\left(\eta^N_n\right)^{\odot q}(F)\right) =
 \EE\left(F\left(\xi^{1}_n,\ldots,\xi^{q}_n\right)
\right)=\PP^{N}_{n,q}(F)
 \end{eqnarray}
As in the case of elementary particle blocks of size one,
the precise analysis of these $q$-tensor product measures
is intimately related to
 the unnormalized pair particle measures  $((\gamma^N_n)^{\odot q},(\gamma^N_n)^{\otimes q})$ defined   for any $F\in \Ba_b(E_n^q)$ by
$$
(\gamma^N_n)^{\odot q}(F)=(\eta^N_n)^{\odot q}(F)\times(\gamma^N_n(1))^q
\quad\mbox{
and}
\quad
(\gamma^N_n)^{\otimes q}(F)=(\eta^N_n)^{\otimes q}(F)\times(\gamma^N_n(1))^q
$$

We consider now the nonnegative measure ${\QQ}^{N}_{n,q}$, indexed by the particle block sizes,
 on the product state spaces
$E_n^q$, and defined   for any $F\in \Ba_b(E_n^q)$ by
$$
{\QQ}^{N}_{n,q}(F):=\EE((\gamma^N_{n})^{\otimes q}(F))
$$
In the further development of this article, we shall prove the
following formula (\ref{gamgam})
\begin{eqnarray*}
\EE((\gamma^N_n)^{\odot q}(F))=\EE((\gamma^N_{n-1})^{\odot q}Q^{\otimes q}_n(F))=
{\QQ}^{N}_{n-1,q}(Q^{\otimes q}_n(F))\label{Gamref}
\end{eqnarray*}
This simple observation indicates that the distribution flow
$({\QQ}^{N}_{n,q})_{n\geq 0}$
has a linear evolution semigroup, for any particle block size $q\leq N$. Much more is true. As another consequence of the results presented in \cite{fk}, the sequence
of distributions ${\QQ}^{N}_{n,q}$ converges, as $N$ tends to infinity, to the distribution of $q$ non absorbed, and independent particles $(X_n^{c,i})_{1\leq i\leq q}$
evolving in the original absorbing medium. That is, we have that
\begin{eqnarray*}
\lim_{N\rightarrow\infty}{\QQ}^{N}_{n,q}(F)&=&\gamma^{\otimes q}_n(F)=
\EE\left(F(X^{c,1}_n,\ldots,X^{c,q}_n)~1_{T^1\geq n}~\ldots~1_{T^q\geq n}\right)
\end{eqnarray*}
for any $F\in \Ba_b(E_n^q)$, and where $T^i$
stands for the random absorption time sequence of the 
chain $(X^{c,i}_k)_{k\geq 0}$, with $1\leq i\leq q$.

This article is mainly concerned with explicit
functional expansions of the
deterministic measures
 $\PP^{N}_{n,q}\in \Pa(E_n^{q})$, and ${\QQ}^{N}_{n,q}\in \Ma(E_n^{q})$, with respect to the precision
 parameter $N$. These (Laurent type) expansions reflect the 
 complete interaction structure of the particle model.
 Roughly speaking, the $k$-th order terms of these
  integral representations represent the 
 $1/N^k$-contributions of mean field particle scenarios with an interaction degree $k$ (see Thm.~\ref{mainthm1}). 
 
To describe more precisely these functional
representations, it is convenient to introduce the abstract definition of the derivative
of a sequence of finite signed measures. 
\begin{defi}
We let $(\Theta^N)_{N\geq 1}\in \Ma(E)^{\NN}$,  be a uniformly bounded sequence of signed measures on a measurable space $(E,\Ea)$, in the sense that $\sup_{N\geq 1}{\|\Theta^N\|_{\rm\tiny tv}}<\infty$. We suppose that
$\Theta^N$ strongly converges to some measure $\Theta\in \Ma(E)$, 
as $N\uparrow\infty$, in the sense that 
 $$
\forall F\in \Ba_b(E)\qquad
\lim_{N\uparrow\infty} \Theta^N(F)= \Theta(F)
$$ 
The discrete derivative of the sequence $(\Theta^N)_{N\geq 1}$ is the sequence of
measures $ (\partial \Theta^N)_{N\geq 1}$ defined by
 $$
 \partial \Theta^N:=N~\left[ \Theta^N-\Theta\right]
 $$
We say that $\Theta^N$ is differentiable, if $ \partial \Theta^N$
is uniformly bounded, and if it 
 strongly converges to some measure $\partial \Theta\in \Ma(E)$, 
as $N\uparrow\infty$.  The discrete derivative of a discrete derivative 
$ \partial \Theta^N$ (of a differentiable sequence) is called the second discrete derivative and it is denoted by 
$$ 
\partial^2 \Theta^N=N\left[\partial\Theta^N-\partial \Theta\right]
$$
The discrete derivative $\partial\Theta^N$ of a differentiable 
sequence  can itself be differentiable.  In this situation, the derivative of the discrete derivative is called the second derivative and it is denoted by $ \partial^2 \Theta=\partial  \left(\partial\Theta \right)$, and so on. 
\end{defi}
At this point, it is convenient to make some comments. Firstly,
the strong form of the weak convergence of measures considered above corresponds
to the s-topology on the set of measures, defined as the coarsest
topology for which all functionals $\mu\in \Ma(E)\mapsto \mu(\varphi)\in\RR$, $\varphi\in\Ba_b(E)$, are continuous. For a more detailed discussion on the s-topology in the context of probability measures, we refer the reader to the article~\cite{jacka}. A sequence
of real numbers can be interpreted as a sequence of measures with a constant value, so that the above definition applies to  any sequence
of real numbers. We also notice that 
the Euler's derivation operators $\partial^{k}$ satisfy the usual linearity
properties. Finally, we notice that a sequence $\Theta^N$ that is
 differentiable up to order $(k+1)$, has the following integral representation
$$
\Theta^N=\sum_{0\leq l\leq k }\frac{1}{N^l}~\partial^l\Theta+\frac{1}{N^{k+1}}~\partial ^{k+1}\Theta^N
$$
with $\sup_{N\geq 1}{\|\partial ^{k+1}\Theta^N\|_{\rm\tiny tv}}<\infty$
, and the convention $\partial^0\Theta=\Theta$, for $l=0$.
\subsection{Forest based expansion formulas}\label{forestexpansion}
 For transparency, we will often use terms
associated with botanical, and genealogical trees. The forthcoming measures expansions will be therefore given in terms of signed measures
 indexed by collections of trees, colored forests and other analogous combinatorial and geometrical objects. More generally, most computations on particle measures will be rephrased in the course of the article in terms of enumerative problems on tree-like objects. In order to avoid lenghtening the present introduction, their definition and study are postponed to the next sections, and we present here our main results in a classical set-theoretic langage.
 
Let $\Aa_{n,q}=_{\mbox{\rm\tiny def.}}([q]^{[q]})^{n+1}$
be the set of $(n+1)$-sequences of mappings ${\bf a}=(a_p)_{0\leq p\leq n}$ from $[q]$ into itself.
$$
[q]\stackrel{ a_0 }{ \longleftarrow }[q]\stackrel{ a_1 }{ \longleftarrow }\cdots\longleftarrow[q]\stackrel{ a_{n-1} }{ \longleftarrow }[q]\stackrel{ a_n }{ \longleftarrow }[q]
$$
Notice that we write in bold the symbols for sequences (of maps, integers...) such as $\bf a$.

We let $\Delta_{n,q}$ be the nonnegative measure valued
 functional on
 $\Aa_{n,q}$ defined by

 \begin{equation}\label{themapDel}
\Delta_{n,q}~:~{\bf a}\in \Aa_{n,q}\mapsto
\Delta^{\bf a}_{n,q}=\left(\eta_0^{\otimes q}D_{a_0}Q_1^{\otimes q}D_{a_1}\ldots Q_n^{\otimes q}D_{a_n}\right)
\in \Ma(E_n^{q})
\end{equation}
The Markov operator $D_b$
can be seen as a coalescent, or a selection type transition.
In this interpretation, the population $(x_n^{b(1)},\ldots,x_n^{b(q)})
$ results from a selection of the individuals with labels in the set
$b([q])$. In addition, arguing as in (\ref{absoreq}), the non negative
integral tensor product operators $Q_n^{\otimes q}$ can be seen
	as the overlapping of an absorption transition from $
	(E_{n-1}\cup\{c\})^q$ into itself, and an exploration transition
	from $
	(E_{n-1}\cup\{c\})^q$ into $
	(E_{n}\cup\{c\})^q$. We then have that
	\begin{equation}\label{coalmarkov}
	\Delta^{\bf a}_{n,q}(F)=\EE\left(F(X^{{\bf a},c,1}_n,\ldots,
	X^{{\bf a},c,q}_n)~1_{T^{{\bf a},1}\geq n}~\ldots~1_{T^{{\bf a},q}\geq n}\right)
	\end{equation}
In the above display, $(X^{{\bf a},c,1}_n,\ldots,
	X^{{\bf a},c,q}_n)$ represents the absorbed Markov chain on $(E_{n}\cup\{c\})^q$ with transitions $Q_n^{\otimes q}D_{a_{n}}$, and initial distribution $\eta_0^{\otimes q}D_{a_{0}}$.

Notice that the measures $\Delta_{n,q}^{\bf a}$ inherit a remarquable invariance property from their set-theoretic definition. Namely, let us introduce the natural left action of the group $\Ga_q^{n+2}$ on $\Aa_{n,q}$ defined for all ${\bf a}\in \Aa_{n,q}$ and all ${\bf s}=(s_0,...,s_{n+1})\in \Ga_q^{n+2}$ by
$${\bf s} ({\bf a}):=( s_0a_0s^{-1}_1,s_1a_1s^{-1}_2,...,s_{n}a_ns^{-1}_{n+1})$$
Then, for any $F\in \mathcal{B}_{b}^{sym}(E^{q}_n)$ we have
\begin{equation}\label{invariance}
\Delta_{n,q}^{\bf b}(F)=\Delta_{n,q}^{\bf s(\bf b)}(F)
\end{equation}
The identity wouldn't hold if $F$ wasn't a symmetric function. However, as already mentioned in Sect.~\ref{particleblock}, this is not a serious restriction as far as the determination of the $q$-tensor occupation measures is concerned.
The action of $\Ga_q^{n+2}$ on $\Aa_{n,q}$ induces a partition of $\Aa_{n,q}$ into orbits or, equivalently, a partition into equivalence classes, where $\bf a$ and $\bf b$ are equivalent, $\bf a\approx\bf b$, if and only if there exists ${\bf s}\in \Ga_q^{n+2}$ such that $\bf a=\bf s(\bf b)$. We write ${\cal F}_{n,q}$ for the set of orbits or equivalence classes. As we shall see later, ${\cal F}_{n,q}$ is nothing but (up to a canonical isomorphism) a remarquable subset of the classical set of forests. If ${\bf a}\in \Aa_{n,q}$, the corresponding equivalence class is written $\overline{\bf a}$. We also write $Stab ({\bf a})$ for the stabilizer of $\bf a$ in $\Ga_q^{n+2}$. According to the class formula, the number of elements in the equivalence class $\overline {\bf a}$, written $\#({\bf a})$ or $\# (\overline{\bf a})$ is given by:
$$\# ({\overline{\bf a}})=\frac{(q!)^{n+2}}{|Stab ({\bf a})|},$$
where we write $|Stab ({\bf a})|$ for the cardinal of $Stab ({\bf a})$ and, more generally, $|S|$ for the cardinal of any set $S$. 

Explicit formulas for $\# ({\overline{\bf a}})$ and for the various quantities associated to the action of $\Ga_q^{n+2}$ on $\Aa_{n,q}$ such as the number of orbits (and much more) will be given later, and form, from the combinatorial point of view, one of the cores of the article. As it will appear, due to its general form, our analysis paves the way for a systematical combinatorial treatment of all mean field path dependent approximations -e.g. without restrictions on the number of particles at each step of the approximation process.

We are now nearly in a position to give our main results 
on Feynman-Kac particle models, but firstly, to simplify the presentation, need to introduce a pair of definitions.  We start with some more or less traditional, and simplifying 
multi index notation.
For any sequences of integers
  ${\bf p}=(p_k)_{0\leq k\leq n}$, and ${\bf l}=(l_k)_{0\leq k\leq n}$,
 we write $\bf p\leq\bf l$ if and only if $p_k\leq l_k$ for all $0\leq k\leq n$. We write $|{\bf p}|$ for $(p_0+...+p_n)$.
Assuming now that $\bf p\leq\bf l$, we use
the multi-index notation
 $$
   ({\bf l})_{\bf p}=\prod_{0\leq k\leq n}(l_k)_{p_k}~,\quad
 {\bf p}!=\prod_{0\leq k\leq n}p_k!~, \quad\mbox{\rm and}\quad
 s({\bf l},{\bf p})=\prod_{k=0}^n s(l_k,p_k)
 $$
where the $s(l_k,p_k)$ are Stirling numbers of the first kind. 
The difference $(\bf p-\bf l)$, and respectively the addition $(\bf p+\bf l)$ of two sequences is the sequence $(p_k-l_k)_{0\leq k\leq n}$, and respectively $(p_k+l_k)_{0\leq k\leq n}$. When no confusions can arise, we write  $\bf N$, and $\bf q$, for the constant sequences $(N)_{0\leq i\leq n}$, and $(q)_{0\leq i\leq n}$.
We also write ${\bf 1}$, and
respectively ${\bf 0}$, for the sequence of unit integers, and respectively null integers. The above definitions are extended to
infinite sequence of integers  ${\bf p}=(p_k)_{k\geq 0}\in \NN^{\NN}$, with a finite
numbers of strictly positive terms. Any function $\alpha:\NN\mapsto\NN$
on the set of integers into itself,
is extended to integer sequences ${\bf p}\in \NN^{\NN}$, by setting
$\alpha({\bf p})=(\alpha(p_k))_{k\geq 0}$.

For ${\bf a}\in \Aa_{n,q}$, we write $|{\bf a}|$ the sequence $(|a_i|)_{0\leq i\leq n}$. The coalescence sequence of $\bf a$ is the sequence ${\bf q}- |{\bf a}|$ of the coalescence numbers of the $a_i$s. The coalescence degree of $\bf a$ is 
$$
coal({\bf a}):=\left(|{\bf q}|-|{\bf a}|\right)=\sum\limits_{i=0}^n(q-|a_i|)
$$ 
The subset of $\Aa_{n,q}$ of sequences $\bf a$ such that $|{\bf a}|\geq\bf q-{\bf p}$ is written $\Aa_{n,q}({\bf p})$. Notice that, if $\bf a\approx\bf b$, $|{\bf a}|=|{\bf b}|$, so that the notions of coalescence degree and coalescence sequence go over to the set ${\cal F}_{n,q}$ of forests. In particular, notation such as $|\overline{\bf a}|$ or $coal(\overline{\bf a})$ is well-defined.  In view of fla~\ref{invariance}, for any choice $({a}_0,...,{a}_{n})$ of a representative of  a forest ${\bf f}$ in ${\cal F}_{n,q}$, we set
$$
\Delta^{\bf f}_{n,q}=\Delta^{\bf a}_{n,q}\quad\mbox{\rm and}\quad
|{\bf f}|=|{\bf a}|
$$
The subset of ${\cal F}_{n,q}$ associated to $\Aa_{n,q}({\bf p})$ is written ${\cal F}_{n,q}({\bf p})$. 

 \begin{theo}\label{mainthm1}
 The non negative measure $ {\QQ}^{N}_{n,q}$ represents the distribution
 of the coalescent, and non absorbed Markov chain introduced in (\ref{coalmarkov}),  associated with
  a sequence   ${\bf A}=({A}_k)_{0\leq k\leq  n}$ of independent random mappings from $[q]$ into itself, with common distribution $\frac{1}{N^q}\sum_{a\in[q]^{[q]}}\frac{(N)_{|a|}}{(q)_{|a|}}\delta_{a}$. That is, for all $F\in \Ba_b^{sym}(E_n^q)$, we have that
$$
 {\QQ}^{N}_{n,q}(F)=\Delta^{\bf A}_{n,q}(F)=\EE\left(F(X^{{\bf A},c,1}_n,\ldots,
	X^{{\bf A},c,q}_n)~1_{T^{{\bf A},1}\geq n}~\ldots~1_{T^{{\bf A},q}\geq n}\right)
	$$
In addition, for any $1\leq q\leq N$, we have the polynomial
expansion
 $$
  {\QQ}^{N}_{n,q}
 =\gamma^{\otimes q}_n+\sum_{1\leq k\leq (q-1)(n+1)}~\frac{1}{N^{k}}~\partial^k  {\QQ}_{n,q}
$$
with the  collection of signed, and weak derivative
measures $\partial^k  {\QQ}_{n,q}$ given by the formula
$$
\partial^k  {\QQ}_{n,q}=
\sum_{{\bf r}<{\bf q}:|{\bf r}|=k}
  \sum_{{{\bf f}}\in {\cal F}_{n,q}({\bf r})}~
  \frac{s(|{\bf f}|,{\bf q}-{\bf r})~
 \#({\bf f})}{({{\bf q}})_{|{\bf f}|}}~\Delta_{n,q}^{{\bf f}}
$$
\end{theo}
We shall present several consequences of this theorem, including explicit descriptions of the first two order terms
in the polynomial expansion, and a new extension of the Wick product
formula to forests.

The functional expansions of the distributions
$\PP^{N}_{n,q}$ rely on the analysis of
tensor product particle measures on path spaces.
These models can be studied following the same
line of arguments as the ones we used for the distributions
${\QQ}^{N}_{n,q}$.
As a parenthesis, we observe that these two measures coincide
for constant potential functions
$$
\left(\forall n\geq 0\quad\forall x_n\in E_n\quad
 G_n(x_n)=1\right)\Longrightarrow
 \left(\forall n\geq 0\quad\forall q\leq N\quad
 \PP^{N}_{n,q}={\QQ}^{N}_{n,q}\right)
$$
To deal with general potential functions,
we need to re-normalize the particle distributions
$\gamma^N_n$, and to
extend the previous functional analysis to
the path-space distribution of the 
particle  measures defined for any integer sequence
${\bf l}=(l_k)_{0\leq k\leq n}$, by the following formula
$$
\gamma_{n}^{{\bf l},N}=(\gamma_0^N)^{\otimes l_0}\otimes(\gamma_1^N)^{\otimes l_1}\otimes
\cdots\otimes (\gamma_n^N)^{\otimes l_n}
$$
The analytic functional structure
of these path space models relies on natural
combinatorial techniques on
colored rooted trees, and forests.
Our result basically read as follows.
\begin{theo}
For any  $q$,
 the sequence of probability measures
 $(\PP^N_{n+1,q})_{N\geq q}$ is differentiable up to any order. Besides, at any order $N\geq q$ we have the closed formula
$$
\PP^N_{n+1,q}=\eta^{\otimes q}_{n+1}+\sum_{1\leq k<\lfloor (N-q)/2\rfloor}~\frac{1}{N^{k}}~\partial^k\PP_{n+1,q}+\frac{1}{N^{\lfloor (N-q)/2\rfloor}}~\partial^{\lfloor (N-q)/2\rfloor}\PP_{n+1,q}^N
$$
for some collection of signed
derivative measures $\partial^k  \PP_{n,q}$,
whose values can be explicitly described in terms of a class of colored  forests, with  maximal coalescent degree $k$.
\end{theo}
As before, we shall derive several consequences of these expansions, including an explicit description of the first order
term, and a Wick product formula on colored forests.
The latter provides sharp estimates for $\LL_q$-mean errors
between the particle  occupation measures $\eta^N_n$ and the
limiting Feynman-Kac measures $\eta_n$. 
Incidentally, combining
this  Wick product formula with the Borel-Cantelli lemma we 
obtain the almost sure convergence result
$$
\lim_{N\rightarrow\infty}(\eta_n^N)^{\odot q}(F_{n,q})=
0
$$
for any bounded symmetric function $F$ on $E^q_n$,
with $
\int_{E_n}F(x^1_n,\ldots,x^q_n)\eta_n(dx_n^q)=0$.
This result is an extension of the law of large numbers
for $U$-statistics obtained by W. Hoeffding~\cite{hoeffding} for independent and identically distributed
random variables to interacting particle models.

\section{Particle measures expansions on forests}\label{sectroot}
\subsection{A preliminary stochastic tensor product formula}
The link between the two, usual and restricted, tensor products measures, $(\eta_n^N)^{\otimes q}$ and $(\eta_n^N)^{\odot q}$
relies in the end on a simple observation, that will appear to be fundamental for all our forthcoming computations. 
Throughout this section, integers $N\geq q\geq 1$ and a mesurable state space $E$ are fixed once for all.

Consider first the surjection
$$\pi :  [q]^{[q]}\times\langle {q,N}\rangle \longrightarrow {[N]^{[q]}}$$
$$(s,a)
\longmapsto as:=a\circ s $$

\begin{lem}\label{lem1}
Let $b \in [N]^{[q]}$, then the cardinal $|\pi^{-1}(\{b\})|$ of $\pi^{-1}(\{b\})$ only depends on the cardinal $|b|$ of the image of $b$. It is given by
$$
|\pi^{-1}(\{b\})|=(N-|b|)_{q-|b|}(q)_{|b|}
$$
\end{lem}

\begin{proof}
Observe indeed that, with the same notations as above, for $b$ and $a$ fixed, the equation $b=as$ can have a solution only if the image of $b$ is contained in the image of $a$. Let us consider such a $a$. Since $a$ is injective, it has a left inverse ${a} '$ for the composition of maps such that moreover ${a'}{a}$ acts as the identity on $[q]$. The equation has therefore a unique solution $s=a'b$. Since the number $n_p$ of injections from $[q]$ to $[N]$ the image of which contains a fixed subset of $[N]$ of cardinal $p<q$ is
$$n_p=(q)_{p}(N-p)_{q-p},$$
the lemma follows.
\end{proof}

The lemma translates almost immediately into a functional relation between the two tensor product measures.
Notice first the identity
$$
D_{a}D_{b}=D_{ab}
$$
that holds for any pair of mappings $ ({ a},{ b})\in([q]^{[q]})^2$. Let $u,v$ be any linear combinations $\sum\limits_{i\in I}\alpha_ia_i$, (resp. $\sum\limits_{j\in J}\beta_ib_i$) of elements of $[q]^{[q]}$, or, in other terms, elements of the monoid algebra of $[q]^{[q]}$. We extend linearly the definition of $D$, so that:
$$D_u:=\sum\limits_{i\in I}\alpha_iD_{a_i}$$
and 
$$D_uD_v=D_{uv}=\sum\limits_{i\in I, j\in J}\alpha_i\beta_jD_{a_ib_j}$$

We
denote now by $m({\bf x})$ the empirical measure associated with
an $N$-uple ${\bf x}=(x^i)_{1\leq i\leq N}\in E^N$ 
$$
m\left({\bf x} \right)=\frac{1}{N}\sum_{i=1}^{N}\delta_{x^i}
$$
For any integer $q\leq N$, we also consider the 
empirical measures on $E^q$ defined by
\begin{eqnarray*}
m({\bf x})^{\otimes q}&=&\frac{1}{N^q}\sum_{a\in [N]^{[q]}}~
\delta_{(x^{a(1)},\ldots,x^{a(q)})}\\
m({\bf x})^{\odot q}&=&\frac{1}{(N)_q}\sum_{a\in\langle q,N\rangle }~
\delta_{(x^{a(1)},\ldots,x^{a(q)})}
\end{eqnarray*}
Although the notation is self-explanatory, notice that, in the sequel, we will write simply ${\bf x}^{a}$ for  $(x^{a(1)},\ldots,x^{a(q)})$ and $\delta_{{\bf x}^a}$ for $\delta_{(x^{a(1)},\ldots,x^{a(q)})}$.

\begin{cor}\label{cor1}
We have
$$
m({\bf x})^{\otimes q}=m({\bf x})^{\odot q}D_{L^N_q}
$$
where
$$
L_{q}^N=\frac{1}{N^q}~\sum_{a\in [q]^{[q]}}~\frac{(N)_{|a |}}{(q)_{|a|}}~{a}
$$

\end{cor}

\begin{proof}

In view of the following identity, that holds
for any function $F\in\mathcal{B}_{b}^{sym}(E^q)$, with $q\leq N$, $$\delta_{{\bf x}^{a}}D_b(F)=D_b(F)(x^{a(1)},...,x^{a(n)})=F(x^{ab(1)},...,x^{ab(n)})=\delta_{{\bf x}^{ab}}(F)$$
the corollary follows from lemma \ref{lem1} and the identity
$$\frac{1}{(N)_q}(N)_p(N-p)_{q-p}=1\quad\mbox{\rm where} \quad N\geq q\geq p$$
\end{proof}

\begin{cor}
The tensor product measure $m({\bf x})^{\otimes q}$ has a Laurent expansion:
$$m({\bf x})^{\otimes q}=m({\bf x})^{\odot q}\cdot\sum\limits_{0\leq k<q}\frac{1}{N^k}D_{\partial^kL_q}$$
where
$$L_q^N=\sum_{0\leq k<q}~\frac{1}{N^k}~\partial^k L_{q}$$
and
$$\partial^k L_{q}=\sum_{q-k\leq p\leq q}~s(p,q-k)~\frac{1}{(q)_p}\sum_{a\in [q]^{[q]}_p}~a
$$
\end{cor}
Notice that, in the expansion of $\partial^k L_{q}$, the double sum involves only elements $a\in [q]^{[q]}$ the coalescence number of which is less or equal to $k$.

The corollary follows from corollary \ref{cor1} and from the Stirling formula
$$(N)_p=\sum\limits_{1\leq k\leq p}s(p,k)N^k$$

For $q>1$, and recalling that $s(q,q-1)=-\left(\begin{array}{c}{q}\\{2}\end{array}\right)$, we get
\begin{eqnarray*}
\partial^0 L_{q}&=& \frac{1}{q !}\sum_{a\in \langle q,q\rangle}~{a}\\
\partial^1 L_{q}&=&-\frac{q(q-1)}{2}~~\frac{1}{q!}\sum_{a\in \langle q,q\rangle}~{a}+\frac{1}{q!}\sum_{a\in [q]^{[q]}_{q-1}}~{a}\\
\partial^2 L_{q}&=&
s(q,q-2)~\frac{1}{q !}\sum_{a\in \langle q,q\rangle}~{a}
-
\frac{
(q-1)(q-2)}{2}~\frac{1}{q !}\sum_{a\in [q]^{[q]}_{q-1}}~{a}
+
\frac{1}{(q)_{q-2}}~\sum_{a\in [q]^{[q]}_{q-2}}~{a}
\end{eqnarray*}

Notice in particular that, since the composition map
$$ \langle q,q\rangle\times \langle q,N\rangle \longrightarrow \langle q,N\rangle$$
is a surjection with all fibers of cardinal $q!$, we have
$$m({\bf x})^{\odot q}D_{\partial^0L_q}=m({\bf x})^{\odot q}$$
so that, asymptotically in $N$
$$m({\bf x})^{\otimes q}\sim m({\bf x})^{\odot q}$$
More precisely, we have, assuming from now on that $\bf x$ is generic (that is, $x_i\not= x_j$ if $i\not= j$):

\begin{cor}
\label{Cor:firstterm}
The following formulas hold for the expansion of $m({\bf x})^{\otimes q}$
$$N~||m({\bf x})^{\otimes q}-m({\bf x})^{\odot q}||_{tv}=2~\frac{N^q-(N)_q}{N^{q-1}}\underset{N \rightarrow +\infty}{\longrightarrow} q(q-1)$$
$$||m({\bf x})^{\odot q}D_{\partial^kL_q}||_{tv}=\sum\limits_{p=q-k}^q~|s(p,q-k)|~S(q,p)$$

\end{cor}
Here, $||~~||_{tv}$ stands fort the total variation, that is, for any linear operator $Q$ on
$\mathcal{B}_{b}(E^q)$,
$$||Q||_{tv}:=\sup\limits_{||f||_\infty =1}|Q(f)|$$
Recall that the Stirling number of the second kind $S(q,p)$ is the number of partitions of $[q]$ into $p$ non-empty subsets, so that
$$|[q]_p^{[q]}|=S(q,p)(q)_p$$
$$
N^q = \sum_{1\leq k \leq q} S(q,k)(N)_k
$$
The corollary follows then from our previous computations and from the following lemma.

\begin{lem}
\label{Lem:tv}
Let $Q=\sum\limits_{0\leq p\leq q}u_p\sum\limits_{a\in [N]_p^{[q]}}\delta_{{\bf x}^a}$, where the $u_p$ are arbitrary real coefficients, then we have
$$||Q||_{tv}=\sum\limits_{0\leq p\leq q}|u_p|~(N)_p~S(q,p)$$
\end{lem}

\begin{proof}
Indeed, we have $||Q||_{tv}\leq\sum\limits_{0\leq p\leq q}|u_p|~(N)_p~S(q,p)$ by direct inspection, since $||\delta_{{\bf x}^a}||_{tv}=1$ for all $a\in [q]^{[q]}$.
Let us write $sgn(k)$ for the sign of $u_k$, and let us introduce the function $\phi\in\mathcal{B}_{b}(E^q)$ defined by:
$$\phi (y_1,...,y_q)=\sum\limits_{0\leq k\leq q}sgn(k)\delta_{\lambda (y_1,...,y_q)}^k$$
where $\lambda (y_1,...,y_q):=|\{y_1,\dots,y_q  \}|$ and $\delta_{\lambda (y_1,...,y_q)}^k$ stands for the Dirac function.
Then, $Q(\phi )=\sum\limits_{0\leq p\leq q}|u_p|~(N)_p~S(q,p)$, and the lemma follows.
\end{proof}

\subsection{Laurent expansions}\label{sectunn}

The goal of this section is to derive a Laurent expansion of the measures ${\QQ}^{N}_{n,q}\in \Ma(E_n^{q})$, ${\QQ}^{N}_{n,q}(F):=\EE ((\gamma_n^N)^{\otimes q}(F))$, $F\in \Ba_b(E_n^q)$,  with
respect to the population size parameter $N$.
 The following proposition is fundamental. Recall from Sect.~\ref{forestexpansion} that we write ${\cal F}_{n,q}$ for the set of equivalence classes in ${\cal A}_{n,q}$ under the action of the permutation group $\Ga_q^{n+2}$.
\begin{prop}\label{propfu}
For any integers $q\leq N$, any time parameter $n\in\NN$ and any $F\in\Ba_b^{sym}(E_n^q)$
 we have the Laurent expansion
 $$
  {\QQ}^{N}_{n,q}(F)
 =\frac{1}{N^{q(n+1)}}~\sum_{{\bf f}\in {\cal F}_{n,q}}~\frac{({\bf N})_{|{\bf f}|}}{({\bf q})_{|{\bf f}|}}~\# ({\bf f})~\Delta_{n,q}^{\bf f}(F)
 $$
 with  the mappings $\Delta_{n,q}$ introduced in
 (\ref{themapDel}).
\end{prop}

\begin{proof}
Indeed, combining the definition of the particle model
with corollary~\ref{cor1}, we first find that
\begin{eqnarray*}
\EE\left(
\left(\gamma^N_n\right)^{\otimes q}(F)~
\left|~
\xi^{(N)}_{n-1}
\right.
\right)&=&\left(\gamma^N_n(1)\right)^q\times
\EE\left(\left(\eta^N_n\right)^{\otimes q}(F)~\left|~
\xi^{(N)}_{n-1}\right.\right)\\
&=&\left(\gamma^N_n(1)\right)^q
\times
\EE\left(
                 \left(\eta^N_n\right)^{\odot q}
D_{L_q^{N} }(F)
~\left|~
\xi^{(N)}_{n-1}
\right.
\right)\\
&=&\left(\gamma^N_n(1)\right)^q
\times\EE\left(
\Phi_n(\eta^N_{n-1})^{\otimes q}D_{L_q^{N} }(F)
~\left|~
\xi^{(N)}_{n-1}
\right.
\right)\\
\end{eqnarray*}
Using the fact that, conditionnaly to $\xi_{n-1}^{(N)}$ and for $G\in\Ba_b^{sym}(E_n^q)$
$$
\Phi_n(\eta^N_{n-1})^{\otimes q}(G)=
\frac{ 
\left(\eta^N_{n-1}\right)^{\otimes q}
(Q_n^{\otimes q}G)
}{
\left(\eta^N_{n-1}\right)^{\otimes q}
(Q_n^{\otimes q}(1))
}=
\frac{ 
\left(\gamma^N_{n-1}\right)^{\otimes q}
(Q_n^{\otimes q}G)
}{
\left(\gamma^N_{n-1}\right)^{\otimes q}
(Q_n^{\otimes q}(1))
}
$$
and 
$$
\left(\gamma^N_{n-1}\right)^{\otimes q}
(Q_n^{\otimes q}(1))=\left(\gamma^N_{n-1}Q_n(1)\right)^q=
\left(\gamma^N_{n-1}(G_{n-1})\right)^q=\left(\gamma^N_{n}(1)\right)^q
$$
we arrive at 
\begin{eqnarray}
\EE\left(
\left(\gamma^N_n\right)^{\otimes q}(F)~
\left|~
\xi^{(N)}_{n-1}
\right.
\right)&=&
\left(\gamma^N_n(1)\right)^q
\times
\EE\left(
\frac{ 
\left(\gamma^N_{n-1}\right)^{\otimes q}
(Q_n^{\otimes q}D_{L_q^N}F)
}{
\left(\gamma^N_{n-1}\right)^{\otimes q}
(Q_n^{\otimes q}(1))
}~\left|~
\xi^{(N)}_{n-1}
\right.
\right)\nonumber\\
&=&
\EE\left(
\left(\gamma^N_{n-1}\right)^{\otimes q}
(Q_n^{\otimes q}D_{L_q^{N} } F)
~\left|~
\xi^{(N)}_{n-1}
\right.
\right)\label{extendpath}
\end{eqnarray}
Integrating over the past, this yields that
$$
\EE\left(
\left(\gamma^N_n\right)^{\otimes q}(F)
\right)=\EE\left(
\left(\gamma^N_{n-1}\right)^{\otimes q}
(Q_n^{\otimes q}D_{L_q^{N} } F)
\right)
$$
Using a simple induction, we readily obtain the formulae
\begin{eqnarray*}
\EE\left(
\left(\gamma^N_n\right)^{\otimes q}(F)
\right)&=&\EE\left(
\left((\gamma_0^N)^{\otimes q}Q_1^{\otimes q}D_{L_q^N }\ldots Q_n^{\otimes q}D_{L_q^N }\right)(F)
\right)\\
& = & \EE \left( \left(\eta_0^{\otimes q}D_{L_q^N }Q_1^{\otimes q}D_{L_q^N }\ldots Q_n^{\otimes q}D_{L_q^N }\right)(F)
\right)\\
&=&\frac{1}{N^{q(n+1)}}~\sum_{{\bf a}\in \Aa_{n,q}}~\frac{({\bf N})_{|{\bf a}|}}{({\bf q})_{|{\bf a}|}}~\Delta^{\bf a}_{n,q}(F)
\end{eqnarray*}
where $\Aa_{n,q}$ is the set of sequences of maps introduced in section \ref{forestexpansion} and where we use the multi-index notation introduced in the same section.
In view of equation (\ref{invariance}), we know that
$$\Delta_{n,q}^{\bf b}(F)=\Delta_{n,q}^{\bf s(\bf b)}(F)$$
whenever $\bf a$ and $\bf b$ belong to the same equivalence class in $\Aa_{n,q}$, and the Proposition follows.
\end{proof}

We are now in position to derive Theorem \ref{mainthm1}
\begin{theo}
\label{mainthm1-bis}
For any time horizon $n\in\NN$, any population size
$N\geq 1$, and any integer $1\leq q\leq N$, we have the Laurent expansion
\begin{equation}\label{polyf}
 {\QQ}^{N}_{n,q}(F)
 =\gamma^{\otimes q}_n+\sum_{1\leq k\leq (q-1)(n+1)}~\frac{1}{N^{k}}~\partial^k  {\QQ}_{n,q}(F)
\end{equation}
with the signed measures $\partial^k  {\QQ}_{n,q}$ such that
$$\partial^k  {\QQ}_{n,q}(F)=
 \sum_{{\bf r}<{\bf q},~|{\bf r}|=k} ~
  \sum_{{{\bf f}}\in {\cal F}_{n,q}({\bf r})}~
  \frac{s(|{\bf f}|,{\bf q}-{\bf r})~
 \#({\bf f})}{({{\bf q}})_{|{\bf f}|}}~\Delta_{n,q}^{{\bf f}}(F)
$$
for all $F\in \Ba_b^{sym}(E_n^q)$.
\end{theo}

\begin{proof}
Indeed, according to Proposition~\ref{propfu}
$$
 {\QQ}^{N}_{n,q}(F)
 =\frac{1}{N^{q(n+1)}}~\sum_{{\bf p}\in[1,q]^{n+1}}
 \frac{({\bf N})_{\bf p}}{({\bf q})_{{\bf p}}}~
  \sum_{{\bf f}\in {\cal F}_{n,q},|{\bf f}|={\bf p}}~\# ({\bf f})~\Delta_{n,q}^{\bf f}(F)
$$
Using the Stirling formula
\begin{eqnarray}
({\bf N})_{\bf p}
&=&\sum_{{\bf l}\leq\bf p}~s({\bf p},{\bf l})N^{|{\bf l}|}\label{SForm}
\end{eqnarray}
we also have that
$$
 {\QQ}^{N}_{n,p}(F)
 =
 \sum_{{\bf l}\leq{\bf q}}
 \sum_{{\bf l}\leq{\bf p}\leq{\bf q}}~s({\bf p},{\bf l})~\frac{1}{N^{|{\bf q}-{\bf l|}}}
 \frac{1}{({\bf q})_{\bf p}}~
  \sum_{{\bf f}\in {\cal F}_{n,q},|{\bf f}|={\bf p}}~\#({\bf f})~\Delta_{n,q}^{\bf f}(F)
$$
from which we conclude that
$$
 {\QQ}^{N}_{n,q}(F)
 =
 \sum_{{\bf r}<\bf q}
 \sum_{\bf q-r\leq p\leq q}~s({\bf p},{\bf q-r})~\frac{1}{N^{|{\bf r}|}}
 \frac{1}{({\bf q})_{\bf p}}~
  \sum_{{\bf f}\in {\cal F}_{n,q}|{\bf f}|={\bf p}}~\#({\bf f})~\Delta_{n,q}^{\bf f}(F)
$$
Observe then that $\# ({\bf f})=(q!)^{n+1}$ if $|{\bf f}|={\bf q}$ since, in that case, ${\bf f}$ is nothing but the equivalence class of all sequences of bijections in $\Aa_{n,q}$.
This ends the proof of the theorem.
\end{proof}

\begin{cor}
Assuming that $||\Delta_{n,q}^{\bf f}||_{tv}\leq 1$ for any ${\bf f}\in{\cal F}_{n,q}$, the following formula holds for the expansion of ${\QQ}_{n,q}^N$
$$\lim\limits_{N\longrightarrow\infty} N~||{\QQ}^{N}_{n,q}-\gamma^{\otimes q}_n||_{tv}\leq (n+1)q(q-1)$$
\end{cor}

\begin{proof}
 We have (see proof of theorem \ref{mainthm1-bis}) that $\# ({\bf f})=(q!)^{n+1}$ if $|{\bf f}|={\bf q}$ and besides, $|\{{\bf r},~|{\bf r}|=1\}|=n+1$, and for any such $\bf r$, $s({\bf q},{\bf q}-{\bf r})=-{q\choose 2}$. In conclusion, the case $|{\bf f}|={\bf q}$ contributes to $(n+1)q(q-1)/2$ to the asymptotic evaluation of $N~||{\QQ}^{N}_{n,q}-\gamma^{\otimes q}_n||_{tv}$.

Now, let us consider a sequence $\bf r$ with $|{\bf r}|=1$, e.g. the one, written ${\bf r}_i$ with $r_i=1$, for a given $0\leq i\leq n$. In that case, there is a unique $\bf f$ with $|{\bf f}|={\bf q - r}$, which is the equivalence class of all sequences of maps $\bf a$ in $\Aa_{n,q}$ with $a_j\in \Ga_q,~j\not= i$ and $|a_i|=q-1$. In particular $\# ({\bf f})=(q!)^{n+1}{q\choose 2}$, so that, on the whole, the case $|{\bf f}|={\bf r}_i$ contributes to $q(q-1)/2$ to the evaluation.

 The corollary follows then from the Laurent expansion of Theorem \ref{mainthm1-bis}.
\end{proof}

 We conclude this section with noting that
the above functional expansions also apply to the dot-tensor
product measures
$$
\left(\gamma^N_n\right)^{\odot q}(F)=
\left(\gamma^N_n(1)\right)^q\times\left(\eta^N_n\right)^{\odot q}(F)
$$
More precisely, by definition of the particle model we have that
\begin{equation}\label{gamgam}
 \EE\left(\left(\gamma^N_n\right)^{\odot q}(F)\right)=
 \EE\left(\left(\gamma^N_{n-1}\right)^{\otimes q}Q_n^{\otimes q}(F)\right)
 = \QQ^{N}_{n-1,q}(Q_n^{\otimes q}F)
\end{equation}

\section{Combinatorial methods for counting forests and jungles}\label{treesfsec}
In the present section, we face the problem of computing the cardinals $\# ({\bf f})$ involved in our Laurent expansions. We also derive various identities relevant for the fine asymptotical analysis of mean field particle interacting models, such as the number of classes of elements in ${\cal F}_{n,q}$ with a given coalescence degree.

As already pointed out in the Introduction, we have choosen the langage of botanical and genealogical trees, both for technical reasons (since forests appear to be the most natural parametrization of elements in ${\cal F}_{n,q}$), and also for the clarity of the exposition.

\subsection{Some terminology on trees and forests}\label{termino}
In this section, we detail the vocabulary of trees, that will be of constant use later in the article. 
\begin{defi}
A tree (respectively a  planar tree) is a (isomorphism class of) finite
non-empty oriented connected (and respectively planar) graph $T$ without loops such that any
vertex of $T$ has at most one outgoing
edge. Paths are oriented from the vertices to the root.
 There is no restriction on the number of incoming edges.
The empty graph $\emptyset$ is viewed as a tree and is called the empty tree. The set of all trees (resp. planar trees) is denoted by $\cal T$ (resp. $\cal PT$). The vertices of a tree without incoming edges are called the leaves; the vertices with both incoming edges and an outgoing edge are called the internal vertices or the nodes; the (necessarily unique) vertex without outgoing edge is called the root.
\end{defi}

\begin{defi}
A forest ${\bf f} $ is a multiset of trees, that is a set of trees, with repetitions of the same tree allowed or, equivalently, an element of the commutative monoid $\langle \Ta\rangle$ on $\cal T$,
with the empty graph  $T_0=\emptyset$ as a unit. Since the algebraic notation is the most convenient, we write
$${\bf f} =T_1^{m_1}...T_k^{m_k}$$
for the forest with the tree $T_i$ appearing with multiplicity $m_i$, $i\leq k$. When $T_i\not= T_j$ for $i\not= j$, we say that $\bf f$ is written in normal form.
A planar forest ${\bf f}'$ is an ordered sequence of planar trees.  Planar forests can be represented by noncommutative monomials (or words) on the set of planar trees. The sets of forests and planar forests are written $\cal F$ and $\cal PF$. 
\end{defi}
 For example, if $T_1'$ and $T_2'$ are two planar trees, ${\bf f}'=T_1'T_2'T_1'$ is the planar forest obtained by left-to-right concatenation of $T_1'$, $T_2'$ and another copy of $T_1'$. 

\ \

\xymatrix{{\bullet}\ar[dr]&&{\bullet}\ar[dl]&&{\bullet}\ar[dl]&&&&{\bullet}\ar[dr]&&{\bullet}\ar[dl]&&{\bullet}\ar[dl]\\
&{\bullet}\ar[dr]&&{\bullet}\ar[dl]&&{\bullet}\ar[dr]&{\bullet}\ar[d]&{\bullet}\ar[dl]&&{\bullet}\ar[dr]&&{\bullet}\ar[dl]&\\
&&{\bullet}&&&&{\bullet}&&&&{\bullet}&&}

\ \par

\ \par

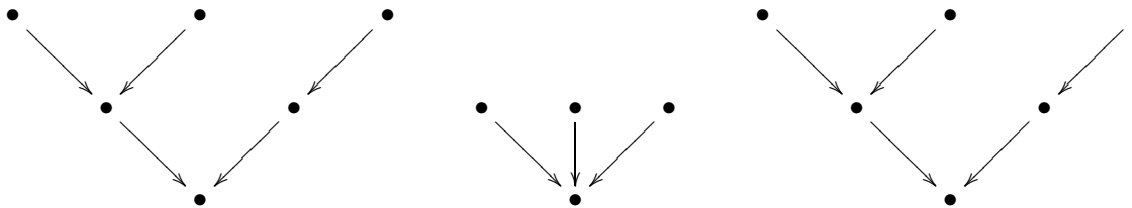
{\sl Figure 1}: a graphical representation of a planar forest ${\bf f}'=T_1'T_2'T_1'$.

\ \par

We say that a tree $T$ is a subtree of $T'$, and we write $T\subset T'$, if the graph of $T$ is a subgraph of the graph of $T'$, and if the root of $T$ is also the root of $T'$. In a more pedantic way, a subtree of a tree is a connected subgraph of $T'$ containing the root. A subforest ${\bf f}\subset {\bf f}'$ of a forest $\bf f'$ is defined accordingly, as a collection of pairwise disjoint subtrees of the trees in $\bf f'$.

The distance between two vertices in a tree is the minimal number of edges of a path joigning them. The height of a vertex is its distance to the root. We also say that a vertex with height $k$ is a vertex at level $k$ in the tree. The height $ht(T)$ of a tree $T$ is the maximal distance between a leaf and the root. The height $ht({\bf f})$ of a forest ${\bf f}$ is the maximal height $ht(T)$ of the trees $T\subset {\bf f}$ in the forest
$$
ht({\bf f})=\sup_{T\subset {\bf f}}{ht(T)}
$$
For the unit empty tree,
we also take the convention $ht(\emptyset)=-1$.

We write $v_k({\bf f})$ for the number of vertices in a forest ${\bf f}$ at level $k\geq 0$. The number of vertices at each level in a forest is encoded in the mapping : 
$$
v~:~ {\bf f}\in \Fa\mapsto v({\bf f})=(v_k({\bf f}))_{k\geq 0}
$$

Notice that  for any pair of forests $({\bf f}, {\bf g})\in\Fa$, we have that
 $$
 v({\bf f} {\bf g})=v({\bf f})+v({\bf g})
 $$
We write ${\bf V}$ for $v({\cal F})$, which coincides with the set of connected initial integer sequences, that is integer sequences ${\bf p}=(p_k)_{k\geq 0}\in \NN^{\NN}$ satisfying the following property
$$
\exists ht({\bf p})\in\NN\quad
\mbox{\rm s.t.}\quad \inf_{k\leq ht({\bf p})}{p_k}>0\quad\mbox{\rm and}\quad
\sup_{k> ht({\bf p})}{p_k}=0
$$

For the unit empty tree $T_0=\emptyset$,
we use the conventions 
$$
v(\emptyset)={\bf 0}=(0,0,\ldots)\quad\mbox{\rm and}\quad
ht({\bf 0})=-1
$$

For any ${\bf p}\in {\bf V}$, 
 we write ${\cal F}_{\bf p}$ for the set of forests ${\bf f}$ such that $v({\bf f})={\bf p}$. Our notation is consistent, since the height $ ht({\bf f})$ of a forest ${\bf f}\in {\cal F}_{\bf p}$ clearly coincides with the height $ht({\bf p})$ of the integer 
sequence ${\bf p}$.
For the null sequence ${\bf 0}$, we have that ${\cal F}_{\bf 0}=\{\emptyset\}$.

When ${\bf p}\in {\bf V}$ is chosen so that $p_0=1$, 
the set ${\cal F}_{\bf p}$ reduces to the set ${\cal T}_{\bf p}$ of all trees $T$ such that $v(T)={\bf p}$. The notation on trees, and forests goes over to planar trees, and planar forests in a self-explanatory way. For example, for any sequence 
${\bf p}\in{\bf V}$, ${\cal PF}_{\bf p}$ is the set of planar forests ${\bf f}$, with $v_k({\bf f})=p_k$ vertices at level $k\geq 0$. 

For any $n\in\left(\NN\cup\{-1\}\right)$, we denote by ${\bf V}_n\subset{\bf V}$
the subset of sequences ${\bf p}$ such that $ht({\bf p})=n$. Notice that
${\bf V}_{-1}=\{{\bf 0}\}$. In this notation, the sets of all forests, and trees with height $n$ are given by the sets 
$$
{\cal F}(n):=\coprod_{{\bf p}\in V_n}\Fa_p,~~{\cal T}(n):=\coprod_{{\bf p}\in V_n,p_0=1}{\cal T}_p
$$

The shift operator
$$
B~:~{\bf p}=(p_k)_{k\geq 0}\in {\bf V}\mapsto
B({\bf p})=(q_l)_{l\geq 0}\in {\bf V},~q_l:=p_{l+1}
$$
induces a canonical bijection, still denoted by $B$, between the set of trees $\Ta_{\bf p}$ (resp. planar trees $\Pa\Ta_{\bf p}$), with ${\bf p}\in {\bf V}_{n+1}$ 
and the set of forests $\Fa_{B({\bf p})}$ (resp. planar forests $\Pa\Fa_{
B({\bf p})}$), obtained by removing the root of the tree
$$
B~:~T\in \Ta_{\bf p}~\left(\mbox{\rm resp.}~\Pa\Ta_{\bf p}\right)\mapsto B(T)\in\Fa_{B({\bf p})}~\left(\mbox{\rm resp.}~\Pa\Fa_{B({\bf p})}\right)
$$

From the graph-theoretic point of view, the number of coalescences $c_k({\bf f})$ at each level $k\geq 0$ in a forest ${\bf f}\in\Fa$ is defined by the  mapping
$$
c~:~{\bf f}\in\Fa\mapsto c({\bf f})=B(v({\bf f}))-|{\bf f}|
$$
The sequence $c({\bf f})$ is called the coalescence sequence of $\bf f$.
In the above displayed formula, 
the sequence $|{\bf f}|=(|{\bf f}|_k)_{k\geq 0}$ represents the number of vertices at level $k\geq 0$ minus the number of leaves at the same level, that is the number of vertices with an ingoind edge, whereas $B(v({\bf f}))_k$ stands, by definition, for the number of edges between the levels $k+1$ and $k$.
Finally observe that for any pair of forests $({\bf f},{\bf g})\in\Fa$, we have
$$
\left(~|{\bf f}{\bf g}|=|{\bf f}|+|{\bf g}|\quad\mbox{\rm and}\quad
B(v({\bf f})+v({\bf g}))=B(v({\bf f}))+B(v({\bf g}))
~\right)\Rightarrow 
c({\bf f}{\bf g})=c({\bf f})+c({\bf g})
$$
We also say that the coalescence order of a vertex in a tree or a forest is its number of incoming edges minus 1. The coalescence degree of a tree or a forest $\bf f$ is the sum of the $c_k({\bf f})$ or, equivalently, the sum of the coalescence orders of its vertices. We say that a tree is \it trivial \rm if its coalescence degree is 0.

\subsection{On the algebraic structure of trees and forests}\label{algb}

In a planar forest $\bf f'$, vertices at the same level $k\geq 0$, are naturally ordered from left to right, and therefore in bijection with $[v_k({\bf f'})]$. Planar forests (resp. planar trees) ${\bf f}'\in \cal PF$ of height $ ht({\bf f}')=(n+1)$ are therefore canonically in bijection with 
sequences of maps ${\bf a}=(a_0,...,a_{n})$, where $a_k$ is a weakly increasing map from $[p_{k+1}]$ to $[p_{k}]$, with $k=0,...,n$, and ${\bf p}\in{\bf V}_{n+1}$ (resp. and $p_0=1$).

To clarify the presentation, we
write $a\uparrow$ to indicate that a given mapping
$a$ between sets of integers is weakly increasing.
 In this notation, for any ${\bf p}\in{\bf V}_{n+1}$ we have that
$$
{\cal PF}_{\bf p}\simeq {\cal IA}_{\bf p}:=\left\{(a_0,...,a_{n})\in \left([p_0]^{[p_1]}\times\ldots\times [p_n]^{[p_{n+1}]}\right)~:~\forall 0\leq k\leq n~~a_k\uparrow\right\}
$$
We write $PF({\bf a})$ (resp. $PT({\bf a})$) for the planar forest associated to a sequence in ${\cal IA}_{\bf p}$ (resp. with $p_0=1$).
In this interpretation, when ${\bf f}'\in{\cal PF}_{\bf p}$ is represented by a sequence of mappings ${\bf a}=(a_0,...,a_{n})\in {\cal IA}_{\bf p}$, the number of coalescent vertices at level $k$ is given by
$$
c_k({\bf f}')=v_{k+1}({\bf f}')-|a_k|
$$
Notice that the definition is coherent with the definition of the coalescence sequence associated to $\bf a$.
In addition, when the integer sequence
 ${\bf p}\in{\bf V}_{n+1}$ is chosen such that $p_0=1$, the mapping
$B$ introduced in the previous Section can be rewritten as follows:
$$
B~:~PT({\bf a})\in{\cal PT}_{\bf p}\mapsto B(PT({\bf a})):=PF(B({\bf a}))
\in{\cal PF}_{B({\bf p})}
$$
with the shift operator $B(a_0,\ldots,a_n):=(a_1,\ldots,a_n)$.
These definitions and notations can be extended to arbitrary sequences of maps and arbitrary trees and forests as follows.

 \begin{defi}
Let ${\bf a}=(a_0,a_1,...,a_n)$ be any sequence of maps $a_k:[p_{k+1}]\rightarrow [p_k]$, with $k=0,\ldots,n$. To a sequence ${\bf a}$ is naturally associated a forest $F({\bf a})$: the one with one vertex for each element of $\coprod_{k=0}^{n+1} [p_k]$, and a edge for each pair $(i,a_k(i)), i\in [p_{k+1}]$. The sequence can also be represented graphically uniquely by a planar graph $J({\bf a})$, where however the edges between vertices at level $k+1$ and $k$ are allowed to cross.
We call such a planar graph, where paths between vertices are entangled, a jungle. The set of jungles is written $\cal J$. 
\end{defi}
 
\ 
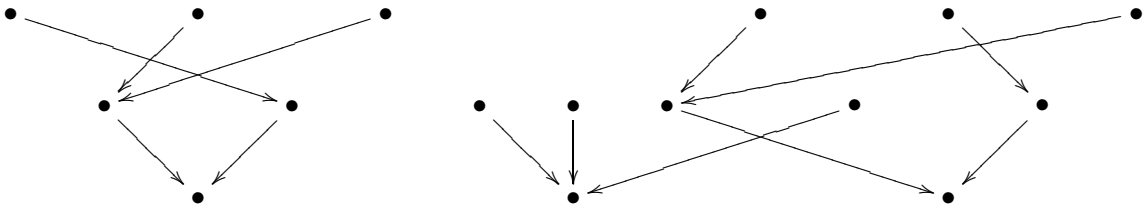
\par

\xymatrix{{\bullet}\ar[drrr]&&{\bullet}\ar[dl]&&{\bullet}\ar[dlll]&&&&{\bullet}\ar[dl]&&{\bullet}\ar[dr]&&{\bullet}\ar[dlllll]\\
&{\bullet}\ar[dr]&&{\bullet}\ar[dl]&&{\bullet}\ar[dr]&{\bullet}\ar[d]&{\bullet}\ar[drrr]&&{\bullet}\ar[dlll]&&{\bullet}\ar[dl]&\\
&&{\bullet}&&&&{\bullet}&&&&{\bullet}&&}

{\it Figure 2:} The entangled graph representation of a jungle with the same underlying graph as the planar forest in Fig. 1.

Notice that planar forests are particular examples of jungles, so that, for any sequence of weakly increasing maps ${\bf a}$ we have $PF({\bf a})=J({\bf a})$. Notation such as ${\cal F}_{\bf p}$, with ${\bf p}\in{\bf V}_{n+1}$, extends from forests to jungles in a self-explanatory way. Jungles $j\in \cal J_{\bf p}$ of height $(n+1)$ are therefore canonically in bijection with sequences of maps $(a_0,...,a_{n})$, where $a_k$ is a map from $[p_{k+1}]$ to $[p_{k}]$, with $k=0,...,n$. More formally, we have that
$$
{\cal J}_{\bf p}\simeq {\cal A}_{\bf p}:=\left([p_0]^{[p_1]}\times\ldots\times [p_n]^{[p_{n+1}]}\right)
$$

 To make these algebraic representations of planar forests and jungles more precise, it is convenient to introduce another round of notation.
For any integers $p,q,r,s\geq 1$, and any pair of maps 
$(\alpha,\beta)\in \left([q]^{[p]}\times [s]^{[r]}\right)$, we let
$(\alpha\vee\beta)$ be the weakly increasing map from $[p+r]$ into $[q+s]$ defined
by
$$
(\alpha\vee\beta)(i)=\left\{
\begin{array}{rcl}
\alpha(i)&\mbox{\rm if}& i\in [p]\\
q+\beta(i)&\mbox{\rm if}& i=p+j~\mbox{\rm with}~j\in [r]
\end{array}
\right.
$$
If $\alpha$ and $\beta$ are weakly increasing, so is $\alpha \vee\beta$.
For null integers $p=0$, or $p=0=q$ we use the conventions $[q]^{[p]}=\{\emptyset\}$, and the composition rules
$
(\emptyset\vee\emptyset)=\emptyset$, 
$(\alpha\vee\emptyset)=\alpha=(\emptyset\vee\alpha)
$.
Observe that for any $(\sigma,\tau)\in \left(\Ga_q\times\Ga_s\right)$, any $(\alpha,\beta)\in \left([q]^{[p]}\times [s]^{[r]}\right)$ and
any $(\rho,\theta)\in \left(\Ga_p\times\Ga_r\right)$ we have the composition formula
\begin{eqnarray}
(\sigma\vee\tau)\circ(\alpha\vee\beta)\circ(\rho\vee\theta)&=&(\sigma\circ\alpha\circ\rho)\vee(\tau\circ\beta\circ\theta)\label{techs0}
\end{eqnarray}
For any $\rho\in\Ga_{p+r}$ with $\alpha$ and $\beta$ as above, we also quote the following factorization rule
$$
(\alpha\vee\beta)\circ\rho\in \left([q]^{[p]}\vee [s]^{[r]}\right) \Longrightarrow
\rho\in (\Ga_p\vee\Ga_r)
$$
To check this final assertion, we suppose there exists some
$i\in[p]$ such that $\rho(i)=p+j$, for some $j\in[r]$. In this case, we would find the following contradiction
$$
(\alpha\vee\beta)(\rho(i))=(\alpha\vee\beta)(p+j)=q+\beta(j)\not\in [q]
$$
By a simple induction argument, for any $\alpha_i\in [q_i]^{[p_i]},i=1,...,m$, we also have that
\begin{equation}\label{factorule}
(\alpha_1\vee\ldots\vee \alpha_m)\circ\rho\in \left([q_1]^{[p_1]}\vee\ldots\vee  [q_m]^{[p_m]}\right) \Longrightarrow
\rho\in (\Ga_{p_1}\vee\ldots\vee\Ga_{p_m})
\end{equation}
for any permutation $\rho\in \Ga_{p_1+\ldots+p_m}$, and $m\geq 1$.

More generally, for any sequences of integers ${\bf p} \in{\bf V}$, and ${\bf r} \in{\bf V}$, and any sequences of maps
 $
  \alpha=(\alpha_k)_{k\geq 0}\in {\cal A}_{\bf p}$, and $
   \beta=(\beta_k)_{k\geq 0}\in {\cal A}_{\bf r}  
  $
we denote by $(\alpha\vee\beta)$ the sequence of maps
 $$
\alpha\vee\beta:= (\alpha_k\vee\beta_k)_{k\geq 0}\in {\cal A}_{\bf p+r}  
$$
Notice that, if $\alpha\in{\cal IA}_{\bf p}$ and $\beta\in{\cal IA}_{\bf q}$, then $\alpha\vee\beta\in{\cal IA}_{\bf p+q}$. 
In the above displayed formula, we have used the conventions
$
\alpha_k=\emptyset
$, and $
\beta_l=\emptyset
$,
for any $k>ht({\bf p})$, and $l>ht({\bf r})$.

Recall that sequences of weakly increasing maps $\alpha\in {\cal IA}_{\bf p}$, and $\beta\in{\cal IA}_{\bf r}$ can be represented graphically by the planar forests $PF(\alpha)\in {\cal PF}_{\bf p}$, and $PF(\beta)\in {\cal PF}_{\bf r}$. In this interpretation, the composition formula on the noncommutative
monoid $ \cal PF$ is now given by 
 \begin{equation}\label{talpha}
PF(\alpha)PF(\beta)=PF(\alpha\vee\beta)\in {\cal PF}_{\bf p+r}
\end{equation}

The same constructions hold for non necessarily increasing maps. In this situation, the sequences of maps $\alpha\in {\cal A}_{\bf p}$, and $\beta\in {\cal A}_{\bf r}$ can be represented graphically by the jungles $J(\alpha)\in {\cal J}_{\bf p}$, and $J(\beta)\in {\cal J}_{\bf r}$, and we have the (noncommutative) product formula on the set $ \cal J$
\begin{equation}\label{talphaj}
J(\alpha)J(\beta)=J(\alpha\vee\beta)\in {\cal J}_{\bf p+r}
\end{equation}

\subsection{Automorphism groups on jungles}\label{autojung}
For any given sequence of integers ${\bf p}\in{\bf V}_{n+1}$, the product permutation group 
$$
{\Ga}_{\bf p}:=\left({\Ga}_{p_0}\times {\Ga}_{p_1}\times ...\times {\Ga}_{p_{n+1}}\right)
$$ 
acts naturally on sequences of maps ${\bf a}=(a_0,a_1,...,a_n)\in
{\cal A}_{\bf p}$, and on jungles $J({\bf a})$ in ${\cal J}_{\bf p}$ by permutation of the vertices at each level. More formally,  for any ${\bf s}=(s_0,...,s_{n+1})\in \Ga_{\bf p}$ this pair of actions is given by
$${\bf s}({\bf a}):=(s_0a_0s_1^{-1},...,s_na_ns_{n+1}^{-1})
\quad\mbox{\rm and}\quad
{\bf s}J({\bf a}):=J({\bf s}({\bf a}))
$$
An automorphism ${\bf s}\in {\Ga}_{\bf p}$ of a given jungle $J({\bf a})\in {\cal J}_{\bf p}$ is a sequence of permutations that
preserves the jungle, in the sense that
${\bf s}J({\bf a})=J({\bf a})$. The set of automorphisms of a given jungle $J({\bf a})\in {\cal J}_{\bf p}$ coincides with the stabilizer
of  $J({\bf a})$ with respect to the group action. 
$$
Stab(J({\bf a}))=\left\{{\bf s}\in {\Ga}_{\bf p}~:~{\bf s}J({\bf a})=J({\bf a})\right\}
$$
By the definition of
a jungle, we also have that 
$$
J({\bf a})=J({\bf b})\Longleftrightarrow {\bf a}={\bf b}
$$
from which we find $Stab(J({\bf a})) $ coincide with the stabilizer 
of  ${\bf a}$ with respect to the group action. The stabilizer theorem also provides a one to one correspondence
$\varphi_{{\bf a}}$
between the orbit set
$$
Orb(J({\bf a})):=\left\{{\bf s}J({\bf a})~:~{\bf s}\in {\Ga}_{\bf p}\right\}\subset {\cal J}_{\bf p}
$$
and the quotient ${\Ga}_{\bf p}/Stab(J({\bf a}))$
$$
\varphi_{{\bf a}}~:~J({\bf b})\in Orb(J({\bf a}))\mapsto
{\bf s}_{({\bf a},{\bf b})}Stab(J({\bf a}))\in {\Ga}_{\bf p}/Stab(J({\bf a}))
$$
for any choice of ${\bf s}_{({\bf a},{\bf b})}\in {\Ga}_{\bf p}$ such that
${\bf s}_{({\bf a},{\bf b})}({\bf a})={\bf b}$.

Notice that if two sequences $\bf a$ and $\bf b$ differ only by the order of the vertices in $J({\bf a})$ and $J({\bf b})$, that is by the action of an element of ${\Ga}_{\bf p}$, then the associated forests are identical: $F({\bf a})=F({\bf b})$. Moreover, the converse is true: if $F({\bf a})=F({\bf b})$, then $J({\bf a})$ and $J({\bf b})$ differ only by the ordering of the vertices, since they have the same underlying non planar graph. In this situation, $\bf a$ and $\bf b$ belong to the same orbit under the action of ${\Ga}_{\bf p}$. In particular, we get a result already mentioned in the Introduction

\begin{lem}

The set of equivalence classes of jungles in ${\cal J}_{\bf p}$ under the action of the permutation groups ${\Ga}_{\bf p}$ is in bijection with the set of forests ${\cal F}_{\bf p}$. In particular, due to the equivalence between the two notions of jungles and sequences of maps, it follows that
the set ${\cal F}_{n,q}$ of equivalence classes in $\Aa_{n,q}$ under the natural action of ${\bf \Ga_q}$ is canonically in bijection with the set of forests with $q$ vertices at each level $0\leq k\leq n+1$. 
\end{lem}
Here, $\bf q$ denotes, once again, the $q$-constant sequence of lenght $n+1$.

Since the root mapping $a_0$ 
of a planar tree is constant and equal to $1$,
we finally observe that for any
${\bf a}=(a_0,\ldots,a_n)\in {\cal IA}_{\bf p}$, we have
$$
Stab(J({\bf a})=Stab(PT({\bf a}))=\{1_1\}\times Stab(B(PT({\bf a})))
$$

\subsection{Orbit sets of planar forests}\label{orbitsec}

From the Feynman-Kac mean-field approximation point of view, jungles in $\Aa_{n,q}$ describe trajectories of families of particles with prescribed coalescence properties. In particular, two trajectories that are equivalent under the action of $\Ga_{\bf q}$ have the same statistical properties.
In the present section, we face the general problem of computing the cardinals $\# ({\bf f})$, defined as the number of jungles in ${\bf f}\in {\cal F}_{\bf p}$,  for some integer sequence ${\bf p}\in{\bf V}_{n+1}$, when $\bf f$ is viewed as an equivalence class according to the previous lemma. Let us write ${\bf f}$ in normal form (with the $T_i$s all distincts) as a commutative monomial of trees:
$${\bf f}=T_1^{m_1}...T_k^{m_k},~m_1+...+m_k=p_0$$
Let us also choose arbitrary representations $T_i'\in{\cal PT}$ of the $T_i$s as planar trees, so that we can view ${\bf f}$ as the forest associated to the planar forest 
$${\bf f}'=(T_{1}')^{m_1}...(T_{k}')^{m_k}
$$ obtained by left-to-right concatenation of $m_1$ copies of $T_1'$, $m_2$ copies of $T_2'$,..., $m_k$ copies of $T_k'$.  We write $T_1''...T_{p_0}''$ for the expansion of the noncommutative monomial 
$$
{\bf f}'=(T_{1}')^{m_1}...(T_{k}')^{m_k}=T_1''...T_{p_0}''
$$ as a product of planar trees (without exponents).

From previous
considerations, we have that
$$
\# ({\bf f})=\# Orb({\bf f'})=\#\left\{ {\bf g}\in {\cal J}_{\bf p}~:~\exists {\bf s}\in\Ga_{\bf p} \quad\mbox{\rm s.t.}\quad {\bf g}={\bf s}{\bf f'} \right\}
$$
Due to the class formula, we also know that
$$\# ({\bf f})=\# Orb({\bf f'})=\frac{|{\Ga}_{\bf p}|}{Stab({\bf f}')}=\frac{{\bf p}!}{Stab({\bf f}')}$$
where we remind that ${Stab({\bf f}')}$ stands for the stabilizer of the jungle ${\bf f'}\in {\cal J}_{\bf p}$. 

The computation of $\# ({\bf f})$ will be done by induction on $(n+1):=ht({\bf f})$. Let us assume that we know $\# ({\bf g})$ for any forest of height less or equal $n$. Notice that, then, we also know $\# ({\bf t})$ for any tree ${\bf t}$ of height $(n+1)$, due to the canonical bijection $B$ between trees and forests: $\# ({\bf t})=\# (B{\bf t})$. From the previous discussion,
the problem amounts to compute the cardinals of the stabilizers ${Stab({\bf f}')}$ inductively with respect to the height parameter. The following technical lemma is instrumental.

\begin{lem}\label{keyle}
There is a natural isomorphism
$$
Stab(T_{1}''...T_{p_0}'')\cap\left(\{1_{p_0}\}\times\Ga_{(p_1,\ldots,p_{n+1})}\right)
\sim Stab(T_1'')\vee\ldots \vee Stab(T_{p_0}'')
$$
where $1_{p_0}$ stands for the identity in $\Ga_{p_0}$.
\end{lem}
\proof
We let ${\bf a}^i=(a^i_0,\ldots,a^i_n)$ be a sequence of weakly increasing maps such that
$T''_i=PT({\bf a}^i)$, with $1\leq i\leq p_0$. In this notation, we have that
\begin{eqnarray*}
{\bf f}'&=&T_1''...T_{p_0}''=PT({\bf a}^1)\ldots PT({\bf a}^{p_0})=
PF({\bf a})
\end{eqnarray*}
where $PF({\bf a})$ stands for the planar
forest associated with the sequence of weakly increasing maps 
$$
{\bf a}=({\bf a}^1\vee\ldots\vee {\bf a}^{p_0})=(a_0,\ldots,a_n)
\quad\mbox{\rm 
with}\quad
\forall 0\leq l\leq n\qquad
a_l= a^1_l\vee \ldots\vee a^{p_0}_l
$$
Notice that the stabilizer of the planar forest $PF({\bf a})$ is given by
$$
Stab(PF({\bf a}))=\left\{{\bf s}\in {\Ga}_{\bf p}~:~{\bf s}({\bf a})={\bf a}\right\}
$$
and for any ${\bf s}=(s_0,\ldots,s_{n+1})\in  {\Ga}_{\bf p}$ we have that
$$
{\bf s}({\bf a})={\bf a}\Leftrightarrow\left(
\forall 0\leq l\leq n\qquad
s_l \circ(a^1_l\vee \ldots\vee a^{p_0}_l)\circ s_{l+1}^{-1}=a^1_l\vee \ldots\vee a^{p_0}_l\right)
$$

Therefore, using the composition formula (\ref{techs0}) and the factorization rule following the formula, it is not difficult to check that any sequence
$$
{\bf s}=(1_{p_0},s_1,\ldots,s_{n+1})\in Stab({\bf f}')
$$
can be represented in an unique way 
$
{\bf s}=\left({\bf s}^{(1)}\vee \ldots\vee {\bf s}^{(p_0)}\right)
$,
in terms of some permutation sequences 
$$
{\bf s}^{(i)}=(1_1,s^{(i)}_1,\ldots,s^{(i)}_{n+1})
$$
between the vertices at each level of the planar tree $PT({\bf a}^{i})$, with $1\leq i\leq p_0$.  In terms of the level parameter  $1\leq l\leq (n+1)$, we also notice that
$$
s_l=s^{(1)}_l\vee \ldots\vee s^{(p_0)}_l
$$ 
To take the final step, we observe that
\begin{eqnarray*}
s_l a_ls^{-1}_{l+1}&
=&\left(s^{(1)}_l\vee \ldots\vee s^{(p_0)}_l\right)\circ
(a_l^1\vee\ldots\vee a^{p_0}_l)\circ \left((s^{(1)}_{l+1})^{-1}\vee \ldots\vee (s^{(p_0)}_{l+1})^{-1}\right)\\
&=&(s^{(1)}_l\circ a_l^1 \circ(s^{(1)}_{l+1})^{-1})
\vee\ldots\vee (s^{(p_0)}_l\circ a^{p_0}_l\circ(s^{(p_0)}_{l+1})^{-1})
\end{eqnarray*}
This implies that
$$
{\bf s} ({\bf a})=\left({\bf s}^{(1)}\vee \ldots\vee {\bf s}^{(p_0)}\right)
({\bf a}^1\vee\ldots\vee {\bf a}^{p_0})
=({\bf s}^{(1)}({\bf a}^1)\vee\ldots\vee {\bf s}^{(p_0)}({\bf a}^{p_0}))
$$
from which we find that ${\bf s} J({\bf a})=J({\bf a})$ if, and only if, we have
\begin{eqnarray*}
J({\bf s}^{(1)}({\bf a}^1))\ldots J({\bf s}^{(p_0)}({\bf a}^{p_0}))&=&
J({\bf a}^1)\ldots J({\bf a}^{p_0})
\end{eqnarray*}
This is clearly equivalent to the fact that
$$
J({\bf s}^{(i)}({\bf a}^i))={\bf s}^{(i)}J({\bf a}^i)=J({\bf a}^i)
$$
for any $1\leq i\leq p_0$.
This yields  that for any $1\leq i\leq p_0$, we have
$$
{\bf s}^{(i)}=(1_1,s^{(i)}_1,\ldots,s^{(i)}_{n+1})\in Stab(PT({\bf a}^{i}))$$
from which we conclude that the mapping
$$
\begin{array}{c}
{\bf s}\in Stab(PT({\bf a}^{1})\ldots PT({\bf a}^{p_0}))\cap\left(\{1_{p_0}\}\times\Ga_{(p_1,\ldots,p_{n+1})}\right)\\
\\
\downarrow\\
 \\
 \left({\bf s}^{(1)}\vee \ldots\vee {\bf s}^{(p_0)}\right) \in Stab(T({\bf a}^1))\vee\ldots \vee  Stab(T({\bf a}^{p_0}))
\end{array}$$
is an isomorphism. This ends the proof of the lemma.
\cqfd
The first coordinate mapping  
$$
\pi~:~{\bf s}=(s_0,\ldots,s_{n+1})\in Stab((T_{1}')^{m_1}...(T_{k}')^{m_k})\mapsto
\pi({\bf s})=s_0\in \left(\Ga_{m_1}\vee\ldots\vee\Ga_{m_k}\right)
$$
is a surjective map (it even has a natural section, the construction of which is omitted). We therefore have that
$$
Ker(\pi)=\pi^{-1}(1_{p_0})\sim Stab(T_1'')\vee\ldots \vee Stab(T_{p_0}'')
\sim Stab(T_1')^{m_1}\times\ldots \times Stab(T_k')^{m_k}
$$
This yields the isomorphism formula
$$
Stab((T_{1}')^{m_1}...(T_{k}')^{m_k})/
\left(Stab(T_1')^{m_1}\times\ldots \times Stab(T_k')^{m_k}\right)
\sim \left(\Ga_{m_1}\vee\ldots\vee\Ga_{m_k}\right)
$$
from which we readily deduce the 
following recursive formula for the computation of $\# (\bf f)$.
\begin{prop}\label{recursion}
We have
$$|Stab((T_{1}')^{m_1}...(T_{k}')^{m_k})|=\prod\limits_{i=1}^k \left(m_i!~|Stab ~(B(T_i'))|^{m_i}\right)$$
\end{prop}

\subsection{An inductive method for counting jungles}\label{closed}
We conclude this section by giving a closed formula for $\# ({\bf f})$. Let us introduce for that purpose some further notations. 
\begin{defi}
Let ${\bf f}$ be a forest written in normal form ${\bf f}=T_1^{m_1}...T_k^{m_k}$. In this situation, we say that the unordered $k$-uplet
$(m_1,...,m_k)$ is the symmetry multiset of the tree $T=B^{-1}({\bf f})$ and write 
$$
{\bf s}(T)={\bf s}(B^{-1}( T_1^{m_1}...T_k^{m_k}))=(m_1,...,m_k)
$$ 
The symmetry multiset of the forest  ${\bf s}({\bf f})=T_1^{m_1}...T_k^{m_k}$  is the disjoint union 
of the symmetry multisets of its trees
$$
{\bf s}(T_1^{m_1}...T_k^{m_k})=\left(\underbrace{{\bf s}(T_1),\ldots,{\bf s}(T_1)}_{m_1~\mbox{\rm terms}},\ldots,\underbrace{{\bf s}(T_k),\ldots,{\bf s}(T_k)}_{m_k~\mbox{\rm terms}}\right)
$$
\end{defi}
For example, the symmetry multiset of the first tree of the forest $\bf f$ displayed in Fig. 1 is $(1,1)$; the symmetry multiset of $B^{-1}({\bf f})$ is $(2,1)$; and the symmetry multiset of $\bf f$ is $(1,1,3,1,1)$.

We also extend the map $B$ from trees to forests to a map between forests and, if ${\bf f}=T_1^{m_1}...T_k^{m_k}$ is a forest, set
$$B({\bf f})=B(T_1)^{m_1}...B(T_k)^{m_k}$$
From the point of view of graphs, the operations amounts to removing all the roots and all the edges that have the root as terminal vertex from the graph defining ${\bf f}$.
The definition extends naturally to planar trees and planar forests. In particular, if $\bf f'$ is a planar forest with $\bf f$ as its underlying forest, we set ${\bf s}({\bf f'}):={\bf s}({\bf f})$.

\begin{theo}\label{numberofjungles}
The number of jungles (or, equivalently, of Feynman-Kac mean-field type trajectories with the same statistics) in ${\bf f}\in {\cal F}_{\bf p}$, with ${\bf p}\in{\bf V}_{n+1}$,
is given by
$$\# ({\bf f})=\frac{{\bf p}!}{\prod\limits_{i={-1}}^{n}{\bf s}(B^i({\bf f}))!}$$
where we use the usual multi-index notation to define ${\bf s}(B^i({\bf f}))!$.
\end{theo}

\proof
It clearly suffices to prove that
$$
|Stab ({\bf f}')|=\prod\limits_{i={-1}}^{n}{\bf s}(B^i({\bf f}'))!
$$
for any planar forest ${\bf f}'$ of height $(n+1)$. We check this assertion by induction
on the height parameter. 
First, we observe that a planar forest ${\bf f}'$ of height $1$ can be represented as a noncommutative monomial $$
{\bf f}'=(T_1')^{m_1}...(T_k')^{m_k}
$$
with different planar trees $T_i'$ of height $1$, and some sequence
of integers  $m_i$. In that case, the planar forests $B(T_i')$ reduce to elementary planar forests with null height, with $S(T_i')$ elementary
planar trees with null height. This yields that
$$
ht(T_i)=1\Longrightarrow |Stab (B(T_i'))|={\bf s}(T_i')!
$$
By proposition \ref{recursion}, we conclude that
$$
|Stab ({\bf f}')| =  \prod_{j=1}^k\left(m_j!~ ({\bf s}(T_j')!)^{m_j}\right)=
 {\bf s}(B^{-1}({\bf f}'))!~{\bf s}({\bf f}')!
$$
This ends the proof of the formula at rank $1$. Suppose now that the 
assertion is satisfied for any planar forest with height $n$.
By proposition \ref{recursion}, for any planar forest ${\bf f}'=(T_1')^{m_1}...(T_k')^{m_k}$, with height 
$(n+1)$, and written in terms of distinct planar trees $T_i$, we have
$$
|Stab ({\bf f}')| =  {\bf s}(B^{-1}({\bf f}'))! \prod_{i=1}^k |Stab (B(T_i'))|^{m_i}
$$
Since the planar forests $B(T_i')$ have height $n$, the induction  hypothesis implies that
$$
 |Stab (B(T_i'))|=\prod_{j=-1}^{n-1}{\bf s}(B^j(B(T_i')))=
 \prod_{j=0}^{n}{\bf s}(B^j(T_i'))!
$$
This yields that
$$
|Stab ({\bf f}')| =  {\bf s}(B^{-1}({\bf f}'))! \prod_{i=1}^k \left(
\prod_{j=0}^{n}{\bf s}(B^j(T_i'))!
\right)^{m_i}={\bf s}(B^{-1}({\bf f}'))!  \prod_{j=0}^{n}\left(\prod_{i=1}^k
{\bf s}(B^j(T_i'))!
\right)^{m_i}
$$
Recalling that
$
B^j({\bf f}')=B^j(T_1')^{m_1}...B^j(T_k')^{m_k}
$, 
we also find that
$$
{\bf s}(B^j({\bf f}'))=(\underbrace{{\bf s}(B^j(T_1')),\ldots,{\bf s}(B^j(T_1'))}_{m_1~\mbox{\rm terms}},\ldots,\underbrace{{\bf s}(B^j(T_k')),\ldots,{\bf s}(B^j(T_k'))}_{m_k~\mbox{\rm terms}})
$$
We conclude that
$$
|Stab ({\bf f}')| ={\bf s}(B^{-1}({\bf f}'))!  \prod_{j=0}^{n}{\bf s}(B^j({\bf f'}))!
$$
This shows that the formula is also true for planar forests with height $(n+1)$. The inductive proof of the theorem is now completed.
\cqfd

\subsection{Wreath product representation}\label{wreathsec}
This section is concerned with an original wreath product interpretation
of the stabilizer of planar forests. To be self-contained, let us recall (see e.g. \cite{hoffman}) that the wreath product, written $\Ga_k<G>$ of an arbitrary group $G$ with the symmetric group of order $k$, $\Ga_k$, is the group that identifies, as a set, with $ \Ga_k\times G^k$, and where the product is defined by:
$$(s,g_1,...,g_k)\odot (s',g_1',...,g_k'):=(ss',g_1g'_{s^{-1}(1)},...,g_kg'_{s^{-1}(k)})$$
In addition, if a group $G$ acts on a set $S$, and if $\Ga_k$ acts by permutation on $S^k$, that is, if for any $s\in \Ga_k$ and any $(x_1,...,x_k)\in S^k$ the action of $s$ on $(x_1,...,x_k)$ is given by
$$s(x_1,...,x_k)=(x_{s^{-1}(1)},...,x_{s^{-1}(k)})$$
we then have, for any $(g_1,...,g_k)\in G^k$
$$s(g_1,...,g_k)s^{-1}(x_1,...,x_k)=s(g_1,...,g_k)(x_{s(1)},...,x_{s(k)})$$
$$=s(g_1x_{s(1)},...,g_kx_{s(k)})$$
$$=(g_{s^{-1}(1)}x_1,...,g_{s^{-1}(k)}x_k)$$
We are now in position to state and prove the main result of this section.
\begin{theo}
For any  planar forest $(T_1')^{m_1}...(T_k')^{m_k}\in{\cal PF}_{\bf p}$, with height 
$(n+1)$ written 
in terms of distinct planar trees $T_i'$, for some sequence of integers 
${\bf p}\in{\bf V}_{n+1}$, there is a natural isomorphism
$$Stab((T_1')^{m_1}...(T_k')^{m_k})\sim \left(\Ga_{m_1}<Stab(T_1')>\right)\vee \ldots\vee\left( \Ga_{m_k}<Stab(T_k')>\right)$$
so that, in particular, the stabilizer of any jungle in the equivalence class of $\bf f$ is conjugated in $\Ga_{\bf p}$ to $ \left(\Ga_{m_1}<Stab(T_1')>\right)\vee \ldots\vee\left( \Ga_{m_k}<Stab(T_k')>\right)$.
\end{theo}
\proof
We write $T_1''...T_{p_0}''$ for the expansion of the non commutative monomial $(T_{1}')^{m_1}...(T_{k}')^{m_k}$ as a product of planar trees (without exponents).
Let us also choose an arbitrary 
$$
{\bf s}=(s_0,...,s_{n+1})\in Stab({\bf f}')
$$ 
This implies in particular that the root permutation 
$s_0\in\Ga_{p_0}$ is of the form
$$
s_0=(s_0^{(1)}\vee\ldots\vee s_0^{(k)})\in\left(\Ga_{m_1}\vee\ldots\vee
\Ga_{m_k}\right)
$$
Also observe that $s_0$ acts on the planar forest $T_1''...T_{p_0}''$  by permutation of the planar trees:
$$s_0[T_1''...T_{p_0}'']:=(T_{s_0^{-1}(1)}''...T_{s_0^{-1}(p_0)}'')$$
This action by permutation can be viewed as a permutation (by blocks, according to the permutation of the planar trees) of the vertices of the planar forest, and this action can therefore be  rewritten uniquely under this process as the action of an element
$$
{\bf s}_0:=(s_0,s_{0,1},\ldots,s_{0,n+1})\in \Ga_{\bf p}
$$ on the vertices of the planar graph of $T_1''...T_{p_0}''$ 
$$
s_0[T_1''...T_{p_0}'']={\bf s}_0(T_1''...T_{p_0}'')
$$
Notice that ${\bf s}_0$ depends on the planar forest
$T_1''...T_{p_0}''$, since it acts non trivially on the vertices at levels $i\geq 1$. In general, the embedding $s_0\in\Ga_{p_0}\hookrightarrow{\bf s}_0\in \Ga_{\bf p}$
is different from the canonical embedding $s_0\in\Ga_{p_0}\hookrightarrow(s_0,1_{p_1},...,1_{p_{n+1}})\in \Ga_{\bf p}$.

We consider then the map (of sets)
$$
\begin{array}{rcl}
\mu~:~Stab({\bf f}')&\mapsto&
 \left(\Ga_{m_1}\vee\ldots\vee
\Ga_{m_k}\right)\times
\left(Stab({\bf f}')\cap \left(\{1_{p_0}\}\times\Ga_{(p_1,\ldots,p_{n+1})}\right)\right)\\
{\bf s}=(s_0,...,s_{n+1})&\mapsto&
\mu({\bf s})=(s_0,{\bf s}{\bf s}_0^{-1})
\end{array}
$$
In lemma~\ref{keyle}, we have proved that
$
\left(Stab(T_1''...T_{p_0}'')\cap \left(\{1_{p_0}\}\times\Ga_{(p_1,\ldots,p_{n+1})}\right)\right)$ is isomorphic to the permutation group $Stab(T_1'')\vee\ldots\vee Stab(T_{p_0}'')$. Therefore, 
the permutation sequence
$$
{\bf s}{\bf s}_0^{-1}=(1_{p_0},s_1s_{0,1}^{-1},\ldots,s_{n+1}s_{0,n+1}^{-1})
$$
can be alternatively represented as
$$
{\bf s}{\bf s}_0^{-1}=\left({\bf s}^{(1)}\vee\ldots\vee {\bf s}^{(p_0)}\right)\in Stab(T_1'')\vee\ldots\vee Stab(T_{p_0}'')
$$
Using these representations, we have that
\begin{eqnarray*}
\mu({\bf s})&=&\left((s_0^{(1)}\vee\ldots\vee s_0^{(k)}), \left({\bf s}^{(1)}\vee\ldots\vee {\bf s}^{(p_0)}\right)\right)\\
&=&
\left(s_0^{(1)},\left({\bf s}^{(1)}\vee\ldots\vee {\bf s}^{(m_1)}\right)\right)
\vee \left(s_0^{(2)},\left({\bf s}^{(m_1+1)}\vee\ldots\vee {\bf s}^{(m_1+m_2)}\right)\right)\\
&&\qquad \qquad\qquad \vee\ldots\vee
\left(s_0^{(k)},\left({\bf s}^{(m_1+\ldots+m_{k-1}+1)}\vee\ldots\vee {\bf s}^{(m_1+\ldots+m_{k-1}+m_k)}\right)\right)
\end{eqnarray*}
On the other hand, we have
\begin{eqnarray*}
\mu({\bf s}{\bf t})&=&(s_0t_0,{\bf st}({\bf s}_0{\bf t}_0)^{-1})\\
&=&
(s_0t_0,{\bf s}{\bf t}{\bf t}_0^{-1}{\bf s}_0^{-1})=
(s_0t_0,({\bf s}{\bf s}_0^{-1}) {\bf s}_0({\bf t}{\bf t}_0^{-1}){\bf s}_0^{-1})
\end{eqnarray*}
and for any level index $1\leq l\leq n$ we find that
\begin{eqnarray*}
s_{0,l}(t_lt_{0,l}^{-1})s_{0,l}^{-1}&=&
s_{0,l}({t}^{(1)}_l\vee\ldots\vee {t}^{(p_0)}_l)s_{0,l}^{-1}\\
&=&
({t}^{((s_0^{(1)})^{-1}(1))}_l\vee\ldots\vee 
{t}^{((s_0^{(1)})^{-1}(m_1))}_l)\\
&&\quad\vee ({t}^{((s_0^{(2)})^{-1}(m_1+1))}_l\vee\ldots\vee {t}^{((s_0^{(2)})^{-1}(m_1+m_2))}_l)\\
&&\quad\vee\ldots\vee({t}^{((s_0^{(k)})^{-1}(m_1+\ldots+m_{k-1}+1))}_l\vee\ldots\vee {t}^{((s_0^{(k)})^{-1}(m_1+\ldots+m_{k-1}+m_k))}_l)
\end{eqnarray*}
This implies that
$$
{\bf s}_0({\bf t}{\bf t}_0^{-1}){\bf s}_0^{-1}=(1_{p_0},r_1,\ldots,r_{n+1})
$$
with a sequence of mappings
$$
(1_{p_0},r_1,\ldots,r_{n+1})
=\left({\bf r}^{(1)}\vee\ldots\vee {\bf r}^{(p_0)}\right)\in Stab(T_1'')\vee\ldots\vee Stab(T_{p_0}'')
$$
given by
$$
\left\{
\begin{array}{rcl}
{\bf r}^{(1)}\vee\ldots\vee{\bf r}^{(m_1)}
&=&{\bf t}^{((s_0^{(1)})^{-1}(1))}\vee\ldots\vee {\bf t}^{((s_0^{(1)})^{-1}(m_1))}\\
{\bf r}^{(m_1+1)}\vee\ldots\vee{\bf r}^{(m_1+m_2)}&=&
 {\bf t}^{((s_0^{(2)})^{-1}(m_1+1))}\vee\ldots\vee {\bf t}^{((s_0^{(2)})^{-1}(m_1+m_2))}\\
 &\vdots&\\
{\bf r}^{(m_1+\ldots+m_{k-1}+1)}\vee\ldots\vee{\bf r}^{(m_1+\ldots+m_k)}&=&
 {\bf t}^{((s_0^{(k)})^{-1}(m_1+\ldots+m_{k-1}+1))}\vee\ldots\vee {\bf t}^{((s_0^{(k)})^{-1}(m_1+\ldots+m_k))} 
\end{array}
\right.
$$
This yields that
\begin{eqnarray*}
({\bf s}{\bf s}_0^{-1})\circ\left({\bf s}_0({\bf t}{\bf t}_0^{-1}){\bf s}_0^{-1}\right)
&=&
\left({\bf s}^{(1)}\vee\ldots\vee {\bf s}^{(p_0)}\right)\circ
\left({\bf r}^{(1)}\vee\ldots\vee {\bf r}^{(p_0)}\right)\\
&=&
\left(({\bf s}^{(1)}\circ {\bf r}^{(1)})\vee\ldots\vee ({\bf s}^{(p_0)}\circ {\bf r}^{(p_0)})\right)
\end{eqnarray*}
We finally obtain the following
formula
\begin{eqnarray*}
\mu({\bf s}{\bf t})&=&
\left(s_0^{(1)}\circ t_0^{(1)},\left(
( {\bf s}^{(1)}\circ {\bf t}^{ ( (s_0^{(1)})^{-1}(1) ) })
\vee\ldots\vee ( {\bf s}^{(m_1)}\circ {\bf t}^{ ( (s_0^{(1)})^{-1}(m_1) ) })\right)\right)\\
&&\vee\left(s_0^{(2)}\circ t_0^{(2)},\left(
( {\bf s}^{(m_1+1)}\circ {\bf t}^{ ( (s_0^{(2)})^{-1}(m_1+1) ) })
\vee\ldots\vee ( {\bf s}^{(m_1+m_2)}\circ {\bf t}^{ ( (s_0^{(2)})^{-1}(m_1+m_2) ) })\right)\right)\\
&&\vee\ldots\\
&=&\left\{\left(s_0^{(1)},\left({\bf s}^{(1)}\vee\ldots\vee {\bf s}^{(m_1)}\right)\right)\odot\left(t_0^{(1)},\left({\bf t}^{(1)}\vee\ldots\vee {\bf t}^{(m_1)}\right)\right)\right\}\\
&&\vee\left\{\left(s_0^{(2)},\left({\bf s}^{(m_1+1)}\vee\ldots\vee {\bf s}^{(m_1+m_2)}\right)\right)\odot\left(t_0^{(2)},\left({\bf t}^{(m_1+1)}\vee\ldots\vee {\bf t}^{(m_1+m_2)}\right)\right)\right\}
\\
&&\vee\ldots
\end{eqnarray*}
and the Theorem follows.
\cqfd

\subsection{Hilbert series method for forests enumeration}\label{hseries}
In the present subsection, we face the problem of computing other cardinals relevant to the analysis of Feynman-Kac mean field particle models, according to our Theorem~\ref{mainthm1}. For example, we want to be able to compute the number of forests in any ${\cal F}_{\bf p}$; the cardinals $|{\cal F}_{n,q}({\bf r})|$; or the number of forests in ${\cal F}_{n,q}$ with a given coalescence degree. Notice, by the way, that the notions of coalescence sequence and degree, as defined in the Introduction for sequences of maps in $\Aa_{n,q}$, and forests in 
${\cal F}_{n,q}$, go over to arbitrary forests and trees. 

 We let $\langle {\cal S}\rangle $ be the  commutative monoid generated by a subset of forests ${\cal S}\in{\cal F}$
 $$
 \langle {\cal S}\rangle=\left\{ f_1^{m_1}\ldots f^{m_k}_k~:~k\geq 0,~\quad \mbox{\rm and}\quad\forall 1\leq i\leq k\quad
 m_i\geq 0 \quad \mbox{\rm and}\quad f_i\in\Sa
 \right\}
 $$
 with the convention $f_1^{m_1}\ldots f^{m_k}_k=\emptyset$, for $k=0$. Notice that there is a canonical map from $\langle {\cal S}\rangle $ to $\cal F$, which is an embedding, e.g. if $\cal S$ is a set of (distinct) trees. We also consider a set of  
variables $\chi(T)$ indexed by the set
of trees $T\in \Ta$, and we let $M[\Ta]$ be the set of
monomials defined as follows 
 $$
M[{\cal T}]=\left\{ \chi(T_1)^{m_1}\ldots \chi(T_k)^{m_k}~:~k\geq 0,~\quad \mbox{\rm and}\quad\forall 1\leq i\leq k\quad
 m_i\geq 0 \quad \mbox{\rm and}\quad T_i\in\Ta
 \right\}
 $$
For any mapping $\chi $ on the set of forests
$$
\chi~:~{\bf f}\in\Fa\mapsto \chi({\bf f})\in M[\Ta]\quad\mbox{\rm such that}\quad
\forall ({\bf f},{\bf g})\in\Fa^2\quad\chi({\bf f}{\bf g})=\chi({\bf f})~\chi({\bf g})
$$
we have the Hilbert series expansions
\begin{equation}\label{hilberteq}
\frac{1}{1-\chi({\bf f})}=\sum_{{\bf g}\in  \langle {\bf f}\rangle }~\chi({\bf g})\quad\mbox{\rm and}\quad
\prod_{{\bf f}\in \Sa}\frac{1}{1-\chi({\bf f})}=
\sum_{{\bf g}\in  \langle {\cal S}\rangle }~\chi({\bf g})
\end{equation}

We let ${\bf V}_{n}^{\star}$ be the subset of all multi-indices ${\bf p}\in {\bf V}$, with height $ht({\bf p})\leq n$
$$
{\bf V}_{n}^{\star}:=\cup_{-1\leq k\leq n}{\bf V}_k
$$
\begin{defi}
The (multidegree) Hilbert series of forests, 
$\HH_{\cal F}^{n}({\bf x})$, is the Hilbert series associated to the partition of the set of forests of height less or equal than $n$ into subsets according to the number of vertices at each level, that is,
$$ \HH_{\cal F}^n({\bf x}):=\sum_{{\bf p}\in {\bf V}_n^{\star}}~|\Fa_{\bf p}|~
{\bf x}^{{\bf p}}=\HH_{\cal F}^{n-1}({\bf x})+\sum_{{\bf p}\in {\bf V}_n}~|\Fa_{\bf p}|~
{\bf x}^{{\bf p}}
$$
where we write ${\bf x}^{\bf p}$ as a shorthand for $x_0^{p_0}x_1^{p_1}\ldots$.
\end{defi}
We let  $\partial_{\bf p}=\partial_{x_0}^{p_0}\partial_{x_1}^{p_{1}}\ldots$, and we consider the mapping $B^{-1}$ from
${\bf V}$ into itself defined  for any $k\geq 0$ by
$$
B^{-1}~:~{\bf p}\in {\bf V}_{k}\mapsto 
B^{-1}({\bf p})=(1,{\bf p})=(1,p_0,p_1,\ldots)\in {\bf V}_{k+1}
$$
Since forests $\bf f$ of height $0$ are characterized by the number of roots in $\bf f$, we clearly have the formula
$$
\left(\forall {\bf p}\in {\bf V}_0\quad |\Fa_{\bf p}|=1\right)
\Longrightarrow
\HH_{\cal F}^0(x_0)=\frac{1}{1-x_0}=\sum_{p\geq 0}~x_0^p
$$

\begin{prop}\label{laprop}

$$
\forall n\geq 1\qquad 
\HH_{\Fa}^n({\bf x})=
           \prod_{{\bf p}\in {\bf V}^{\star}_{n-1}}\left(\frac{1}{1-{\bf x}^{B^{-1}({\bf p})}}\right)^{|\Fa_{\bf p}|}
           =
\HH_{\cal F}^{n-1}({\bf x})\times\prod_{{\bf p}\in {\bf V}_{n-1}}\left(\frac{1}{1-{\bf x}^{B^{-1}({\bf p})}}\right)^{\partial_{\bf p}\HH_{\cal F}^{n-1}({\bf 0})}
$$
In addition, the generating function associated with the
number of forests with multidegree ${\bf p}\in {\bf V}_n$,
and built with trees with height $n$   is given by  
$$\forall n\geq 1\qquad 
\sum_{{\bf p}\in {\bf V}_n}~|\langle\Ta(n)\rangle\cap \Fa_{\bf p}|~
{\bf x}^{{\bf p}}=\prod_{{\bf p}\in {\bf V}_{n-1}}\left(\frac{1}{1-{\bf x}^{B^{-1}({\bf p})}}\right)^{\partial_{\bf p}\HH_{\cal F}^{n-1}({\bf 0})}
$$
\end{prop}

Before getting into the proof of this result,
we notice that the formula stated in the above proposition makes the Hilbert series computable at any finite vertex order, and at any height
using any formal computation software. The first two orders can be handly computed. For $n=0$, we have already seen that
$$
\HH_{\cal F}^0({\bf x})=\frac{1}{1-x_0}=\sum_{p\geq 0}~x_0^p\Rightarrow\forall {\bf p}\in{\bf V}_0\quad
\partial_{\bf p}\HH_{\cal F}^{0}({\bf 0})=|\Fa_{\bf p}|=1
$$
By the recursion formula we find that
\begin{eqnarray*}
\HH_{\cal F}^1({\bf x})&=&\frac{1}{1-x_0}\times
\prod_{n\geq 1}\frac{1}{1-x_0x_1^{n}}\\
&=&\prod_{n\geq 0}\left(\sum_{k_n\geq 0}~(x_0~x_1^n)^{k_n}\right)=\sum_{(k_n)_{n\geq 0}\in\NN^{\NN}}~x_0^{\sum_{n\geq 0}k_n}~x_1^{\sum_{n\geq 0}nk_n}\\
&=&\sum_{p_0\geq 0}~
x_0^{p_0}+\sum_{{\bf p}\in{\bf V}_1}~|\Fa_{\bf p}|~
x_0^{p_0}x_1^{p_1}
\end{eqnarray*}
This readily yields that for any $ {\bf p}\in{\bf V}_1$, we have
$$
 |\Fa_{\bf p}|=\#\left\{
(k_n)_{n\geq 0}\in\NN^{\NN}~:~\sum_{n\geq 0}k_n=p_0\quad\mbox{\rm and}\quad
\sum_{n\geq 0}~n~k_n=p_1
\right\}
$$
More generally, the $(n+1)$-th order Hilbert series are given by the formula
\begin{eqnarray*}
\HH_{\cal F}^{n+1}({\bf x})&=&\prod_{{\bf p}\in {\bf V}^{\star}_{n}}\left(\frac{1}{1-{\bf x}^{B^{-1}({\bf p})}}\right)^{| \Fa_{\bf p}|}
\end{eqnarray*}
Recalling that for any $m\geq 1$, we have the following
formal series expansion
$$
\left(\frac{1}{1-u}\right)^m=\sum_{k\geq 0}\frac{(m-1+k)!}{(m-1)!~k!}~
u^k
$$
we find that
\begin{eqnarray*}
\HH_{\cal F}^{n+1}({\bf x})&=&\prod_{{\bf p}\in {\bf V}^{\star}_{n}}
\sum_{k_{\bf p}\geq 0}\frac{(| \Fa_{\bf p}|-1+k_{\bf p})!}{(| \Fa_{\bf p}|-1)!~k_{\bf p}!}~
{\bf x}^{k_{\bf p}~B^{-1}({\bf p})}\\
&=&\sum_{k\in\NN^{{\bf V}^{\star}_n}}
\left(\prod_{{\bf p}\in{\bf V}_n^{\star}}
\frac{(| \Fa_{\bf p}|-1+k({\bf p}))!}{(| \Fa_{\bf p}|-1)!~k({\bf p})!}\right)~
{\bf x}^{\sum_{{\bf p}\in{\bf V}_n^{\star}}k({\bf p})~B^{-1}({\bf p})}
\end{eqnarray*}
This implies that that for any $ {\bf p}\in{\bf V}_{n+1}$, we have
\begin{equation}\label{reccc}
 |\Fa_{\bf p}|=
 \sum_{\stackrel{k:{\bf q}\in{\bf V}^{\star}_n\mapsto k({\bf q})\in \NN}{
 \sum_{{\bf q}\in{\bf V}_n^{\star}}k({\bf q})~B^{-1}({\bf q}) ={\bf p}}
 }\left(\prod_{{\bf q}\in{\bf V}_n^{\star}}
\frac{(| \Fa_{\bf q}|-1+k({\bf q}))!}{(| \Fa_{\bf q}|-1)!~k({\bf q})!}\right)
\end{equation}
Now, we come to the proof of proposition~\ref{laprop}.\\
{\bf Proof of proposition~\ref{laprop}:}
We consider the set $\Sa(n)$ of all trees with height less or equal to $n$
\begin{eqnarray*}
\Sa(n)=\Ta^{\star}(n)&:=&\cup_{-1\leq k\leq n} \Ta(k)=\{\emptyset\}
\cup\left(
\cup_{{\bf p}\in {\bf V}_{n}^{\star}:p_0=1}\Ta_{\bf p}\right)\\
&=&\{\emptyset\}
\cup\left(
\cup_{0\leq k\leq n}
\cup_{{\bf p}\in {\bf V}_{k}:p_0=1}\Ta_{\bf p}\right)
\end{eqnarray*}
The set of forests generated by $\Sa(n)$ coincide with the set  $\Fa^{\star}(n)$ of
all the forests with height less or equal to $n$. That is, we have that
$$
\langle {\cal S}(n)\rangle=\Fa^{\star}(n):=\cup_{-1\leq k\leq n}\Fa(k)=
\cup_{-1\leq k\leq n}\cup_{{\bf p}\in {\bf V}_k}\Fa_{\bf p}
$$
To take the final step, we consider the mapping
$$
\forall {\bf p}\in {\bf V}_n\quad
\forall {\bf f}\in \Fa_{\bf p}\qquad
\chi({\bf f})={\bf x}^{v({\bf f})}=x_0^{p_0}x_1^{p_1}\ldots x_n^{p_n}
$$
From previous decompositions, and using the Hilbert series expansions (\ref{hilberteq}), we find that
\begin{eqnarray*}
\sum_{{\bf f}\in  \langle {\cal S}(n)\rangle }~{\bf x}^{v({\bf f})}&=&
\sum_{-1\leq k\leq n}\sum_{{\bf p}\in {\bf V}_k}~|\Fa_{\bf p}|~
{\bf x}^{{\bf p}}=\sum_{{\bf p}\in {\bf V}_n^{\star}}~|\Fa_{\bf p}|~
{\bf x}^{{\bf p}}
\end{eqnarray*}
Recalling that
$$
B~:~{\bf p}\in {\bf V}_{k}\rightarrow B({\bf p})\in {\bf V}_{k-1}\quad
\mbox{\rm and}\quad
|\Ta_{\bf p}|=|B(\Ta_{\bf p})|=|\Fa_{B({\bf p})}|
$$
we also have that
\begin{eqnarray*}
\prod_{{\bf f}\in \Sa(n)}
\frac{1}{1- {\bf x}^{ v({\bf f}) } 
           }&=&\frac{1}{1-x_0}~\prod_{{\bf p}\in {\bf V}_{n}^{\star}:p_0=1}
           \left(\frac{1}{1-{\bf x}^{{\bf p} }}\right)^{|\Ta_{\bf p}|}\\
&=&\frac{1}{1-x_0}~
           \prod_{0\leq k\leq n}\prod_{{\bf p}\in {\bf V}_k:p_0=1}\left(\frac{1}{1-{\bf x}^{{\bf p} }}\right)^{|\Ta_{\bf p}|}\\
           &=&\frac{1}{1-x_0}~
           \prod_{0\leq k\leq n}\prod_{{\bf p}\in {\bf V}_k:p_0=1}\left(\frac{1}{1-x_0~B({\bf x})^{B({\bf p}) }}\right)^{|B(\Ta_{\bf p})|}
           \end{eqnarray*}
           with the shift operator
           $
           B((x_k)_{k\geq 0})=(x_{k+1})_{k\geq 0}
           $.
This implies that           
\begin{eqnarray*}
\prod_{{\bf f}\in \Sa(n)}
\frac{1}{1- {\bf x}^{ v({\bf f}) } 
           }
           &=&
           \prod_{0\leq k\leq n}\prod_{{\bf p}\in {\bf V}_{k-1}}\left(\frac{1}{1-{\bf x}^{B^{-1}({\bf p})}}\right)^{|\Fa_{\bf p}|}
           =
           \prod_{{\bf p}\in {\bf V}^{\star}_{n-1}}\left(\frac{1}{1-{\bf x}^{B^{-1}({\bf p})}}\right)^{|\Fa_{\bf p}|}
           \end{eqnarray*}
           The end of the proof follows then from the Hilbert series formula (\ref{hilberteq}).
\cqfd

Recall now that the signed measures in the Laurent expansion of the distributions ${\QQ}_{n,q}^N$ are indexed by the coalescence degrees of forests. To compute or estimate these Laurent expansions, we are therefore interested in  a more precise Hilbert series, namely the one taking into account, besides the multidegrees of forests, their coalescence numbers.
\begin{defi}
We denote by $\Fa_{{\bf p}}[{\bf q}]:=\Fa_{\bf p}\cap c^{-1}({\bf q})
$ the set of forests in $\Fa_{{\bf p}}$ with a prescribed coalescence sequence ${\bf q}$ (we also use the convention $\Fa_{{\bf 0}}[{\bf 0}]=\{\emptyset\}$). 

We let
$\CC_{\cal F}^{n}({\bf x},{\bf y})$ 
be the (multidegree) Hilbert series of forests associated to the partition of the set of forests of height less or equal than $n$, 
and prescribed coalescence sequences  into subsets parametrized by multidegrees, that is
$$
\CC_{\cal F}^{n}({\bf x},{\bf y}):=\sum_{{\bf p}\in {\bf V}_n^{\star}}~
\sum_{{\bf q}\in C({\bf p}) }|\Fa_{\bf p}[{\bf q}]|~
{\bf x}^{{\bf p}}~{\bf y}^{{\bf q}}
$$
In the above display, $C({\bf p}) $ stands for
the set of coalescence multidegrees 
$$
\forall {\bf p}\in {\bf V}\qquad C({\bf p}):=c(\Fa_{\bf p})
=\left\{{\bf q}\in\NN^{\NN}~:~\mbox{\rm and}~~
{\bf q}\leq(B({\bf p})-1)_+\right\}
 $$
 with $(B({\bf p})-1)_+=((p_k-1)_+)_{k\geq 1}
$, for any
 $
 {\bf p}=(p_k)_{k\geq 0}
 \in{\bf V}$.
\end{defi}
\begin{prop}\label{laprop2}
The multidegree Hilbert series of coalescent forests
$\CC_{\Fa}^n({\bf x},{\bf y})$ satisfies the recursive formulae
\begin{eqnarray*}
\CC_{\cal F}^{n}({\bf x},{\bf y})&=&\CC_{\cal F}^{n-1}({\bf x},{\bf y})~
~\prod_{{\bf p}\in{\bf V}_{n-1}}\prod_{{\bf q}\in C(B^{-1}({\bf p}))}
\left(
\frac{1}{
1-{\bf x}^{ B^{-1}({\bf p}) }~
{\bf y}^{\bf q}
                                     }
                                     \right)^{
\left|\partial_{\bf p}\partial_{B({\bf q})}\CC_{\cal F}^{n-1}({\bf 0},{\bf 0})
\right|
}
\end{eqnarray*}
\end{prop}
Using the same lines of arguments as before, this proposition can
used to derive a recursive formula for the explicit combinatorial calculation of the
number of forests prescribed heights, and coalescence multi-indices. Notice that the first two orders are given by
$$
\CC_{\cal F}^{0}({\bf x},{\bf y})=\frac{1}{1-x_0}
\quad\mbox{\rm 
and}\quad
\CC_{\cal F}^{1}({\bf x},{\bf y})=
\prod_{p_0\geq 0}\frac{1}{(1-x_0x_1^{p_0}~y_0^{(p_0-1)_+})}
$$
Now, we come to the proof of proposition~\ref{laprop2}.\\
\proof
We use the same notation as the one used in the proof of proposition~\ref{laprop}. Firstly, we notice that
$
\Fa_{\bf p}=\cup_{{\bf q}\in C({\bf p})}\Fa_{\bf p}[{\bf q}]
$,
and
$ht({\bf q})\leq \left(ht({\bf p})-1\right)$, for any $ {\bf q}\in C({\bf p})$, and where
 $
ht({\bf q}):=\inf{\left\{
 n\geq 0:\sup_{m>n}{q_m}=0
 \right\}}
 $, for any ${\bf q}=(q_m)_{m\geq 0}$.
Also observe that
 $$
\Fa_{B^{-1}({\bf p})}=\Ta_{(1,{\bf p})}\Longrightarrow  
C(B^{-1}({\bf p}))=c(\Ta_{(1,{\bf p})})=\{p_0-1\}\times C({\bf p})
 $$
We next consider the mapping $\chi$  
$$
\forall {\bf p}\in {\bf V}_n\quad
\forall {\bf f}\in \Fa_{\bf p}[{\bf q}]\qquad
\chi({\bf f})={\bf x}^{v({\bf f})}~{\bf y}^{c({\bf f})}=\left(x_0^{p_0}x_1^{p_1}\ldots x_n^{p_n}\right)~\left(y_0^{q_0}y_1^{q_1}\ldots y_n^{q_{n-1}}
\right)
$$
Using the Hilbert series expansions (\ref{hilberteq}), we readily check that
$$
\prod_{{\bf f}\in {\cal S}(n)}\frac{1}{1-{\bf x}^{v({\bf f})}~{\bf y}^{c({\bf f})}}=\prod_{{\bf f}\in {\cal S}(n)}\left( \sum_{{\bf g}\in  \langle {\bf f}\rangle }~{\bf x}^{v({\bf g})}~{\bf y}^{c({\bf g})}\right) =
\sum_{{\bf g}\in  \langle {\cal S}(n)\rangle }~{\bf x}^{v({\bf g})}~{\bf y}^{c({\bf g})}
$$
Arguing as in the proof of proposition~\ref{laprop}, we find that
\begin{eqnarray*}
\sum_{{\bf g}\in  \langle {\cal S}(n)\rangle }~{\bf x}^{v({\bf g})}~{\bf y}^{c({\bf g})}&=&1+\sum_{{\bf p}\in{\bf V}_0}|\Fa_{\bf p}|~{\bf x}^{{\bf p}}+\sum_{k=1}^n\sum_{{\bf p}\in{\bf V}_k}\sum_{{\bf q}\in C({\bf p})}~\left|\Fa_{\bf p}[{\bf q}]\right|~
{\bf x}^{\bf p}~{\bf y}^{\bf q}\\
&=&\sum_{{\bf p}\in{\bf V}_n^{\star}}\sum_{{\bf q}\in C({\bf p})}~\left|\Fa_{\bf p}[{\bf q}]\right|~
{\bf x}^{\bf p}~{\bf y}^{\bf q}
\end{eqnarray*}
In much the same way, we have that
\begin{eqnarray*}
\prod_{{\bf f}\in {\cal S}(n)}\frac{1}{1-{\bf x}^{v({\bf f})}~{\bf y}^{c({\bf f})}}
&=&\frac{1}{1-x_0}~\prod_{k=1}^n\prod_{{\bf p}\in{\bf V}_k:p_0=1}\prod_{{\bf q}\in C({\bf p})}\left(\frac{1}{1-{\bf x}^{\bf p}~{\bf y}^{\bf q}}\right)^{\left|\Fa_{\bf p}[{\bf q}]\right|}
\end{eqnarray*}
Therefore,
$$
\begin{array}{l}
\prod_{{\bf f}\in {\cal S}(n)}\frac{1}{1-{\bf x}^{v({\bf f})}~{\bf y}^{c({\bf f})}}
\\
=\prod_{k=0}^n\prod_{{\bf p}\in{\bf V}_k:p_0=1}\prod_{{\bf q}\in C({\bf p})}
\left(
\frac{1}{
1-
x_0
B({\bf x})^{ B({\bf p}) }~
y_0^{(p_1-1)}~ B({\bf y})^{
B({\bf q})}
                                     }
                                     \right)^{
\left|
\Fa_{(1,
           B({\bf p}))
       }  [(p_1-1,B({\bf q}))]
\right|
}
\end{array}
$$
This yields that
\begin{eqnarray*}
\prod_{{\bf f}\in {\cal S}(n)}\frac{1}{1-{\bf x}^{v({\bf f})}~{\bf y}^{c({\bf f})}}
&=&\prod_{k=0}^n\prod_{{\bf p}\in{\bf V}_{k-1}}\prod_{{\bf q}\in C(B^{-1}({\bf p}))}
\left(
\frac{1}{
1-{\bf x}^{ B^{-1}({\bf p}) }~
{\bf y}^{\bf q}
                                     }
                                     \right)^{
\left|
\Fa_{
           B^{-1}({\bf p})
       }  [{\bf q}]
\right|
}\\
&=&\prod_{{\bf p}\in{\bf V}_{n-1}^{\star}}\prod_{{\bf q}\in C(B^{-1}({\bf p}))}
\left(
\frac{1}{
1-{\bf x}^{ B^{-1}({\bf p}) }~
{\bf y}^{\bf q}
                                     }
                                     \right)^{
\left|
\Fa_{
           {\bf p}
       }  [B({\bf q})]
\right|
}
\end{eqnarray*}
The proposition follows.
\cqfd

Notice that the Hilbert series technique can be developed to any order of refinement. For example, it could be used to take into account, besides the number of vertices or of coalescences at each level, the cardinals $\# ({\bf f})$, so that the coefficients of the Laurent expansion of the measures ${\QQ}_{n,q}^N$ could be read, in the end, entirely on the corresponding Hilbert series. We leave the task of expressing the recursive formula to the interested reader, and simply point out that the technique allows an easy, systematic, recursive, computation of the coefficients of the expansion of the measures ${\QQ}_{n,q}^N$ at any order, both in $n$ and $q$. The observation can be usefull, especially in view of the systematic development of numerical schemes and numerical approximations based on Feynman-Kac particle models.

\subsection{Some forests expansions}\label{someforests}

In the present section, and the forthcoming one, we take advantage of the langage of trees and of the results obtained on their statistics to compute the first orders of the Laurent functional representation of ${\QQ}_{n,q}^N$ and, respectively, to derive a natural generalization to Feynman-Kac particle models of the classical Wick product formula. 

Let us fix $n$ and $q\geq 4$, so that the notation $\bf q$ denotes, once again, the constant sequence of lenght $n$ associated to $q$. 
As we have already pointed out,
there is only one forest in ${\cal F}_{\bf q}$ without coalescence, which is the product of $q$ trivial trees of height $n$. There is only one forest in ${\cal F}_{\bf q}$ with only one coalescence, at level $i$, that will be written ${\bf f}_{1,i}$. Its only non trivial tree is the tree with one coalescence at level $i$ and two leaves at level $n$. 
There are two forests with coalescence degree 2 and the two coalescences at levels $i$, written ${\bf f}_{2,i}^1$ and ${\bf f}_{2,i}^2$. The notation ${\bf f}_{2,i}^1$ denotes the forest with only one non trivial tree with coalescence degree 2, a vertex with coalescence order 2 at level $i$ and its three leaves at level $n$. The notation ${\bf f}_{2,i}^2$ denotes the forest with two non trivial trees with coalescence degree 1 and the coalescence at level $i$. There are four forests with coalescence degree 2 and the two coalescences at levels $i<j$, written ${\bf f}_{2,i,j}^{k}$, $k=1...4$.
The forest ${\bf f}_{2,i,j}^{1}$ has one non trivial tree with coalescences at levels $i$ and $j$ and its tree leaves at level $n$. The forest ${\bf f}_{2,i,j}^{2}$ has one non trivial tree with coalescences at levels $i$ and $j$ and its three leaves at the levels $j,n,n$. The forest ${\bf f}_{2,i,j}^{3}$ has two non trivial trees with one coalescence at level $i$, resp. $j$ and their two leaves at level $n$. The forest ${\bf f}_{2,i,j}^{4}$ has two non trivial trees with one coalescence at level $i$, resp. $j$ and their two respective leaves at level $j,n$, resp. $n,n$. 

Expanding the formulae for $\partial^i{\QQ}_{q,n}$ and using the formulae obtained in Thm \ref{numberofjungles} for the cardinals $\# ({\bf f})$ and using that $s(q,q-2)= \left(\begin{array}{c}
{q}\\{3}\end{array}\right)~\frac{3q-1}{4}$ (see for instance \cite{comtet}, on  page 63), we get the following result.

\begin{cor}\label{3order}
The first three order terms in the polynomial
functional representation of ${\QQ}_{n,q}^N$ (\ref{polyf}) are given by the following
formulae
\begin{eqnarray*}
\partial^0  {\QQ}_{n,q}&=&\gamma_n^{\otimes q}\\
\partial^1  {\QQ}_{n,q}&=&\frac{q(q-1)}{2}
\sum_{0\leq k\leq n} \left(
 \Delta_{n,q}^{{\bf f}_{1,k}}-\gamma^{\otimes q}_n\right) \\
\partial^2  {\QQ}_{n,q}&=&
\frac{q!}{(q-3)!~3!}~\sum_{0\leq k\leq n}
 \left(\Delta_{n,q}^{{\bf f}_{2,k}^1} 
+ \frac{3}{4}~(q-3)~\Delta_{n,q}^{{\bf f}_{2,k}^2} 
-
 \frac{3}{2}~(q-1)
 ~\Delta_{n,q}^{{\bf f}_{1,k}}
+
 \frac{(3q-1)}{4}~\gamma^{\otimes q}_n
 \right)\\
 &&+
\left( \frac{q(q-1)}{2}\right)^2~\sum_{0\leq k<l\leq n} 
\left(
\gamma^{\otimes q}_n-\left( \Delta_{n,q}^{{\bf f}_{1,l}}+ \Delta_{n,q}^{{\bf f}_{1,k}}\right)\right)\\
 &&
+ \frac{q(q-1)}{2}~\sum_{0\leq k<l\leq n} 
\left(
\Delta_{n,q}^{{\bf f}_{2,k,l}^2}   
+(q-2)~(
\Delta_{n,q}^{{\bf f}_{2,k,l}^1}+\Delta_{n,q}^{{\bf f}_{2,k,l}^4}    )+\frac{(q-2)(q-3)}{2}
\Delta_{n,q}^{{\bf f}_{2,k,l}^3}
\right)
\end{eqnarray*}
\end{cor}

\subsection{A Wick product formula on forests}\label{wicksetcion}

Let now $\Ba_0^{  sym}(E^q_n)\subset\Ba_b(E^q_n)$ be the set of symmetric functions $F$ on $E^q_n$ such that 
$$(D_{1_{q-1}}\otimes \gamma_n) (F)(x_1,...,x_{q-1})=
\int~F(x_1,\ldots,x_{q-1},x_n)~\gamma_n(dx_n)=0
$$

Notice that $\Ba_0^{sym}(E^q_n)$ contains
the set of functions $F=(P)_{\rm\tiny sym}$, with $P\in \mbox{\rm Poly}(E^q_n)$,  where 
 $\mbox{\rm Poly}(E^q_n)\subset \Ba_b(E^q_n)$ stands for the 
subset of polynomial functions of the form
$$
P=\sum_{{ a}\in I}~c({ a})~f^{ a}\quad\mbox{\rm with}\quad f^{ a}=
 (f^{{a}(1)}\otimes\ldots\otimes f^{{a}(q)})
$$
In the above display, $I$ is a 
 finite subset of $\NN^{[q]}$, $c\in \RR^{I}$,
and the elementary functions $f^{p}\in\Ba_b(E_n)$ are choosen  such
that
$\gamma_n(f^{p})=0
$.
For instance, we can take
$$
f^{p}=(g^{p}-\eta_n(g^{p}))\quad\mbox{\rm with}\quad 
g^{p}\in\Ba_b(E_n)
$$ 

Let us recall that, when explicited in the particle models setting, the traditional Wick product formula reads, for any integer $q$ and any $F\in{\cal B}_0^{sym}(E_0^q)$
$$\partial^{i}{\QQ}_{0,q}(F)=0,\ i<\frac{q}{2}$$
and, if $q$ is even,
$$\partial^{\frac{q}{2}} {\QQ}_{0,q} (F)=
\frac{q!}{2^{q/2}~(q/2)!}~\Delta^{\bf f}_{0,q} (F)
$$
where $\bf f$ is the forest in ${\cal F}_{0,q}$ containing $\frac{q}{2}$ copies of the tree of unit height, with two vertices at level 1, and $\frac{q}{2}$ copies of the tree with the root as unique vertex.
For symmetric tensor product functions  
$F=(f^1\otimes\ldots\otimes f^q)_{\rm\tiny sym}
$, 
associated with a collection of functions $f^i\in {\cal B}_b(E_0)$, such that $\eta_0(f^i)=\gamma_0(f^i)=0$, for any $1\leq i\leq q$, we readily find that
$$
\Delta^{\bf f}_{0,q} (F)= \frac{2^{q/2}~(q/2)!}{q!}
\sum_{{\cal I}_q}~\left(\prod_{\{i,j\}\in {\cal I}_q }\eta_0(f^i f^j)\right)
$$
as soon as $q$ is even. In the above displayed formula, ${\cal I}_q$ ranges over all partitions of $[q]$ into pairs. In a more probabilistic language,
the above formula can be interpreted as the $q$-th order central moment
$$
\partial^{\frac{q}{2}} {\QQ}_{0,q} (F)=\EE(W_0(f^1)\ldots W_0(f^q))=
\sum_{{\cal I}_q}~\left(\prod_{\{i,j\}\in {\cal I}_q }\EE(W_0(f^i)W_0(f^j))\right)
$$
of a Gaussian field  $W_0$ on the Banach space of functions $\Ba_b(E_0)$, such that for any pair of functions $(\varphi,\psi)\in \Ba_b(E_0)$,
$
\EE(W_0(\varphi))=0$, and $\EE(W_0(\varphi)W_0(\psi))=\eta_0(\varphi\psi)
$. In terms of the Laplace moment generating function
$$
{\bf x}=(x^i)_{1\leq i\leq q}\mapsto
\exp{\left(\frac{1}{2} {\bf x}^{\prime}C_0(f){\bf x}\right)}=\frac{1}{((2\pi)^q|C_0(f)|)^{1/2}} \int_{\RR^q} d{\bf y}
\exp{\left(-\frac{1}{2} {\bf y}^{\prime}C_0(f)^{-1}{\bf y}+{\bf y}^{\prime}{\bf x}\right)}
$$
with the gaussian covariance matrix 
$C_0(f)=(C_0(f^i, f^j))_{1\leq i,j\leq q}$, we also have that
 $$
\partial^{\frac{q}{2}} {\QQ}_{0,q} (F)=\frac{\partial ^q}{\partial {\bf x}^q}\exp{\left(\frac{1}{2} {\bf x}^{\prime}C_0(f){\bf x}\right)} _{\left|{\bf x}={\bf 0}\right.}\quad
\mbox{\rm 
with} \quad \frac{\partial ^q}{\partial {\bf x}^q}:=\frac{\partial}{\partial {x}^1 }\ldots\frac{\partial}{\partial x^q}
$$

In the present section, we will show that the Wick formula generalizes to forests in ${\cal F}_{n,q}$ of arbitrary height.
Let us start by listing various straightforward properties of trees and forests.
A tree $T$ with coalescence degree $d$ has exactly $(d+1)$ leaves. A forest with coalescence degree $d$ has at most $d$ non trivial trees, and the equality holds if and only if all its non trivial trees have coalescence degree $1$. In particular, if a forest in ${\cal F}_{n,q}$ has coalescence degree $d$, it has at most $(2d)$ leaves belonging to non trivial trees so that, if $d<\frac{q}{2}$, there is at least one vertex at level $n+1$ belonging to a trivial tree (that is, a tree with coalescence number 0).

The same reasoning shows that, when $d=\frac{q}{2}$, a forest in ${\cal F}_{n, q}$ with coalescence degree $d$, and coalescence sequence $\bf r$, does not contain a trivial tree of height $(n+1)$ if and only if it is the forest ${\bf f}_{\bf r}:=T_0^{r_0}U_0^{r_0}...T_{n}^{r_{n}}U_{n}^{r_{n}}$, where we write $T_k$ for the unique tree of coalescence degree 1 with a coalescence at level $k$ and its two leaves at level $(n+1)$, and where we write $U_k$ for the trivial tree of height $k$.\label{discusswick}

We conclude this series of remarks by noting that, if ${\bf f}\in {\cal F}_{n,q}$ can be written as the product (or disjoint union) of a forest $\bf g$ in ${\cal F}_{n,q-1}$ with $U_{n+1}$, the trivial tree of height $(n+1)$, then, for any $F\in \Ba_b^{sym}(E^q_n)$ we have, by definition of the measures $\Delta^{\bf f}$
$$\Delta^{\bf f}_{n,q}(F)=\Delta^{\bf g}_{n,q}(D_{1_{q-1}}\otimes \gamma_n)(F)$$
which is equal to 0 if $F\in \Ba_0^{sym}(E^q_n)$. 

We are now in position to derive the forest Wick formula.

\begin{theo}\label{thwick}
For any even integer $q\leq N$ and  any symmetric function
$F\in \Ba_0^{ sym}(E^q_n)$, we have 
\begin{equation}\label{wick1}
\forall k<q/2\qquad\partial^k  {\QQ}_{n,q}
(F)=0\quad\mbox{\rm and}
\quad
\partial^{q/2} {\QQ}_{n,q} (F)=\sum_{{\bf r}< {\bf q},|{\bf r}|=\frac{q}{2}}~
\frac{q!}{2^{q/2}~{\bf r}!}~\Delta^{{\bf f}_{\bf r}}_{n,q} F
\end{equation}
For odd integers $q\leq N$, the partial measure valued
derivatives $\partial^k$ are the null measure on $\Ba_0^{sym}(E^q_n)$, up to any order $k\leq \lfloor q/2\rfloor$.
\end{theo}

We close this section with a gaussian field interpretation of the 
Wick formula (\ref{wick1}). We further assume that $q$ is an even integer. We consider a collection of independent
gaussian fields $(W_k)_{0\leq k\leq n}$ on the Banach spaces
$(\Ba_b(E_k))_{0\leq k\leq n}$, with for any $(\varphi_k,\psi_k)\in\Ba_b(E_k)^2 $, and $0\leq k\leq n$
$$
\EE(W_k(\varphi_k))=0\quad\mbox{\rm and}\quad\EE(W_k(\varphi_k)W_k(\psi_k))=\gamma_k(\varphi_k\psi_k)
$$
We also introduce the centered gaussian field $V_n$ on $\Ba_b(E_n)$
defined for any $\varphi_n\in\Ba_b(E_n)$ by the following formula 
$$
V_{n}(\varphi_n)=\sum_{0\leq k\leq n}~\sqrt{\gamma_k(1)}~W_k(Q_{k,n}(\varphi_n))
$$

Let $(\varphi_i)_{1\leq i\leq q}\in {\cal B}_b(E_n)^q$ be a collection
of functions such that $\gamma_n(\varphi_i)=0$, for any $1\leq i\leq q$.
For the tensor product function  $F=\frac{1}{q!}\sum_{\sigma\in \Ga_q}(\varphi_{\sigma(1)}\otimes\ldots \otimes \varphi_{\sigma(q)})$, one can check that
\begin{eqnarray*}
\Delta^{{\bf f}_{\bf r}}_{n,q}(F)&=&\frac{2^{q/2}{\bf r}!}{q!}\sum_{I\in \Ia}\prod_{0\leq k\leq n}
\sum_{J_k\in \Ia_k}
\left\{\gamma_k(1)^{r_k}\prod_{(i,j)\in J_k}\gamma_k(
Q_{k,n}(\varphi_i)
Q_{k,n}(\varphi_{j}))
\right\}
\end{eqnarray*}
In the above displayed formula, the first sum is over the set $\Ia$ of all partitions  $I=(I_k)_{0\leq k\leq n}$ of $[q]$ into
$(n+1)$ blocks with cardinality
$|I_k|=(2r_k)$, and the second sum ranges over the set $\Ia_k$ of all partitions $J_k$ of the sets $I_k$ into pairs, with $0\leq k\leq n$. By definition of the gaussian fields $(W_k)_{0\leq k\leq n}$, and due to the classical Wick formula, we find that
\begin{eqnarray*}
\frac{q!}{2^{q/2}{\bf r}!}\Delta^{{\bf f}_{\bf r}}_{n,q}(F_n)&=&\sum_{I\in \Ia}\prod_{0\leq k\leq n}
\EE\left(\prod_{i\in I_k}\sqrt{\gamma_k(1)}~W_k(Q_{k,n}(\varphi_i))\right)\\
&=&\EE\left(\sum_{I\in \Ia}\prod_{0\leq k\leq n}
\prod_{i\in I_k}\sqrt{\gamma_k(1)}~W_k(Q_{k,n}(\varphi_i))\right)
\end{eqnarray*}
from which we arrive at the following formula
$$
\sum_{{\bf r}< {\bf q},|{\bf r}|=\frac{q}{2}}~
\frac{q!}{2^{q/2}~{\bf r}!}~\Delta^{{\bf f}_{\bf r}}_{n,q}(F_n)
=
\EE\left( \sum_{{\bf r}:|2{\bf r}|=q} \sum_{I}\prod_{0\leq k\leq n}
\prod_{i\in I_k}\sqrt{\gamma_k(1)}~W_k(Q_{k,n}(\varphi^i_n))\right)
$$
Recalling that all gaussian fields $(W_k)_{0\leq k\leq n}$ are independent and centered, we prove that
$$
\partial^{q/2} {\QQ}_{n,q} (F_n)=\EE\left(\prod_{1\leq i\leq q}
\left(\sum_{0\leq k\leq n} \sqrt{\gamma_k(1)}~W_k(Q_{k,n}(\varphi^i_n))\right)\right)=\EE\left(\prod_{1\leq i\leq q}V_{n}(\varphi_n^i)\right)
$$
Written in a more synthetic way, we have proved the following formula
$$
\partial^{q/2} {\QQ}_{n,q} (F_n)=\EE\left(V_n^{\otimes q}(F_n)\right)
$$
This result can alternatively be derived combining the $\LL_q$-mean error
estimates presented in \cite[Thm.7.4.2]{fk}, with
the multidimensional
central limit theorems presented in~\cite[Prop.9.4.1]{fk}. More precisely, the $q$-dimensional particle random fields  $\left(V_{n}^N(\varphi_n^i)\right)_{1\leq i\leq q}:=\left(\sqrt{N}\gamma^N_n(\varphi^i_n)
\right)_{1\leq i\leq q}$, converge in law to 
$(V_{n}(\varphi_n^i))_{1\leq i\leq q}$. By the continuous mapping theorem, combined with simple uniform integrability arguments, one checks that 
$$
N^{q/2}~{\QQ}_{n,q}^N (F_n)=
\EE\left((V_n^N)^{\otimes q}(F_n)\right)
\stackrel{N\uparrow\infty}{\longrightarrow}\partial^{q/2} {\QQ}_{n,q} (F_n)=\EE\left(V_n^{\otimes q}(F_n)\right)
$$

\section{Extension to path-space models}\label{extensions}
In the present section, we extend our previous analysis to the statistical study of path spaces. Due to the mean-field definition of the trajectories, it will appear soon that distributions on the path space are parametrized by 
forests with a slightly more complex structure than the classical ones. Namely, we have to introduce colored trees and forests (black and white, with a particular structure deviced to reflect the geometry of paths). These new objects are introduced and studied in the first two subsections; applications to the path space and to propagation of chaos properties are postponed to the last three sections.  

\subsection{Colored trees, forests, and jungles}\label{tropsect}

\begin{defi}
A colored tree is a tree with colored vertices, with two distinguished colors, say black and white. Only black vertices may have an 
ingoing edge. That is, equivalently, all white vertices are leaves -notice that the converse is not true in general.
A colored forest is a multiset of colored trees.
The sets of colored trees and colored forests are denoted respectively by $\overline{\cal T}$ and $\overline{\cal F}$.
\end{defi}

\ \par

\xymatrix{&\circ\ar[dr]&\bullet\ar[d]&\bullet\ar[dl]&&\circ\ar[d]\\
\circ\ar[dr]&&\bullet\ar[dl]&\circ\ar[dr]&\circ\ar[d]&\bullet\ar[dl]\\
&\bullet&&&\bullet&&}

\ \par

{\it Figure 3}: A colored forest

\ \par

Notice that a colored forest can also be viewed as a commutative monomial over the set of colored trees.
In this interpretation, the generating series techniques that we have developped to deal with the enumeration of forests will apply to colored forests. The computation of the corresponding series is left to the interested reader. 

Most of the notions associated to trees and forests go over in a straightforward way to colored forests and colored trees. 
As a general rule, we will write a line over symbols associated to colored trees and colored forests.   For instance, we write 
$\overline{v}_k(\overline{T})=({v}_k(\overline{T}),{v}_k'(\overline{T}))$ for the number ${v}_k(\overline{T})$, resp. ${v}_k'(\overline{T})$, of white, resp. black vertices in the colored tree
$\overline{T}$, at each level $k\geq 0$.

We also let $\overline{\bf V}$ be the set of all sequences of pair of integers $\overline{\bf p}=(\overline{p}_k)_{k\geq 0}\in (\NN^2)^{\NN}$, with $\overline{p}_k=({p}_k,{p}_k')$ for every $k\geq 0$,
 satisfying the following property
$$
\exists \overline{ht}(\overline{\bf p})\in\NN\quad
\mbox{\rm s.t.}\quad \inf_{k\leq  \overline{ht}(\overline{\bf p})}{p_k+p'_k}>0\quad\mbox{\rm and}\quad
\sup_{k>  \overline{ht}(\overline{\bf p}) }{p_k+p'_k}=0
$$
For any $n\in\NN$, we denote by $\overline{\bf V}_n\subset\overline{\bf V}$
the subset of sequences $\overline{\bf p}$ such that the height of $\overline{\bf p}$, $\overline{ht}(\overline{\bf p})$, is equal to $n$. Finally, 
we let  ${\cal \overline{ T}}_{\overline{\bf p}}$ be the set of colored trees $\overline T$ with
$v_k(\overline T)=p_k$ white vertices, and $v_k'(\overline T)=p_k'$ black vertices, at each level $k\geq 0$. Since $\overline T$ is a colored tree, this implies that $p_0=0$ and $p_0'=1$, excepted if $\overline{ht}(\overline{\bf p})=0$. In that case, $p_0$ may also be equal to 1 (and then $p_0'=0$).
Notice that the
set of black and white vertices at each level $k$ is in bijection with the disjoint
union $$
[\overline{p}_k]=[(p_k,p'_k)]:=[p_k]\scoprod [p'_k] 
$$
of the sets $[p_k]$ and $[p'_k]$.

For $\overline{\bf p}\in \overline{\bf V}_{n+1}$, let us call, by analogy with the uncolored case, planar colored trees $\overline{T}$ of type $\overline{\bf p}$ any sequence
of maps $(\overline{\alpha}_0,\ldots,\overline{\alpha}_n)$,
where
$$
\overline{\alpha}_k=(\alpha_k\tcoprod\alpha'_k)\in 
[p'_k]^{[p_{k+1}]}\tcoprod[p'_k]^{[p'_{k+1}]} =:[p'_k]^{[(p_{k+1},p'_{k+1})]}
$$
is the disjoint union of the pair of weakly increasing maps 
$\alpha_k\in [p'_k]^{[p_{k+1}]}$, and $\alpha_k\in [p'_k]^{[p_{k+1}']}$.

In the above definition, the parameter level $k$ runs from
$0$, to $n$, and we have used the conventions
$
(\emptyset\tcoprod\emptyset)=\emptyset$, and $(\alpha\tcoprod\emptyset)=\alpha$.
Writing ${\cal P{\overline T}}$ for the set of planar colored trees and ${\cal P\overline{ T}}_{\overline{\bf p}}$ for the set of planar colored trees of type $\overline{\bf p}$, from the previous discussion, for any $\overline{\bf p}\in \overline{\bf V}_{n+1}$ with $p_0=0, ~p_0'=1$, we have that
$$
\begin{array}{l}
{\cal P\overline{ T}}_{\overline{\bf p}}\\
\simeq 
\left\{
(\overline{\alpha}_0,\ldots,\overline{\alpha}_n)\in [1]^{[\overline{p}_1]}\times [p'_1]^{[\overline{p}_2]}\times \ldots\times[p'_n]^{[\overline{p}_{n+1}]}~:~\overline{\alpha}_k=(\alpha_k \tcoprod \alpha'_k)~\mbox{\rm with}~\alpha_k\uparrow~\mbox{\rm and}~\alpha_k'\uparrow\right\}
\end{array}
$$

A planar colored forest is, by definition, a sequence (or a non commutative monomial) of planar colored trees. 

Since to any sequence $(\overline{\alpha}_k^1,...,\overline{\alpha}_k^n)$, $\overline{\alpha}_k^i=(\alpha_k^i\coprod \alpha_k^{i'})\in [p_k^{i'}]^{[(p_{k+1}^i,p_{k+1}^{i'})]}$ is associated a unique element 
$(\alpha_k^1\coprod ...\coprod\alpha_k^n;\alpha_k^{1'}\coprod ...\coprod\alpha_k^{n'})$ of
$[p_k^{1'}]\coprod ...\coprod [p_k^{n'}]^{(\coprod\limits_{i=1}^n[p_{k+1}^i])\coprod (\coprod\limits_{i=1}^n[p_{k+1}^{i'}])}$,
planar colored forests are in bijection with
sequences of disjoint union of weakly increasing maps, in the sense that (with self-explaining notations)
$$
{\cal P\overline{ F}}_{\overline{\bf p}}\simeq {\cal I\overline{A}}_{\overline{\bf p}}
:=\left\{
(\overline{a}_0,\ldots,\overline{a}_n)
\in \overline{\cal A}_{\overline{\bf p}}~:~
\overline{a}_k=(a_k\tcoprod a'_k)~\mbox{\rm with}~a_k\uparrow~\mbox{\rm and}~a_k'\uparrow
\right\}
$$
for any $\overline{\bf p}\in \overline{\bf V}_{n+1}$, with
$$
 \overline{\cal A}_{\overline{\bf p}}:=\left(
[p'_0]^{[\overline{p}_1]}\times [p'_1]^{[\overline{p}_2]}\times\ldots\times[p'_n]^{[\overline{p}_{n+1}]}\right)
$$

\begin{defi}

We associate with a sequence of mappings 
$
{\overline{\bf a}}=(\overline{a}_0,\ldots,\overline{a}_n)\in \overline{\cal A}_{\overline{\bf p}}$, with $\overline{\bf p}\in \overline{\bf V}_{n+1}$, a  colored forest $\overline{F}({\overline{a}})$ of height $(n+1)$. The vertex set of the colored forest $\overline{F}({\overline{a}})$ is given by the
disjoint union set  $\left( \coprod_{k=0}^{n+1}[p_k]\right) \coprod \left( \coprod_{k=0}^{n+1}[p_k']\right)$.
The edges correspond to the pairs $(i,a_k(i))$, and $(i',a_k'(i'))$, with the pair of white and black colored vertices   
$(i,i')\in \left([p_{k+1}]\coprod [p_{k+1}']\right)$. 
The sequence can be represented graphically by an entangled planar graph with colored vertices, or colored jungle , written $\overline{J}({\overline{a}})$.
Notice that
edges between two levels are allowed to cross. At a given level $k$, all white vertices are drawn left to the black vertices. The set of colored jungles is written $\overline{\cal J}$.
\end{defi}

\ \par

\xymatrix{&&\circ\ar[dr]&\bullet\ar[dl]&\bullet\ar[d]\\
\circ\ar[dr]&\circ\ar[drr]&\bullet\ar[dl]&\bullet\ar[d]&\bullet\ar[dlll]\\
&\bullet&&\bullet&&}

\ \par
{\it Figure 3}: A colored jungle

\ \par

Once again, the notation we use on colored forests extends to colored jungles in a self-explanatory
way. For instance, $\overline{\cal J}_{\overline{\bf p}}$,  $\overline{\bf p}\in \overline{\bf V}_{n+1}$, stands for the set of colored jungles $\overline{\bf j}$ with
$v_k(\overline{\bf j})=p_k$ white vertices, 
and $v_k'(\overline{\bf j})=p_k'$ black vertices,
at each level $0\leq k\leq (n+1)$. For any  $\overline{\bf p}\in \overline{\bf V}_{n+1}$, we also have that
$
\overline{\cal J}_{\overline{\bf p}}\simeq \overline{\cal A}_{\overline{\bf p}}
$.

Our next objective is to make more precise the algebraic picture
of these objects. For any integers $p,p',q,q',r,s$,
and any disjoint maps 
$$
\overline{a}=(a\tcoprod a')\in [r]^{[(p,p')]}\quad\mbox{\rm and}\quad \overline{b}=(b \tcoprod b')\in [s]^{[(q,q')]}
$$
we let $(\overline{a}\vee \overline{b})$ be the disjoint union of maps
$$
(a\tcoprod a')\vee (b \tcoprod b')=\left(a\vee b\right)\tcoprod\left(a'\vee b'\right)\in [r+s]^{[(p+p',q+q')]}
$$
 Notice that the property of being weakly increasing is preserved by this construction.
More generally, for any $\overline{\bf p} \in\overline{\bf V}$, and $\overline{\bf q} \in\overline{\bf V}$, and any 
sequences 
$\overline{\bf a}=(\overline{a}_k)_{k\geq 0}\in {\cal \overline{A}}_{\overline{\bf p}}$, and
$\overline{\bf b}=(\overline{b}_k)_{k\geq 0}\in {\cal \overline{A}}_{\overline{\bf q}}$,    we set
\begin{equation}\label{factor}
\overline{\bf a}\vee\overline{\bf b}=(\overline{a}_k\vee \overline{b}_k)_{k\geq 0}\in {\cal
\overline{A}}_{\overline{\bf p}+ \overline{\bf q} }
\end{equation}
In the above displayed formula, we have used the conventions
$\overline{\bf p}+\overline{\bf q}=(p_k+q_k,p_k'+q_k')_{k\geq 0}$,
$\overline{a}_k=\emptyset
$, and $
\overline{b}_l=\emptyset
$,
for any $k>\overline{ht}(\overline{\bf p})$, and $l>\overline{ht}(\overline{\bf q})$.
Notice that, if $\overline{\bf a}\in{\cal I{\overline{A}}}_{\overline{\bf p}}$ and $\overline{\bf b}\in{\cal I{\overline{A}}}_{\overline{\bf q}}$, then $\overline{\bf a}\vee \overline{\bf b}\in{\cal I{\overline{A}}}_{\overline{\bf p}+\overline{\bf q}}$.

The main advantage of this algebraic structure lies in the following
non commutative composition formula in the set of colored jungles
$$
\overline{J}({\overline{\bf a}})\overline{J}({\overline{\bf b}})=\overline{J}({\overline{\bf a}}\vee {\overline{\bf b}})\in \overline{\cal J}_{\overline{\bf p}+ \overline{\bf q}}
$$
Notice that any pair of permutations $(s,s')\in\left(\Ga_p\times\Ga_{p'}\right)$ can be seen as a pair of
mappings $(s,s')\in\left([p+p']^{[p]}\times[p+p']^{[p']}\right)$. In this interpretation, the disjoint union map 
$$
\overline{s}=(s\tcoprod s')\in[p+p']^{[(p,p')]}$$ 
can be seen as a permutation map of the set $[(p,p')]=
\left([p]\coprod[p']\right)$ that preserves the sets $[p]$, and $[p']$.
The reverse assertion is also true. The set of these permutations 
is denoted by
$$
\overline{\Ga}_{(p,p')}:=\Ga_p\scoprod\Ga_{p'}
$$
We end this section with a technical lemma that permits to
extend in a natural way 
the jungle automorphism group analysis
developed
in section~\ref{treesfsec} to their colored analogue.
\begin{lem}
For any collection of integers $p,p',q,q',r,r',s,s'$, for and
sequence of 
 permutation mappings $\overline{\rho}\in \overline{\Ga}_{(p,p')}$,
$\overline{\theta}\in \overline{\Ga}_{(q,q')}$,
$\overline{\sigma}\in\overline{\Ga}_{(r,r')}$,
$\overline{\tau}\in \overline{\Ga}_{(s,s')}$, and for
any pair of disjoint mappings $\overline{a}\in[r']^{[(p,p')]}$, and
$\overline{b}\in[s']^{[(q,q')]}$ we have the composition
formula
$$
\left(\overline{\sigma}\vee\overline{\tau}\right)\circ
\left(\overline{a}\vee\overline{b}\right)\circ
\left(\overline{\rho}\vee\overline{\theta}\right)=
\left(\overline{\sigma}\circ \overline{a}\circ \overline{\rho}\right)\vee\left(\overline{\tau}\circ \overline{b}\circ \overline{\theta}\right)
$$
In addition, for any $\overline{\lambda}\in\overline{\Ga}_{(p+q),(p'+q')}$ we have that
$$
\left(\overline{a}\vee\overline{b}\right)\circ
\overline{\lambda}\in \left([r']^{[(p,p')]}\vee [s']^{[(q,q')]}\right)\Longrightarrow
\overline{\lambda}\in\left(\overline{\Ga}_{(p,p')}\vee\overline{\Ga}_{(q,q')}\right)
$$
\end{lem}
The two assertions are a direct consequence of the definition of the $\vee$ product. In the first equation above, it should be understood that $\overline{a}\in[r']^{[(p,p')]}$ is viewed as a map in $\overline{a}\in ([r]\coprod [r'])^{[(p,p')]}$ through the embedding of $[r']$ into $[r]\coprod [r']$. The same observation will be used without further notice in the forthcoming paragraph.

\subsection{Automorphism groups on colored jungles}\label{autotrop}
As for the expansion of measures in section~\ref{sectroot}, we will
need a definition of equivalence classes of colored jungles, a parametrization of the classes and a computation of their cardinals. 
Thanks to the algebraic model of colored trees designed in section~\ref{tropsect}, most of
  the constructions in section~\ref{sectroot}, 
and section~\ref{treesfsec} can be extended without further work
to colored jungles. Since all constructions
can be mimicked, we only outline the main ideas of the generalization from jungles to colored jungles.

We let $\overline{\bf p}\in \overline{V}_{n+1}$ 
be a sequence of
$(n+1)$ integer pairs
$\overline{p}_k=(p_k,p'_k)$, with $0\leq k\leq (n+1)$. 
We associate with $\overline{\bf p}$ the product
permutation group
$$
\overline{\Ga}_{\overline{\bf p}}=\overline{\Ga}_{\overline{p}_0}\times\overline{\Ga}_{\overline{p}_1}\times\ldots\times \overline{\Ga}_{\overline{p}_{n+1}}
$$
This group acts naturally on sequences of maps $\overline{\bf a}=(\overline{a}_0,\overline{a}_1,...,\overline{a}_n)\in
\overline{\cal A}_{\overline{\bf p}}$, and on jungles $\overline{J}(\overline{\bf a})$ in $\overline{\cal J}_{\overline{\bf p}}$ by permutation of the colored vertices at each level. More formally,  for any $\overline{\bf s}=(\overline{s}_0,...,\overline{s}_{n+1})\in \overline{\Ga}_{\overline{\bf p}}$ this pair of actions is given by
$$\overline{\bf s}(\overline{\bf a}):=(\overline{s}_0\overline{a}_0\overline{s}_1^{-1},...,\overline{s}_n\overline{a}_n\overline{s}_{n+1}^{-1})
\quad\mbox{\rm and}\quad
\overline{\bf s}\overline{J}(\overline{\bf a}):=\overline{J}(\overline{\bf s}(\overline{\bf a}))
$$
Two colored jungles in the same orbit under the action of $\overline{\Ga}_{\overline{\bf p}}$ have the same underlying colored forest. Conversely, if two colored jungles in $\overline{\cal J}_{\overline{\bf p}}$ have the same underlying colored forests, they differ only by a permutation of the vertices of their colored graphs that preserves the colors of the vertices, and therefore are in the same orbit under the action of  $\overline{\Ga}_{\overline{\bf p}}$.
In other terms, we have the following lemma.
\begin{lem}
Equivalence classes of colored jungles in $\overline{\cal J}_{\overline {\bf p}}$ under the action of the permutation groups $\overline{\Ga}_{\overline {\bf p}}$ are in bijection with colored forests in $\overline{\cal F}_{\overline{\bf p}}$
\end{lem}

Notice the following corollary of the lemma: if ${\overline{\bf f}}\in\overline{\cal F}_{\overline{\bf p}}$, with $\overline{\bf p}\in \overline{V}_{n+1}$, is a colored forest, then we can define unambiguously 
$$|{\overline{\bf f}}|=(|{\overline{\bf f}}|_k)_{k\geq 0}=(|\overline{a}_0|,...,|\overline{a}_{n}|,(0,0),(0,0),\ldots)$$
for any choice $(\overline{a}_0,...,\overline{a}_{n})$ of a representative of 
${\overline{\bf f}}$ in $\overline{\cal J}_{\overline{\bf p}}$.

As for usual trees and forests, we write $B$ for the map from colored trees to colored forests defined by removing the root of a colored tree, and we write also, as in section~\ref{autojung}
, $B$ for the induced map from the set of
colored forests into itself.

Let us conclude this section by enumerating the number $\# (\overline{\bf f})$ of colored jungles associated to a given colored forest $\overline{\bf f}\in \overline{\cal F}_{\overline{\bf p}}$, with  $\overline{\bf p}\in \overline{V}_{n+1}$. The process is as in section~\ref{autojung}, and the result follows ultimately from the class formula
$$
\# (\overline{\bf f})=\frac{\overline{\bf p}!}{
|Stab (\overline{\bf f})|}
$$
In the above displayed formula,
$\overline{\bf p}!$ stands for 
 the multi-index factorial  
$\overline{\bf p}!=\left(p_0!p_0'!...p_{n+1}!p_{n+1}'!\right)$, and we have written abusively $|Stab (\overline{\bf f})|$, for the cardinal of the stabilizer in $\overline{\Ga }_{\overline{\bf p}}$ of any representative $\overline{\bf f}'$ of $\overline{\bf f}$, where $\overline{\bf f}$ is viewed as an equivalence class of colored jungles. 

Let us assume that $\overline{\bf f}$ can written, as a monomial over the set of colored trees as
$$\overline{\bf f}=\prod\limits_{i=1}^k\overline{T_i}^{m_i}$$
where the $T_i$s are pairwise distincts and $m_1\geq ...\geq m_k$. We write $\overline{T_i}'$ for a set of representatives of the $\overline{T_i}$ viewed as equivalence classes of colored jungles.

As in section~\ref{orbitsec}, we shall write ${S}(B^{-1}(\overline{\bf f}))$ for the unordered $k$-uplet $(m_1,...,m_k)$, called the symmetry multiset of the colored tree $B^{-1}(\bf{\overline f})$, and extend the notation to colored forests so that
${\bf s}(\overline{\bf f})$ is the disjoint union of the ${S}(\overline{T_i})$. The proof in section~\ref{orbitsec} goes then over without changes, excepted for the replacement of forests, planar forests and planar trees by their colored analogues. For instance, we have the recursion formula
$$
|Stab (\overline{\bf f})|=\prod\limits_{i=1}^k\left(
m_i!~|Stab~B(\overline{T_i'})|^{m_i}\right)
$$
Expliciting the recursion in the Proposition, and following the wreath product analysis designed in section~\ref{wreathsec}, we deduce the following theorem.
\begin{theo}
The number of colored jungles with underlying colored forest $\overline{\bf f}\in \overline{\cal F}_{\overline{\bf p}}$, with  $\overline{\bf p}\in \overline{V}_{n+1}$, is given by
$$\# (\overline{\bf f})=\frac{\overline{\bf p}!}{\prod\limits_{i=-1}^{n}{\bf s}(B^i(\overline{\bf f}))!}
$$
where we use the usual multi-index notation to define 
${\bf s}(B^i(\overline{\bf f}))!$.
In addition, we have the wreath product representation formula 
$$Stab (\overline{\bf f}')\sim \Ga_{m_1}<Stab ~B(\overline{T}_1')>\vee\ldots\vee \Ga_{m_k}<Stab ~B(\overline{T}_k')>$$
\end{theo}

\subsection{Feynman-Kac semigroups}\label{fkpath}
In this section, the parameter $n\geq 0$ represents a fixed
time horizon. For any  sequence
of integers ${\bf q}=(q_0,\ldots,q_n)\in \NN^{n+1}$, and any 
$-1\leq m\leq n$, we set 
\begin{equation}\label{qprime}
{q}'_{m}:=q_{m+1}+{q}'_{m+1}=
\sum_{m< k\leq n}q_k
\end{equation}
Notice that ${q}'_{m}=|[q_{m+1}]\tcoprod[{q}'_{m+1}]|$. We  associate with $\bf q$ the unnormalized
Feynman-Kac measures 
$$
\gamma_{n}^{\bf q}=\gamma_0^{\otimes q_0}\otimes
\cdots\otimes \gamma_n^{\otimes q_n}
\in 
\Ma(\Eb^{\bf q}_{n})
\quad
\mbox{\rm where}
\quad
\Eb^{\bf q}_{n}=E^{q_{0}}_{0}\times\cdots\times
E^{q_{n}}_{n}
$$
is equipped with the 
tensor product sigma field. Points in $\Eb^{\bf q}_{n}$ are indexed as follows
$$
\left((x_0^1,\ldots,x_0^{q_0}),\ldots,
(x_{n}^1,\ldots,x_{n}^{q_{n}})
\right)
$$

\begin{defi}
We let $(\gamma_{p}^{\bf q})_{0\leq p\leq n}$, and
$(\Qb_{p}^{\bf q})_{0\leq p\leq n}$, be the collection of measures
and integral operators defined by
\begin{eqnarray*}
\gamma_{p}^{\bf q}&=& \left(\gamma_0^{\otimes q_0}\otimes\ldots\otimes\gamma_{p-1}^{\otimes q_{p-1}}\right)\otimes\gamma_p^{\otimes ({q}_{p}+{q}'_{p})}~\left(=\gamma_{p}^{(q_0,...,q_{p-1},q_p+{q}'_{p})}\right)\\
\Qb_{p}^{\bf q}&=&
\left(D_{1_{q_{0}}}\otimes
\cdots \otimes D_{1_{q_{p-2}}}\right)\otimes Q_p^{(q_{p-1},q'_{p-1})}
\end{eqnarray*}
with the operators $Q_p^{(q_{p-1},q'_{p-1})}$
from $E_{p-1}^{q_{p-1}+q'_{p-1}}$ into $
E_{p-1}^{q_{p-1}}\times E_p^{q'_{p-1}}$ 
defined by the tensor product formula
$$
Q_p^{(q_{p-1},q'_{p-1})}=D_{1_{q_{p-1}}}\otimes Q_p^{\otimes q'_{p-1}}
$$
\end{defi}
Notice that, for any $0\leq p\leq n$,
 $\gamma_{p}^{\bf q}$ is a positive measure
on the product space
$$
\Eb^{\bf q}_{p}=_{\mbox{\rm\tiny def.}}E^{q_{0}}_{0}\times\cdots\times E^{q_{p-1}}_{p-1}\times
E^{q_p+{q}'_{p}}_{p}
$$
and $\Qb_{p}^{{\bf q}}$ is a positive integral operator
from $\Eb^{{\bf q}}_{p-1}$ into $\Eb^{\bf q}_{p}$. \begin{lem}
For any sequence of integers ${\bf q}\in\NN^{n+1}$, and any $0\leq p\leq n$, 
we have 
$$
\gamma_{n}^{\bf q}=\gamma_{p}^{\bf q}\Qb_{p,n}^{\bf q}
$$
with the semigroup $(\Qb_{p_1,p_2}^{{\bf q}})_{1\leq p_1\leq p_2\leq n}$ defined
by
$$
\Qb_{p_1,p_2}^{\bf q}=\Qb_{p_1+1}^{{\bf q}}\left(\Qb_{p_1+1,p_2}^{{\bf q}}\right)~\left(=
\Qb_{p_1+1}^{{\bf q}}\ldots\Qb_{p_2-1}^{{\bf q}}\Qb_{p_2}^{{\bf q}}\right)
$$
and the convention $\Qb_{p_1,p_1}^{{\bf q}}={\rm Id}$, the identity
operator, for $p_1=p_2$.
\end{lem}
\proof
We start with observing that
\begin{eqnarray*}
\gamma_{n}^{{\bf q}}&=&\gamma_0^{\otimes q_0}\otimes\gamma_1^{\otimes q_1}\otimes
\cdots\otimes\gamma_{n-1}^{\otimes q_{n-1}}\otimes \left(\gamma_{n-1}^{\otimes q_n}Q^{\otimes q_n}_n\right)\\
&=&\gamma_0^{\otimes q_0}\otimes\gamma_1^{\otimes q_1}\otimes
\cdots\otimes\gamma_{n-2}^{\otimes q_{n-2}}\otimes \left(\gamma_{n-1}^{\otimes (q_{n-1}+q_{n})}(D_{1_{q_{n-1}}}\otimes
Q^{\otimes q_n}_n)\right)
\end{eqnarray*}
Therefore, we find that
$
\gamma_{n}^{{\bf q}}
=\gamma_{n-1}^{{\bf q}}
\Qb_{n}^{{\bf q}}
$.
Using a simple induction, 
the lemma follows.
\cqfd

\subsection{Unnormalized particle measures}\label{unnpath}
In this section, we derive a functional representation and a Laurent expansion of particle tensor product measures on the path space similar to the ones obtained in Section 2.2. 

In the further development of this section, the time horizon $n\geq 0$ is a fixed parameter,  and we let  ${\bf q}=(q_0,\ldots,q_n)\in \NN^{n+1}$ be a  given sequence
of $(n+1)$ integers. Before we turn to the study of the particle tensor product measures
$$
\gamma_{n}^{{\bf q},N}=(\gamma_0^N)^{\otimes q_0}\otimes(\gamma_1^N)^{\otimes q_1}\otimes
\cdots\otimes (\gamma_n^N)^{\otimes q_n}
$$
and the associated path space measures $\overline{{\QQ}}_{n,{\bf q}}^N\in{\cal M}({\bf E}_n^{\bf q})$ defined by
$$\overline{{\QQ}}_{n,{\bf q}}^N:~F\in \Ba_b({\bf E}_n^{\bf q})\longmapsto \overline{{\QQ}}_{n,{\bf q}}^N(F)=\mathbb{E}\left( \gamma_{n}^{{\bf q},N}(F)\right)\in\RR$$
let us introduce some useful definitions. 
We associate with ${\bf q}$ the pair of integer sequences $\overline{\bf q}$,  and
${\bf q'}$
defined by
$$
\overline{\bf q}:=(q_k,q'_k)_{-1\leq k\leq n}\quad\mbox{\rm and}\quad
{\bf q'}:=(q'_k)_{-1\leq k\leq  n}
$$
with the convention $q_{-1}=0$, and the integer sequence ${\bf q'}$
introduced in (\ref{qprime}).
Notice that $q'_n=0$,  $|{\bf q}|=q_0+q'_0$, $\overline{ht}(\overline{\bf q})=(n+1)$, and 
$$
|{\bf q'}|=\sum_{-1\leq k\leq n}q'_{k}=\sum_{0\leq k\leq n}(q_k+q'_{k})
=\sum_{0\leq k\leq n} (k+1)q_k
$$

\begin{defi}
We let $\Delta_{n,{\bf q}}$
be the nonnegative measure valued functional  on colored jungles 
$
\overline{J}({\overline{\bf a}})\in \overline{\cal J}_{\overline{\bf q}}$, with $
\overline{\bf a}=(\overline{a}_0,\overline{a}_1,...,\overline{a}_n)\in \overline{\cal A}_{\overline{\bf q}}$, and
defined by $$\Delta_{n,{\bf q}}: \overline{J}({\overline{\bf a}})\longmapsto\Delta_{n,{\bf q}}^{\overline{\bf a}}:=\eta_0^{\otimes |{\bf q}|}
{\bf D}^{\bf q}_{0,\overline{a}_0} {\bf Q}_1^{\bf q}{\bf D}^{\bf q}_{1,\overline{a}_1}  \ldots {\bf Q}_n^{\bf q}{\bf D}^{\bf q}_{n,\overline{a}_{n}}\in {\cal M}({\bf E}_n^{\bf q})$$
with the Markov transitions ${\bf D}^{\bf q}_{k,\overline{a}_k}$, resp.
$D_{\overline{a}_k}$, with $\overline{a}_k=({a}_k\tcoprod {a}'_k)\in [q'_{k-1}]^{[q_k]\tcoprod[q'_k]}$,  from ${\bf E}_k^{\bf q}$, resp.
$E_k^{q_k+q'_{k}}$, into itself, and defined by
$$
{\bf D}^{\bf q}_{k,\overline{a}_k}=D_{1_{q_0}}\otimes\ldots\otimes D_{1_{q_{k-1}}}\otimes D_{\overline{a}_k}
$$
and for any function $F\in \Ba_b(E_k^{q_k+q_k'})$ by
$$
D_{\overline{a}_k}(F)
(
(x_{k}^{1},\ldots,x_{k}^{q_k})
,
(
y_{k}^{1},\ldots,y_{k}^{{q}'_{k}})
)
=F(
(x_{k}^{a_k(1)},\ldots, x_{k}^{a_k(q_k)})
,
(
y_{k}^{a'_k(1)},\ldots,y_{k}^{a'_k({q}'_{k})})
)
$$
\end{defi}
\begin{defi}
We write $\Ba_b^{sym}({\bf E}_n^{\bf q})$ for the elements $F$ of $\Ba_b({\bf E}_n^{\bf q})$ that are symmetric in the first $q_0$ variables, the next $q_1$ variables,..., and the last $q_n$ variables. That is we have that
$$
\forall (\sigma_0,\ldots,\sigma_{n})\in \left(
\Ga_{q_0}\times\ldots\times\Ga_{q_n}\right)\qquad F=\left(D_{\sigma_0}\otimes\ldots\otimes D_{\sigma_n}\right)F
$$
We consider the subset $\Ba_0^{sym}({\bf E}_n^{\bf q})\subset
\Ba_b^{sym}({\bf E}_n^{\bf q})$ of all functions $F$ such that 
$$
\left(D_{1_{q_0}}\otimes\ldots\otimes D_{1_{q_{p-1}}}\otimes
\left(D_{1_{q_p-1}}\otimes\gamma_p\right)\otimes D_{1_{q_{p+1}}}\otimes\ldots\otimes D_{1_{q_n}}\right)F
=0
$$
for any $0\leq p\leq n$.  
\end{defi}

Notice that $\Ba_0^{sym}({\bf E}_n^{\bf q})$
contains the tensor product  $$
\left(
\Ba_0^{sym}(E_0^{q_0})\otimes\ldots\otimes\Ba_0^{sym}(E_n^{q_n})\right)=\left\{
F_0\otimes
\ldots\otimes F_n~:~\forall 0\leq p\leq n\quad F_p\in\Ba_0^{sym}({E}_p^{q_p}) \right\}$$ 
of the sets 
$\Ba_0^{sym}({E}_p^{q_p})$ introduced in the beginning of section~\ref{wicksetcion}, with $0\leq p\leq n$.
 For the same reasons as in equation \ref{invariance}, we also have the symmetry invariance property
$$\forall F\in \Ba_b^{sym}({\bf E}_n^{\bf q})\quad
\forall \overline{\bf s}\in \overline{\Ga}_{\overline{\bf q}}\quad\forall 
\overline{\bf f}\in\overline{J}_{\overline{\bf q}}
\qquad
\Delta_{n,{\bf q}}^{\overline{\bf f}}(F)=\Delta_{n,{\bf q}}^{\overline{\bf s} (\overline{\bf f})}(F)$$
By construction, we also have that
$$
\forall F\in \Ba_0^{sym}({\bf E}_n^{\bf q})\qquad\Delta_{n,{\bf q}}^{\overline{\bf f}}(F)=0
$$
as soon as the colored forest $\overline{\bf f}$ contains at least one trivial colored tree with a white leaf. 

Recall that a tree, or a colored tree is said to be trivial if its coalescence sequence is the null sequence of integers. Notice also that for any ${\bf l}<{\bf q'}$, the set of colored forests 
with exactly $l_k(\in [0,q'_{k-1}[)$ coalescent  edges at level $k$, is given by
$$
\overline{\cal F}_{\overline{\bf q}}[{\bf l}]:=\left\{
\overline{\bf f}\in\overline{\cal F}_{\overline{\bf q}}~
:~|\overline{\bf f}|={\bf q'}-{\bf l}\right\}
$$
Thus, for any ${\bf r}<{\bf q'}$, the set of colored forests $\overline{\cal F}_{\overline{\bf q}}({\bf r})$
with less than $r_k$ coalescent  edges at level $k$ is given by
$$
\overline{\cal F}_{\overline{\bf q}}({\bf r})
:=
\cup_{{\bf l}\leq\bf r}{\overline{\cal F}_{\overline{\bf q}}[{\bf l}]}
$$
The coalescence degree of a colored forest $\overline{\bf f}\in \overline{\cal F}_{\overline{\bf q}}[{\bf l}]$, is the sum $|{\bf q'}-{\bf l}|$ of the coalescence orders
of its vertices.  

Notice that a colored forest with coalescence degree $d$ has at most $d$ non trivial colored trees. 
In addition, a colored forest with a coalescence degree $d$, has at most
$(2d)$ leaves belonging to non trivial trees. 
Since a colored forest in $\overline{\Fa}_{\overline{\bf q}}$ has exactly $|{\bf q}|$ white leaves, if $|{\bf q}|>(2d)$, then it contains at least one trivial colored tree with a white leaf.

Next, we discuss the situation where $|{\bf q}|$ is an even integer, and we characterize the subset of forests (in $\overline{\Fa}_{\overline{\bf q}}$), with a coalescence degree $d=|{\bf q}|/2$, without any trivial colored tree with a white leaf. This charaterization follows the same lines of arguments as the ones presented on page~\pageref{discusswick}. For any integers $0\leq k\leq l\leq m\leq n$, we let
$\overline{T}_{k,l,m}$ be the unique colored tree, with a single coalescence at level $k$, and white leaves at the levels $(l+1)$,
and $(m+1)$. Notice that
colored forests with $(2r_k)$ white leaves, no black leaves, and coalescent degree $r_k$, with  $r_k$ different pairs of coalescent edges at level $k$ are necessarily of the form\label{defitreet}
$$
\overline{\bf f}_k^{({\bf t}_k)}=\prod_{k\leq l\leq m\leq n}\overline{T}_{k,l,m}^{~t_{k,l,m}}
$$
for some families of integers ${\bf t}_k=(t_{k,l,m})_{k\leq l\leq m\leq n}$ such that $|{\bf t}_k|=\sum\limits_{k,l,m}t_{k,l,m}=r_k$. Therefore, a colored forest  in $ \overline{\cal F}_{\overline{\bf q}}$, with coalescence degree $|{\bf q}|/2$,
a coalescence sequence 
${\bf r}$ such that $|{\bf r}|=|{\bf q}|/2$, without a trivial tree with a white leaf, has necessarily the following form
$$
\overline{\bf f}^{(\bf t)}=\overline{\bf f}_0^{({\bf t}_0)}U_0^{r_0}\ldots\overline{\bf f}_n^{({\bf t}_n)}U_n^{r_n}
$$
for some sequences of integers ${\bf t}=({\bf t}_k)_{0\leq k\leq n}$, such that
$|{\bf t}_k|=r_k$, for any $0\leq k\leq n$. In the above displayed formula, $U_k$ denotes the unique 
trivial tree with a single black leaf at level $k$, with $0\leq k\leq n$.
We write $\bf t!$ for $\prod\limits_{k\leq l\leq m\leq n}t_{k,l,m}!$.

We are now ready to extend the Laurent expansions, and the Wick formula presented in theorem~\ref{mainthm1-bis}, and theorem~\ref{thwick} to particle models in path spaces.

\begin{theo}\label{thwick2}
For any ${\bf q}=(q_0,\ldots,q_n)$, with $|{\bf q}|\leq N$, and any $F\in \Ba_b^{sym}({\bf E}_n^{\bf q})$, we have the Laurent expansion
\begin{equation}\label{polyfpath}
 \overline{{\QQ}}^{N}_{n,{\bf q}}(F)
 =(\gamma^{\bf q}_n+\sum_{1\leq k\leq |({\bf q'}-{\bf 1})_+|}~\frac{1}{N^{k}}~\partial^k  \overline{{\QQ}}_{n,{\bf q}})(F)
\end{equation}
with  the signed measures $\partial^k  \overline{{\QQ}}_{n,{\bf q}}$ defined by
$$\partial^k  \overline{{\QQ}}_{n,{\bf q}}(F)=
 \sum_{{\bf r}<{\bf q'},~|{\bf r}|=k} ~
  \sum_{{\overline{\bf f}}\in  \overline{\cal F}_{\overline{\bf q}}({\bf r}) }~
  \frac{s(|\overline{\bf f}|,{\bf q'}-{\bf r})~
 \#(\overline{\bf f})}{({{\bf q'}})_{|\overline{\bf f}|}}~\Delta_{n,{\bf q}}^{\overline{\bf f}}(F)
$$
In addition, for any even integer $|{\bf q}|\leq N$, and any symmetric function
$F\in \Ba_0^{ sym}({\bf E}^{\bf q}_n)$, we have 
$$
\forall k<|{\bf q}|/2\qquad\partial^k  \overline{{\QQ}}_{n,{\bf q}}
(F)=0\quad\mbox{\rm and}
\quad
\partial^{|{\bf q}|/2} \overline{{\QQ}}_{n,{\bf q}} (F)=
\sum_{{\bf r}< {\bf q'},|\bf r|=|{\bf q}|/2}~\sum_{\|{\bf t}\|={\bf r}}~\frac{{\bf q}!}{~2^{\delta({\bf t})}{\bf t}!}
\Delta^{\overline{\bf f}^{(\bf t)}}_{n,{\bf q}} F$$
with the integer sequence 
$\|{\bf t}\|=(|{\bf t}_k|)_{0\leq k\leq n}$, and the sum of the diagonal terms $\delta({\bf t})=\sum_{0\leq k\leq l\leq n} t_{k,l,l}$.
For odd integers $|{\bf q}|\leq N$, the partial measure valued
derivatives $\partial^k$ are the null measure on $\Ba_0^{sym}({\bf E}^{\bf q}_n)$, up to any order $k\leq \lfloor |{\bf q}|/2\rfloor$.
\end{theo}
Before getting into the details of the proof, we mention 
that the Wick formula stated above has a natural interpretation in  terms
of the gaussian fields $(V_k)_{0\leq k\leq n}$ introduced in section~\ref{wicksetcion}. Following the discussion given in that
section, we get that
$$
\partial^{|{\bf q}|/2} \overline{{\QQ}}_{n,{\bf q}} (F)=\EE\left(\left(V_0^{\otimes q_0}\otimes\ldots\otimes V_n^{\otimes q_n}\right)(F)\right)
$$
for any tensor product function $F$ of the following form 
$
F=F_{0,q_0}\otimes\ldots\otimes F_{n,q_n}
$, with
$F_{k,q_k}=\frac{1}{q_k!}\sum_{\sigma_k\in \Ga_{q_k}}(\varphi^{\sigma_k(1)}_k\otimes\ldots \otimes \varphi^{\sigma_k(q_k)}_k)$, 
$\varphi^i_k\in\Ba_b(E_k)$, and ${\gamma}_k(\varphi^i_k)=0$, for any $1\leq i\leq q_k$, and any $0\leq k\leq n$.
\proof
By definition of the particle model, and arguing as in the proof of (\ref{extendpath}) we find that
\begin{eqnarray*}
\EE\left(\left[
(\gamma_{n-1}^N)^{\otimes q_{n-1}}\otimes
(\gamma_{n}^N)^{\otimes q_n}\right](\varphi)
\left|~
\xi^{(N)}_{n-1}
\right.
\right)&=&\left[(\gamma_{n-1}^N)^{\otimes q_{n-1}}\otimes
\left\{
(\gamma_{n-1}^N)^{\otimes q_{n-1}}Q_n^{\otimes q_n}D_{L^N_{q_n}}\right\}\right]
(\varphi)\\
&=&\left[
(\gamma_{n-1}^N)^{\otimes (q_{n-1}+q_n)}\right]Q^{(q_{n-1},q_n)}_n\left(D_{1_{q_{n-1}}\tcoprod L^N_{q_n}}\right)(\varphi)\\
&=&\left[
(\gamma_{n-1}^N)^{\otimes (q_{n-1}+q'_{n-1})}\right]Q^{(q_{n-1},q'_{n-1})}_n\left(D_{1_{q_{n-1}}\tcoprod L^N_{q_{n}}}\right)(\varphi)
\end{eqnarray*}
for any $\varphi\in \Ba_b^{sym}(E_{n-1}^{q_{n-1}}\times E_n^{q_n})$.
This yields that
$$
\EE\left(
\gamma_{n}^{{\bf q},N}(F)~
\left|~
\xi^{(N)}_{n-1}
\right.
\right)
=
\gamma_{n-1}^{{\bf q},N}
\Qb^{{\bf q}}_n
{\bf D}^{\bf q}_{n,L^N_{q_n}}(F)=\gamma_{n-1}^{{\bf q},N}
\Qb^{{\bf q}}_n
{\bf D}^{\bf q}_{n,L^N_{q_n+q'_n}}(F)
$$
from which we readily conclude that
\begin{eqnarray*}
\EE\left(
\gamma_{n}^{{\bf q},N}(F)
\right)&=&\EE\left(
\gamma_{n-1}^{{\bf q},N}\left[
\Qb^{{\bf q}}_n
{\bf D}^{\bf q}_{n,L^N_{q_n+q'_n}}(F)\right]\right)
\end{eqnarray*}
A simple induction yields that
\begin{eqnarray*}
\EE\left(
\gamma_{n}^{{\bf q},N}(F)
\right)
&=&\EE\left(
\eta_{0}^{\otimes q_0+q'_0}
{\bf D}^{\bf q}_{0,L^N_{q_0+q'_0}}
\Qb^{{\bf q}}_1
{\bf D}^{\bf q}_{1,L^N_{q_1+q'_1}}
\ldots
\Qb^{{\bf q}}_n
{\bf D}^{\bf q}_{n,L^N_{q_n+q'_n}}(F)\right)\\
&=&\frac{1}{N^{|{\bf q'}|}}~\sum_{\overline{\bf a}\in \overline{\Aa}_{\overline{\bf q}}}~\frac{({\bf N})_{|\overline{\bf a}|}}{({\bf q'})_{|\overline{\bf a}|}}~
\Delta_{n,{\bf q}}^{\overline{\bf a}}(F)=\frac{1}{N^{|{\bf q'}|}}\sum\limits_{\overline{\bf f}\in \overline{\cal F}_{\overline{\bf q}}}\frac{({\bf N})_{\overline{\bf f}}}{({\bf q'})_{\overline{\bf f}}}~\# (\overline{\bf f})~\Delta_{n,{\bf q}}^{\overline{\bf f}}(F)
\end{eqnarray*}
from which we find the following formula
$$
\overline{{\QQ}}_{n,{\bf q}}^N(F)=\frac{1}{N^{|{\bf q'}|}}
\sum_{{\bf 1}\leq {\bf p}\leq {\bf q'}}\frac{({\bf N})_{|{\bf p}|}}{({\bf q'})_{|{\bf p}|}}~\sum\limits_{\overline{\bf f}\in \overline{\cal F}_{\overline{\bf q}}:|\overline{\bf f}|={\bf p}}\# (\overline{\bf f})~\Delta_{n,{\bf q}}^{\overline{\bf f}}(F)
$$
Using the Stirling formula (\ref{SForm}), we readily check that
\begin{eqnarray*}
\overline{{\QQ}}_{n,{\bf q}}^N&=&
\sum_{{\bf 1}\leq {\bf l}\leq {\bf q'}}
\sum_{{\bf l}\leq\bf p\leq  {\bf q'}}~s({\bf p},{\bf l}) 
\frac{1}{N^{|{\bf q'}-{\bf l}|}}
\frac{1}{({\bf q'})_{|{\bf p}|}}~\sum\limits_{\overline{\bf f}\in \overline{\cal F}_{\overline{\bf q}}:|\overline{\bf f}|={\bf p}}\# (\overline{\bf f})~\Delta_{n,{\bf q}}^{\overline{\bf f}}
\end{eqnarray*}
From previous computations, we conclude that
\begin{eqnarray*}
\overline{{\QQ}}_{n,{\bf q}}^N
&=&\sum_{{\bf r}< {\bf q'}}
\frac{1}{N^{|{\bf r}|}}
~\sum\limits_{\overline{\bf f}\in \overline{\cal F}_{\overline{\bf q}}({\bf r})}~ s(|\overline{\bf f}|,{\bf q'}-{\bf r}) 
\frac{1}{({\bf q'})_{|\overline{\bf f}|}}~
~\#(\overline{\bf f})~\Delta_{n,{\bf q}}^{\overline{\bf f}}
\end{eqnarray*}
Finally, we notice that $\overline{\cal F}_{\overline{\bf q}}({\bf 0})$
reduces to the single class of all sequences of bijections in $\overline{\Aa}_{\overline{\bf q}}$.
The end the proof of the first assertion is now clear. To end the proof of the theorem, we notice that (with $\bf t$ as in the expansion of $\partial^{{|{\bf q}|}/2}\overline{\QQ}_{n,{\bf q}}(F)$)
$$
\#(\overline{\bf f}^{(\bf t)})=\frac{\overline{\bf q}!}{\prod_{0\leq k\leq n}\left[~r_k!~\left(\prod_{k\leq l\leq m\leq n}t_{k,l,m}!\right)~\left(\prod_{k\leq l\leq n}2^{t_{k,l,l}}\right)\right]}
$$
Since we have ${\bf s}(|\overline{\bf f}^{(\bf t)}|,{\bf q'}-{\bf r})=1$, and
$({\bf q'})_{|\overline{\bf f}^{(\bf t)}|}=\prod_{-1\leq k<n}(q'_k)_{q'_k-r_{k+1}}$, for any ${\bf r}=(r_k)_{0\leq k\leq n}$, we conclude that
$$
\frac{{\bf s}(|\overline{\bf f}^{(\bf t)}|,{\bf q'}-{\bf r})}{({\bf q'})_{|\overline{\bf f}^{(\bf t)}|}}~\#(\overline{\bf f}^{(\bf t)})=\frac{{\bf q}!}{~2^{\delta({\bf t})}~{\bf t}!}
$$
with ${\bf t}!=\prod_{0\leq k\leq l\leq m\leq n}t_{k,l,m}!$. The end of the proof of the theorem is now straightforward.\cqfd
\subsection{Propagations of chaos type expansions}\label{finalsec}
This section is essentially concerned with applications of the
differential forest expansion machinery developed 
earlier to propagations of chaos properties of interacting 
particle models. In order to state, and prove the main 
results of this section, we need to introduce some notation.
We shall work throughout with a fixed time horizon $n\geq 0$, and a
constant particle block size $q\leq N$. We let $\NN^{n+1}_q$ be the set of integer sequences ${\bf p}=(p_k)_{0\leq k\leq n}\in\NN^{n+1}$
such that  $|{\bf p}|=q$.
We associate to a given ${\bf p}=(p_k)_{0\leq k\leq n}\in \NN^{n+1}$ the pair of integer sequences ${\bf p'}$ and ${\bf p}+q$ defined by
$$
\forall 0\leq k\leq n\qquad p'_k:=\sum_{k<l\leq n}p_l
\quad\mbox{\rm and}\quad  {\bf p}+q:=(p_0,\ldots,p_{n-1},p_{n}+q)
$$
We also denote by $(\overline{G}_k)_{0\leq k\leq n}$, and  
$\overline{G}_{n}^{\otimes {\bf p}}$ the collection of functions 
defined by
$$
 \overline{G}_k:=\frac{1}{\gamma_{k}(G_k)}\times \left({\eta_{k}(G_k)}-G_k\right)
\quad\mbox{\rm and}\quad
\overline{G}_{n}^{\otimes {\bf p}}:=\overline{G}_{0}^{\otimes p_0}\otimes\ldots\otimes \overline{G}_{n}^{\otimes p_n}
$$
Finally, we consider the integral operators $S^{({\bf p}+q)}_{n}$
from ${\bf E}^{{\bf p}+q}_{n}$ into ${E}^{q}_{n+1}$, defined for any
function $F_{n+1}\in\Ba_b(E^{q}_{n+1})$ by the following formula
\begin{eqnarray*}
S^{({\bf p}+q)}_{n}(F_{n+1})&:=&
\overline{G}_{0}^{\otimes p_0}\otimes\ldots\otimes
\overline{G}_{n-1}^{\otimes p_{n-1}}\otimes \left(\overline{G}_{n}^{\otimes p_n}\otimes Q_{n+1}^{\otimes q}\overline{F}_{n+1}\right)_{\rm\tiny sym}
\end{eqnarray*}
with
$$
\overline{F}_n:=\frac{1}{\gamma_{n}(1)^q}\left(F_n-\eta^{\otimes q}_{n}(F_n)\right)$$
\begin{prop}\label{lpb}
For any $q\leq N$, and any $n\geq 0$, we have the polynomial
decompositions
$$
E_{q,n}^N:=\EE\left(\left(1- {\gamma_{n}^N(G_n)}/{\gamma_{n}(G_n)}\right)^q\right)=\sum_{q/2\leq k\leq (n+1)(q-1)}\frac{1}{N^k}~\partial^k
E_{q,n}
$$
The derivatives of order $q/2\leq k\leq (n+1)(q-1)$ are given by the 
following formula
$$
  \partial^k E_{q,n}= \sum_{{\bf p}\in \NN^{n+1}_q}\sum_{k\leq |({\bf p'}-1)_+|}~
  \frac{q!}{{\bf p}!}~\partial^k\overline{{\QQ}}_{n,{\bf p}}
  \left(\overline{G}_{n}^{\otimes {\bf p}}\right)
$$

\end{prop}
\proof
We first check the following decomposition
\begin{equation}\label{normd}
1- {\gamma_{n}^N(G_n)}/{\gamma_{n}(G_n)}=\sum_{0\leq p\leq n}~\gamma^N_p( \overline{G}_p)
\end{equation}
To prove this formula, we simple notice that
 \begin{eqnarray*}
 \gamma^N_{n}(G_n)-\gamma_{n}(G_n)&=&\prod_{0\leq p\leq n}\eta^N_p(G_p)-
 \prod_{0\leq p\leq n}\eta_p(G_p)\\
 &=&\gamma^N_n( G_n-\eta_n(G_n))+\left[   \gamma^N_{n-1}(G_{n-1})-\gamma_{n-1}(G_{n-1})\right]\times \eta_n(G_n)\\
 &=&\ldots\\
 &=&\sum_{0\leq p\leq n}~\gamma^N_p(G_p-\eta_p(G_p))\times\left[\prod_{p+1\leq k\leq n}
\eta_k(G_k)\right]
 \end{eqnarray*}
 with the convention $\prod_{\emptyset}=1$. 
Finally, combining the multinomial decomposition
\begin{eqnarray}
 \left( 1-{\gamma_{n}^N(G_n)}/{\gamma_{n}(G_n)}\right)^q
 &=& \sum_{|{\bf p}|=q}~\frac{q!}{{\bf p}!}~
 \gamma^{{\bf p},N}_n\left(\overline{G}_{n}^{{\bf p}}\right)
\label{multidecomp}
 \end{eqnarray}
with the Wick expansion stated in Thm.\ref{thwick2}
we conclude  that
$$
E_{q,n}^N=\sum_{{\bf p}\in \NN^{n+1}_q}~\frac{q!}{{\bf p}!}~
 \overline{{\QQ}}^{N}_{n,{\bf p}}\left(\overline{G}_{n}^{{\bf p}}\right)=
 \sum_{{\bf p}\in \NN^{n+1}_q} \sum_{q/2\leq k\leq |({\bf p'}-1)_+|}~\frac{1}{N^k}~\frac{q!}{{\bf p}!}~\partial^k\overline{{\QQ}}_{n,{\bf p}}\left(\overline{G}_{n}^{{\bf p}}\right)
$$
This completes the proof of the proposition. 
\cqfd
\begin{theo}\label{theoterminal}
For any $n\geq 0$, the sequence of probability measures
 $(\PP^N_{n+1,q})_{N\geq q}$ is differentiable up to any order
with $\partial^0\PP_{n+1,q}=\eta_{n+1}^{\otimes q}$, and the partial derivatives given by
the following formula
$$
 \partial^k\PP_{n+1,q}
=
\sum_{{\bf p}\in \cup_{0\leq l<2k}\NN^{n+1}_l}~\frac{((q-1)+|{\bf p}|)!}{(q-1)!~{\bf p}!}~\partial^k\overline{{\QQ}}_{n,{\bf p}+q}
S^{{\bf p}+q}_{n}
$$
At any order $N\geq q$, we have the exact formula
$$
\PP^N_{n+1,q}=\eta^{\otimes q}_{n+1}+\sum_{1\leq k<\lfloor (N-q)/2\rfloor}~\frac{1}{N^{k}}~\partial^k\PP_{n+1,q}+\frac{1}{N^{\lfloor (N-q)/2\rfloor}}~\partial^{\lfloor (N-q)/2\rfloor}\PP_{n+1,q}^N
$$
\end{theo}

\proof
We let $F_{n+1}$ be a bounded measurable symmetric function 
on ${\bf E}^q_{n+1}$,  such that $\eta_{n+1}^{\otimes q}(F_{n+1})=0$.
By definition of the particle model, we have that
\begin{eqnarray*}
\EE\left(
\left(\eta^N_{n+1}\right)^{\odot q}(F_{n+1})~|~\xi_{n}^{(N)}\right)&=&\gamma^N_{n}(G_n)^{-q}\times
\left(\gamma^{N}_{n}\right)^{\otimes q}
Q_{n+1}^{\otimes q}\left(F_{n+1}\right)\\
&=&\left(1-u^N_{n}\right)^{-q}\times
\left(\gamma^{N}_{n}\right)^{\otimes q}~
Q_{n+1}^{\otimes q}\left(\overline{F}_{n+1}\right)
\end{eqnarray*}
 with the sequence of random variables 
 $
u^N_n:=\left(1-{\gamma^N_{n}(G_n)}/{\gamma_{n}(G_n)}\right)
$.
On the other hand,  
for any $u\in\RR-\{1\}$, and $m\geq 1$, we have the more of less
well known decomposition
 $$
 \frac{1}{(1-u)^{q+1}}=\sum_{0\leq k\leq m}~(q+k)_k~\frac{u^k}{k!}+u^m~\sum_{1\leq k\leq q+1}~\left(\begin{array}{c}(q+1)+m\\ k+m\end{array}\right)~
 \left(\frac{u}{1-u}\right)^k
 $$
with $(q+k)_k={(q+k)!}/{q!}$. The proof of this formula is 
essentially based on the fact that for any $n\geq 0$, and $u\not=1$, we have
$$
f(u)=\frac{1}{1-u}=\sum_{0\leq k\leq n}~u^k~+\frac{u^{n+1}}{1-u}
\mbox{and}
\frac{\partial^n f}{\partial u^n}=n!~f^{n+1}
$$
This implies that for any $m> n$, we have
$$
\frac{\partial^n f}{\partial u^n}=
\frac{n!}{(1-u)^{n+1}}=\sum_{n\leq k\leq m}~(k)_{n}~u^{k-n}+
\frac{\partial^n }{\partial u^n}\left(\frac{u^{m+1}}{1-u}
\right)
$$
Recalling the Leibniz binomial derivation formula
$$
\frac{\partial^n }{\partial u^n}(fg)=\sum_{0\leq k\leq n}~
\left(\begin{array}{c}
 n\\
 k\end{array}\right)~\frac{\partial^{n-k} }{\partial u^{n-k}}(f)~\frac{\partial^k }{\partial u^k}(g)
$$
we find that
\begin{eqnarray*}
\frac{1}{n!}~\frac{\partial^n }{\partial u^n}\left(\frac{u^{m+1}}{1-u}
\right)
&=&\sum_{0\leq k\leq n}~
\left(\begin{array}{c}
 n\\
 k\end{array}\right)~(m+1)_k~u^{(m+1)-k}~
 \frac{(n-k)!}{n!}
~\frac{1}{(1-u)^{(n-k)+1}}
\\
&=&u^{m-n}~\sum_{0\leq k\leq n}~
\left(\begin{array}{c}
 m+1\\
 k\end{array}\right)~\frac{u^{(n-k)+1}}{(1-u)^{(n-k)+1}}
\end{eqnarray*}
and
$$
\frac{1}{n!}~\frac{\partial^n }{\partial u^n}\left(\frac{u^{m+1}}{1-u}
\right)
=u^{m-n}~\sum_{0\leq k\leq n}~
\left(\begin{array}{c}
 m+1\\
 n-k\end{array}\right)~\frac{u^{k+1}}{(1-u)^{k+1}}
$$
Thus,
$$
\frac{1}{(1-u)^{n+1}}=\sum_{0\leq k\leq (m-n)}~
 \frac{(k+n)!}{n!}
~\frac{u^{k}}{k!}+
u^{m-n}~\sum_{0\leq k\leq n}~
\left(\begin{array}{c}
 m+1\\
 n-k\end{array}\right)~\frac{u^{k+1}}{(1-u)^{k+1}}
$$
We prove the desired decomposition
with replacing $m$ by $(m+q)$, and $n$ by $q$.
Using the above decomposition we find that
$$
\begin{array}{l}
\EE\left(
\left(\eta^N_{n+1}\right)^{\odot q}(F_{n+1})~|~\xi_{n}^{(N)}\right)\\
\\
=\displaystyle\sum_{0\leq k\leq m}~((q-1)+k)_k~\frac{1}{k!}
\left(1-\frac{\gamma^N_{n}(G_n)}{\gamma_{n}(G_n)}\right)^k~\left(\gamma_{n}^N\right)^{\otimes q}(Q_{n+1}^{\otimes q}\overline{F}_{n+1})+R^{q,N}_{m,n}({F}_{n+1})
\end{array}
$$
with the remainder term
$$
\begin{array}{l}
R^{q,N}_{m,n}(F_{n+1})\\
\\=\left(1-\frac{\gamma^N_{n}(G_n)}{\gamma_{n}(G_n)}\right)^{m+1}
\left(\gamma_{n}^N\right)^{\otimes q}(Q_{n+1}^{\otimes q}\overline{F}_{n+1})
~\sum_{1\leq k\leq q}~\left(\begin{array}{c}q+m\\ k+m\end{array}\right)~
 {\left(1-\frac{\gamma^N_{n}(G_n)}{\gamma_{n}(G_n)}\right)^{k-1}}/{
 \left(\frac{\gamma^N_{n}(G_n)}{\gamma_{n}(G_n)}\right)^{k}}
\end{array}
$$
We also have from proposition~\ref{lpb} 
$$
\sup_{N\geq m+1}{\sqrt{N}
\EE\left(\left|\gamma^N_{n}(G_n)-\gamma_{n}(G_n)\right|^{m+1}\right)^{\frac{1}{m+1}}}<\infty
$$
Finally, using the multinomial decomposition (\ref{multidecomp}),
we conclude  that
$$
\begin{array}{l}
\displaystyle\sum_{0\leq k\leq m}~((q-1)+k)_k~\frac{1}{k!}
\left(1-\frac{\gamma^N_{n}(G_n)}{\gamma_{n}(G_n)}\right)^k~\left(\gamma_{n}^N\right)^{\otimes q}(Q_{n+1}^{\otimes q}\overline{F}_{n+1})\\
\\
=\displaystyle\sum_{0\leq k\leq m}~((q-1)+k)_k~
\displaystyle\sum_{{\bf p}\in\NN^{n+1}_k}~\frac{1}{{\bf p}!}~\gamma_{n}^{({\bf p}+q),N}\left(
S^{({\bf p}+q)}_{n}(F_{n+1})\right)
\end{array}
$$
This yields, for any
$m+q\leq N$, the functional expansion
$$
\PP_{n+1,q}^N
=
\displaystyle\sum_{|{\bf p}|\leq m}~\frac{((q-1)+|{\bf p}|)!}{(q-1)!~{\bf p}!}~\overline{{\QQ}}_{n,{\bf p}+q}^{N}
S^{{\bf p}+q}_{n}+\overline{R}^{q,N}_{m,n}
$$
with a remainder measure $\overline{R}^{q,N}_{m,n}$, such that
$
\sup_{N\geq 1}{N^{(m+1)/2}\left\|\overline{R}^{q,N}_{m,n}\right\|_{\rm\tiny tv}}<\infty
$. The last assertion follows from the regularity hypothesis \ref{potfun} on the potential functions.
This implies that for any $k<(m+1)/2$ we have
$$
\partial^k\PP_{n+1,q}
=
\displaystyle\sum_{|{\bf p}|\leq m}~\frac{((q-1)+|{\bf p}|)!}{(q-1)!~{\bf p}!}~\partial^k\overline{{\QQ}}_{n,{\bf p}+q}
S^{{\bf p}+q}_{n}
$$
Notice that $k$-th order derivative
measure
$\partial^k\overline{{\QQ}}_{n,{\bf p}+q}$ only involves colored 
forests 
with less than $k$ coalescent branches, from the original
root, up to the final level.
If $|{\bf p}|\geq (2k)$, then these colored forests contain at least one elementary tree with a white leaf.
 By definition of the operator $S^{{\bf p}+q}_{n}$ we find that
$$
\partial^k\overline{{\QQ}}_{n,{\bf p}+q}
S^{{\bf p}+q}_{n}(F_{n+1})=0
$$
This yields that
$$
\PP^N_{n+1,q}=\eta^{\otimes q}_{n+1}+\sum_{1\leq k< ((N-q)+1)/2}~\frac{1}{N^{k}}~\partial^k\PP_{n+1,q}+\RR^N_{n+1,q}
$$
with a remainder measure $\RR^N_{n+1,q}$ such that $
\sup_{N\geq 1}{N^{((N-q)+1)/2 }\left\|\RR^N_{n+1,q}
\right\|_{\rm\tiny tv}}<\infty
$. We end of the proof of the second assertion of the theorem using the fact that
$$
\lfloor (N-q)/2\rfloor\leq (N-q)/2\leq ((N-q)+1)/2
$$
The Theorem follows.
\cqfd

We end this section, with a series of some direct consequences of the above theorem: 
\begin{itemize}
\item The differential forest expansions presented in  theorem~\ref{theoterminal} allow us to deduce precise strong propagation of chaos estimates. For instance, 
for any $N\geq (q+7)$ (so that $\left(\lfloor (N-q)/2\rfloor-1\right)\geq 2$) we have
$$
\sup_{N\geq q+7}{N^2~\left\|\PP_{n+1,q}^N-\eta^{\otimes q}_{n+1}-\frac{1}{N} ~\partial^1\PP_{n+1,q}\right\|_{\rm\tiny tv}}<\infty
$$
with a first order partial derivative given by the formula
$$
 \partial^1\PP_{n+1,q}(F_{n+1})
=\partial^1 {\QQ}_{n,q}
Q^{\otimes q}_{n+1}(\overline{F}_{n+1})+q~\sum_{0\leq m\leq n}\partial^1\overline{{\QQ}}_{n,{\bf q}_m}
S^{{\bf q}_m}_{n}(F_{n+1})
$$
with the sequence of integers ${\bf q}_m=(1_m(k)+q1_n(k))_{0\leq k\leq n}$. The first term in the right hand side in the above displayed
formula has been treated in corollary~\ref{3order}, on page~\pageref{3order}, and we have that
$$
\begin{array}{l}
\partial^1 {\QQ}_{n,q}
Q^{\otimes q}_{n+1}(\overline{F}_{n+1})\\
\\
=\displaystyle\frac{q(q-1)}{2}~\displaystyle\sum_{0\leq k\leq n}\Delta^{{\bf f}_{1,k}}_{n,q}Q^{\otimes q}_{n+1}(\overline{F}_{n+1})\\
=\displaystyle\frac{q(q-1)}{2}~\displaystyle\sum_{0\leq k\leq n}\gamma_k(1)\displaystyle\int_{E_k^{q-1}}\gamma_k^{(q-1)}(d(x^2,\ldots,x^{q}))~Q^{\otimes q}_{k,n+1}(\overline{F}_{n+1})(x^2,x^2,x^3,\ldots,x^q)
\end{array}
$$
Each of the terms $\partial^1\overline{{\QQ}}_{n,{\bf q}_m}
S^{{\bf q}_m}_{n}(F_{n+1})$ only involves the colored forests 
$$
\forall 0\leq k\leq m\qquad \overline{\bf f}_{k,m}:=U_kT_{k,m,n}U_{n+1}^{q-1}
$$ 
associated with the  trees $T_{k,m,n}$, and $U_l$, 
introduced on page \pageref{defitreet}. After some elementary
manipulations, we find that
$$
\begin{array}{l}
\partial^1\overline{{\QQ}}_{n,{\bf q}_m}S^{{\bf q}_m}_{n}(F_{n+1})
\\
\\
=q~\displaystyle\sum_{0\leq k\leq  m}\Delta^{\overline{\bf f}_{k,m}}
_{n,{\bf q}_m}S^{{\bf q}_m}_{n}(F_{n+1})\\
=q~\displaystyle\sum_{0\leq k\leq  m}\gamma_k(1)~\displaystyle\int_{E^q_k}~\gamma_k^{\otimes q}(d(x^1,\ldots,x^q))~Q_{k,m}(\overline{G}_m)(x^1)~Q^{\otimes q}_{k,n+1}(\overline{F}_{n+1})(x^1,\ldots,x^q)
\end{array}
$$
\item A Wick formula derives from theorem~\ref{thwick2}.
More precisely, for any $F\in\Ba_0^{\rm\tiny sym}(E^{q}_{n+1})$, 
the partial derivatives $\partial^k\PP_{n+1,q}(F)$
are null up to any order $k\leq \lfloor q/2\rfloor$, and for  any even integer $q$, we have
$$
\partial^{q/2}\PP_{n+1,q}(F)=\gamma_{n+1}(1)^{-q}~\partial^{q/2}{\QQ}_{n,q}
Q^{\otimes q}_{n+1}(F)
$$
as soon as $N\geq 2(q+2)$ (so that $\left(\lfloor (N-q)/2\rfloor-1\right)\geq q/2$). To prove this claim, we notice that for any $k<q/2$, and any $F\in\Ba_0^{\rm\tiny sym}(E^{q}_{n+1})$, we have
$$
\forall {\bf p}\in \cup_{0\leq l<2k}\NN^{n+1}_l\qquad
k<q/2\leq  (|{\bf p}|+q)/2\quad\mbox{\rm and therefore}\quad
\partial^k\overline{{\QQ}}_{n,{\bf p}+q}
S^{{\bf p}+q}_{n}(F)=0
$$
This yields that $\partial^k\PP_{n+1,q}(F)=0$, for any even integer $q$,
and any $k<q/2$. In the case $k=q/2$, we have
$
k=q/2=(|{\bf p}|+q)/2
$ if, and only if, ${\bf p}$ coincide with the null sequence ${\bf 0}$. 
\item We let $\Vert \mu\Vert_{\Ba_{0,1}^{\rm\tiny sym}}=\sup_{F\in\Ba_{0,1}^{\rm\tiny sym}}\left|\mu(F)\right|$ be the Zolotarev seminorm on $\Ma(E^q_{n+1})$ associated with
the collection of functions 
$$
\Ba_{0,1}^{\rm\tiny sym}:=\left\{F\in \Ba_{0}^{\rm\tiny sym}(E^{q}_{n+1})~:~\Vert F\Vert\leq 1\right\}
$$
For any even integer $q$ such that $(q+2)\leq N/2$, we have
$$
\sup_{N\geq q+7}{N^{1+q/2}\left\|\PP_{n+1,q}^N-\frac{1}{N^{q/2}} ~\frac{1}{
\gamma_{n+1}(1)^{q}}~\partial^{q/2}{\QQ}_{n,q}
Q^{\otimes q}_{n+1}
\right\|_{\Ba_{0,1}^{\rm\tiny sym}}}<\infty
$$
\item Combining 
the  Wick formula stated above with the Borel-Cantelli lemma, we 
obtain for all $q\geq 4$ the almost sure convergence result
$$
\lim_{N\rightarrow\infty}(\eta_n^N)^{\odot q}(F)=
0 ~~\mbox{p.s.}
$$
for any bounded symmetric function $F\in \Ba_0^{\rm\tiny sym}(E^{q}_{n})$. This result is an extension of the law of large numbers
for $U$-statistics obtained by W. Hoeffding~\cite{hoeffding} for independent and identically distributed
random variables to interacting particle models. 

\item We mention that the same lines of arguments used in the proof of  
theorem~\ref{theoterminal} show that the sequence of probability measures
$$
\tilde\PP^N_{n+1,q}~:~F\in\Ba_b(E_{n+1}^q)\mapsto
\tilde\PP^N_{n+1,q}(F):=\EE\left(\left(\eta^N_{n+1}\right)^{\otimes q}(F)\right)
$$
  is differentiable up to order $\lfloor (N-q)/2\rfloor$, 
with $\partial^0\tilde\PP_{n+1,q}=\eta_{n+1}^{\otimes q}$, and the partial derivatives given for any $1\leq k<\lfloor (N-q)/2\rfloor$ by
the following formula
$$
 \partial^k\tilde\PP_{n+1,q}(F)
=
\sum_{{\bf p}\in \cup_{0\leq l<2k}\NN^{n+1}_l}~\frac{((q-1)+|{\bf p}|)!}{(q-1)!~{\bf p}!}~\partial^k\overline{{\QQ}}_{n+1,({\bf p},q)}
\left(\overline{G}_n^{\bf p}\otimes \overline{F}\right)
$$
In the same way, for any $F\in\Ba_0^{\rm\tiny sym}(E^{q}_{n+1})$,
and $N\geq 2(q+2)$,   
the partial derivatives $\partial^k\QQ_{n+1,q}(F)$
are null up to any order $k<q/2$, and for any even integer $q$,  we have that
$$
\partial^{q/2}\tilde\PP_{n+1,q}(F)=\gamma_{n+1}(1)^{-q}~\partial^{q/2}{\QQ}_{n+1,q}(F)
$$
\item Finally, the Wick formula stated above allows to deduce sharp
$\LL_q$-mean error bound. To see this claim, we simply observe that
$$
F=(f-\eta_{n+1}(f))^{\otimes q}\quad\mbox{\rm with}\quad f\in\Ba_b(E_{n+1})\Longrightarrow
\tilde\PP^N_{n+1,q}(F)=\EE((\eta_{n+1}^N)^{\otimes q}(f-\eta_{n+1}(f))^{\otimes q})
$$
and for any even integer $q$ such that $(q+2)\leq N/2$, we have
the inequality
$$
\sup_{N\geq q+7}{N^{1+q/2}\left\|\tilde\PP_{n+1,q}^N-\frac{1}{N^{q/2}} ~\frac{1}{
\gamma_{n+1}(1)^{q}}~\partial^{q/2}{\QQ}_{n+1,q}
\right\|_{\Ba_{0,1}^{\rm\tiny sym}}}<\infty
$$

\end{itemize}

 \end{document}